\documentclass[11pt,twoside,leqno]{aomamlt2e} 
  \pageno{643}
\received{December 4, 2002}
  \renewcommand{\thesection}{\arabic{section}}
\makeatletter
\renewcommand{\@seccntformat}[1]{\csname
the#1\endcsname.\hspace{0.5em}\setcounter{Subsec}{0}\setcounter{Subsubsec}{0}}\makeatother

\newcommand\sdemo[1]{\demo{\scshape #1}}

   \DeclareMathOperator{\lca}{lca}
  \DeclareMathOperator{\diam}{diam}
\newcommand{\A}{\mathcal{A}}
\newcommand{\C}{\mathcal{C}}
\newcommand{\E}{\mathbb{E}}

\newcommand{\M}{\mathcal{M}}
\newcommand{\N}{\mathcal{N}}
\newcommand{\PP}{\mathcal{P}}
\newcommand{\R}{\mathbb{R}}

\newcommand{\F}{\mathcal{F}}
\newcommand{\HH}{\mathcal{H}}

\newcommand{\e}{\varepsilon}

\newcommand{\hst}{{\rm HST}}
\newcommand{\um}{{\rm UM}}
\newcommand{\eq}{{\rm EQ}}
\newcommand{\comp}{{\rm comp}}
\newcommand{\XX}{\psi} 
\newcommand{\xx}{{\psi}} 
\newcommand{\lvs}{{\rm lvs}}

\newcommand{\classeqv}[2]{{{#1}^{\stackrel{#2}{\hookleftarrow}}}}

\newcommand{\mommit}[1]{}
 
\theoremstyle{plain}
  \newtheorem{lemma}{Lemma}[section]
  \newtheorem{theorem}[lemma]{Theorem}
  \newtheorem{corollary}[lemma]{Corollary}
  \newtheorem{claim}[lemma]{Claim}
  \newtheorem{proposition}[lemma]{Proposition}

  \theorembodyfont{\upshape}
  \newtheorem{definition}[lemma]{{\it Definition}}

  \newtheorem{remark}[lemma]{{\it Remark}}
  \newtheorem{notation}[lemma]{{\it Notation}}

\begin{document}
\currannalsline{162}{2005} 

 \title{On metric Ramsey-type phenomena}

 \acknowledgements{}
\twoauthors{Yair Bartal, Nathan Linial, Manor Mendel,}{Assaf Naor}

 \institution{Institute of Computer Science, Hebrew University,  Jerusalem,
Israel
\\
\email{yair@cs.huji.ac.il}
\\
\vglue-9pt
Institute of Computer Science, Hebrew University, Jerusalem, Israel
\\
\email{nati@cs.huji.ac.il}\\
\vglue-9pt
Institute of Computer Science, Hebrew University, Jerusalem, Israel
\\
\current{Computer Science Division, The Open University
of Israel, Ra'anana, Israel}
\\
\email{mendelma@gmail.com}
\\
\vglue-9pt
Theory  Group, Microsoft Research, Redmond, WA
\\
\email{anaor@microsoft.com}}

 \shorttitle{On metric Ramsey-type phenomena}

  \centerline{\bf Abstract}
\vglue12pt

The main question studied in this article may be viewed as a
nonlinear analogue of Dvoretzky's theorem in Banach space theory
or as part of Ramsey theory in combinatorics. Given a finite
metric space on $n$ points, we seek its subspace of largest
cardinality which can be embedded with a given distortion in
Hilbert space. We provide nearly tight upper and lower bounds on
the cardinality of this subspace in terms of $n$ and the desired
distortion. Our main theorem states that for any $\epsilon>0$,
every $n$ point metric space contains a subset of size at least
$n^{1-\epsilon}$ which is embeddable in Hilbert space with
$O\left(\frac{\log(1/\epsilon)}{\epsilon}\right)$ distortion. The
bound on the distortion is tight up to the $\log(1/\epsilon)$
factor. We further include a comprehensive study of various other
aspects of this problem.

\vglue12pt

\centerline{\bf Contents}
 \begin{small}
\def\sni#1{\smallbreak\noindent{#1}.}
\def\ssni#1{\vglue-1pt\noindent\hskip13pt  {#1}.}
\vglue-4pt

\sni{1} Introduction
\ssni{1.1} Results for arbitrary metric spaces
\ssni{1.2} Results for special classes of metric spaces
\sni{2} Metric composition
\ssni{2.1} The basic definitions
\ssni{2.2} Generic upper bounds via metric composition
\sni{3} Metric Ramsey-type theorems 
\ssni{3.1} Ultrametrics and hierarchically well-separated trees
\ssni{3.2} An overview of the proof of Theorem~1.3
\ssni{3.3} The weighted metric Ramsey problem and its relation to metric composition
\ssni{3.4} Exploiting metrics with bounded aspect ratio
\ssni{3.5} Passing from an ultrametric to a $k$-HST
\ssni{3.6} Passing from a $k$-HST to metric composition
\ssni{3.7} Distortions arbitrarily close to $2$
\sni{4} Dimensionality based upper bounds
\sni{5} Expanders and Poincar\'e inequalities
\sni{6} Markov type, girth and hypercubes
\ssni{6.1} Graphs with large girth
\ssni{6.2} The discrete cube

\end{small}

\section{Introduction}\label{introduction}

The philosophy of modern Ramsey theory states that large systems
necessarily contain large, highly structured sub-systems. The
classical Ramsey coloring theorem \cite{ramsey}, \cite{grs} is a prime
example of this principle: Here ``large'' refers to the
cardinality of a set, and ``highly structured'' means being
monochromatic.

Another classical theorem, which can be viewed as a Ramsey-type
phenomenon, is Dvoretzky's theorem on almost spherical sections of
convex bodies. This theorem, a cornerstone of modern Banach space
theory and convex geometry, states that for all $\epsilon>0$,
every $n$-dimensional normed space $X$ contains a $k$-dimensional
subspace $Y$ with $d(Y,\ell_2^k)\le 1+\epsilon$, where $k\ge
c(\epsilon)\log n$. Here $d(\cdot,\cdot)$ is the Banach-Mazur
distance, which is defined for two isomorphic normed spaces
$Z_1,Z_2$ as:
$$
d(Z_1,Z_2)=\inf\{\|T\|\cdot\|T^{-1}\|; T\in{\rm GL}(Z_1,Z_2)\}.
$$
Dvoretzky's theorem is indeed a Ramsey-type theorem, in which
``large'' is interpreted as high-dimensional, and ``highly
structured'' means close to Euclidean space in the Banach-Mazur
distance.

Dvoretzky's theorem was proved in \cite{dvo}, and the estimate
$k\ge c(\e)\log n$, which is optimal as a function of $n$, is due
to Milman~\cite{milman}. The dimension of almost spherical
sections of convex bodies has been studied in depth by Figiel,
Lindenstrauss and Milman in~\cite{flm}, where it was shown that
under some additional geometric assumptions, the logarithmic lower
bound for $\dim(Y)$ in Dvoretzky's theorem can be improved
significantly. We refer to the books \cite{milschechtbook}, \cite{pisier}
for good expositions of Dvoretzky's theorem, and to
\cite{milschecht1}, \cite{milschecht2} for an ``isomorphic'' version of
Dvoretzky's theorem.

The purpose of this paper is to study nonlinear versions of
Dvoretzky's theorem, or viewed from the combinatorial perspective,
metric Ramsey-type problems.
    In spite of the similarity of these problems, the
results in the metric setting differ markedly from those for the
linear setting.

Finite metric spaces and their embeddings in other metric spaces
have been intensively investigated in recent years. See for
example the surveys \cite{indyk}, \cite{linial}, and the book
\cite{matbook} for an exposition of some of the results.

Let $f:X\to Y$ be an embedding of the metric spaces $(X,d_X)$ into
$(Y,d_Y)$. We define the {\em distortion} of $f$ by
$$
\mathrm{dist}(f)=\sup_{\substack{x,y\in X\\x\neq
y}}\frac{d_Y(f(x),f(y))}{d_X(x,y)}\cdot \sup_{\substack{x,y\in X\\x\neq
y}}\frac{d_X(x,y)}{d_Y(f(x),f(y))}.
$$
We denote by $c_Y(X)$ the least distortion with which $X$ may be
embedded in $Y$. When $c_Y(X)\le \alpha$ we say that $X$
$\alpha$-embeds into $Y$ and denote
$X\stackrel{\alpha}{\hookrightarrow}Y$. When there is a bijection
$f$ between two metric spaces $X$ and $Y$ with
$\mathrm{dist}(f)\le \alpha$ we say that $X$ and $Y$ are
$\alpha$-equivalent. For a class of metric spaces $\M$, $c_\M(X)$
is the minimum $\alpha$ such that $X$ $\alpha$-embeds into some
metric space in $\M$. For $p\ge 1$ we denote $c_{\ell_p}(X)$ by
$c_p(X)$. The parameter $c_2(X)$ is known as the {\em Euclidean
distortion} of $X$. A fundamental result of Bourgain
\cite{bourgainembedding} states that $c_2(X)=O(\log n)$ for every
$n$-point metric space $(X,d)$.

A metric Ramsey-type theorem states that a given metric space
contains a large subspace which can be embedded with small
distortion in some ``well-structured" family of metric spaces
(e.g., Euclidean). This can be formulated using the following
notion:

\begin{definition}[Metric Ramsey functions]
Let $\M$ be some class of metric spaces. For a metric space $X$,
and $\alpha \geq 1$, $R_\M(X;\alpha)$ denotes the {\it largest size} of
a subspace $Y$ of $X$ such that $c_\M(Y) \leq \alpha$.

Denote by $R_{\M}(\alpha,n)$ the largest integer $m$ such that any
$n$-point metric space has a subspace of size $m$ that
$\alpha$-embeds into a member of $\M$. In other words, it is the
infimum over $X$, $|X|=n$, of $R_\M(X;\alpha)$.

It is also useful to have the following conventions: For
$\alpha=1$ we allow omitting $\alpha$ from the notation. When $\M
= \{ X \}$, we write $X$ instead of $\M$. Moreover when
$\M=\{\ell_p\}$, we use $R_p$ rather than $R_{\ell_p}$.

In the most general form, let $\N$ be a class of metric spaces and
denote by $R_\M(\N;\alpha,n)$ the largest integer $m$ such that
any $n$-point metric space in $\N$ has a subspace of size $m$ that
$\alpha$-embeds into a member of $\M$. In other words, it is the
infimum over $X \in \N$, $|X|=n$, of $R_\M(X;\alpha)$.
\end{definition}

\vglue-12pt
\Subsec{Results for arbitrary metric
spaces}
This paper provides several results concerning metric Ramsey
functions. One of our main objectives is to provide bounds on the
{\em Euclidean Ramsey Function}, $R_2(\alpha,n)$.

The first result on this problem, well-known as a nonlinear
version of Dvoretzky's theorem, is due to Bourgain, Figiel and
Milman~\cite{bfm}:
\begin{theorem}[\cite{bfm}] \label{thm:bfm86}
For any $\alpha>1$ there exists $C(\alpha)>0$ such that $R_2(\alpha,n)  
\geq
C(\alpha) \log n$. Furthermore{\rm ,} there exists $\alpha_0>1$ such that
$R_2(\alpha_0,n)=O(\log n)$.
\end{theorem}

While Theorem~\ref{thm:bfm86} provides a tight characterization of
$R_2(\alpha,n) = \Theta(\log n)$ for values of $\alpha \le
\alpha_0$ (close to $1$), this bound turns out to be very far from
the truth for larger values of $\alpha$ (in fact, a careful
analysis of the arguments in \cite{bfm} gives $\alpha_0\approx
1.023$, but as we later discuss, this is not the right threshold).

Motivated by problems in the field of Computer Science, more
researchers \cite{krr}, \cite{bkrs}, \cite{bbm} have investigated metric Ramsey
problems. A close look (see \cite{bbm}) at the results of
\cite{krr}, \cite{bkrs} as well as \cite{bfm} reveals that all of these
can be viewed as based on Ramsey-type theorems where the target
class is the class of ultrametrics (see
\S\ref{section:lower-preliminaries} for the definition).

The usefulness of such results for embeddings in $\ell_2$ stems
from the well-known fact \cite{lemin} that ultrametrics are
isometrically embeddable in $\ell_2$. Thus, denoting the class of
ultrametrics by $\um$, we have that $R_2(\alpha,n) \ge
R_{\um}(\alpha,n)$.

The recent result of Bartal, Bollob\'as and Mendel \cite{bbm}
shows that for large distortions the metric Ramsey function
behaves quite differently from the behavior expressed by
Theorem~\ref{thm:bfm86}. Specifically, they prove that
$R_2(\alpha,n) \ge R_{\um}(\alpha,n) \ge \exp \left((\log
n)^{1-O(1/\alpha)} \right)$ (in fact, it was already implicit in
\cite{bkrs} that a similar bound holds for a particular $\alpha$).
The main theorem in this paper is:

\begin{theorem}[Metric Ramsey-type theorem]\label{thm:lowerlarge}
For every $\e >0${\rm ,} any\break $n$\/{\rm -}\/point metric space has a subset of size
$n^{1-\e}$ which embeds in Hilbert space with distortion
$O\left(\frac{\log(1/\e)}{\e}\right)$. Stated in terms of the
metric Ramsey function{\rm ,} there exists an absolute constant $C>0$
such that for every $\alpha>1$ and every integer~$n$\/{\rm :}\/
$$
  R_2(\alpha,n) \ge R_{\um}(\alpha,n) \ge  
n^{1-C\frac{\log(2\alpha)}{\alpha}}.
$$
\end{theorem}

We remark that the lower bound above for $R_{\um}(\alpha,n)$ is
meaningful only for large enough $\alpha$. Small distortions are
dealt with in Theorem~\ref{thm:phase} (see also
Theorem~$3.26$).

The fact that the subspaces obtained in this Ramsey-type theorem
are ultrametrics in not just an artifact of our proof. More
substantially, it is a reflection of new embedding techniques that
we introduce. Indeed, most of the previous results on embedding
into $\ell_p$ have used what may be called Fr\'echet-type
embeddings: forming coordinates by taking the distance from a
fixed subset of the points. This is the way an arbitrary finite
metric space is embedded in $\ell_\infty$ (attributed to
Fr\'echet). Bourgain's embedding~\cite{bourgainembedding} and its
generalizations~\cite{matexpander} also fall in this category of
embeddings. However, it is possible to show that Fr\'echet-type
embeddings are not useful in the context of metric Ramsey-type
problems. More specifically, we show in~\cite{blmn2} that such
embeddings cannot achieve bounds similar to those of
Theorem~\ref{thm:lowerlarge}.

Ultrametrics have a useful representation as {\em hierarchically
well-separated trees} (HST's). A $k$-HST is an ulrametric where
vertices in the rooted tree are labelled by real numbers. The
labels decrease by a factor $\geq k$ as you go down the levels
away from the root. The distance between two leaves is the label
of their lowest common ancestor. These {\em decomposable} metrics
were introduced by Bartal \cite{bartal1}. Subsequently, it was
shown (see \cite{bartal1}, \cite{bartal2}, \cite{frt}) that any $n$-point metric
can be $O(\log n)$-probabilistically embedded\footnote{A metric
space can be $\alpha$-probabilistically embedded in a class of
metric spaces if it is $\alpha$-equivalent to a convex combination
of metric spaces in the class, via a noncontractive Lipschitz
embedding \cite{bartal2}.} in ultrametrics. This theorem has found
many unexpected algorithmic applications in recent years, mostly
in providing computationally efficient approximate solutions for
several NP-hard problems (see the survey \cite{indyk} for more
details).

The basic idea in the proof of Theorem~\ref{thm:lowerlarge} is to
iteratively find large subspaces that are hierarchically
structured, gradually improving the distortion between these
subspaces and a hierarchically well-separated tree. These
hierarchical structures are naturally modelled via a notion (which
is a generalization of the notion of $k$-HST) we call \emph{metric
composition} closure. Given a class of metric spaces $\M$, we obtain a
metric space in the class $\comp_k(\M)$  by taking a
metric space $M \in \M$ and replacing its points with copies of
metric spaces from $\comp_k(\M)$ dilated so that there is a factor
$k$ gap between distances in $M$ and distances within these
copies.

Metric compositions are also used to obtain the following bounds
on the metric Ramsey function in its more general form:

\begin{theorem}[Generic bounds on the metric Ramsey  
function]\label{thm:property} Let
$\cal C$\break be a proper class of finite metric spaces that is closed
under\/{\rm : (i)}\/ Isometry\/{\rm , (ii)}\/ Passing to a subspace{\rm , (iii)} Dilation.
Then there exists $\delta<1$ such that $R_{\cal C}(n)\le n^\delta$
for infinitely many values of $n$.
\end{theorem}

In particular we can apply Theorem~\ref{thm:property} to the class
$\cal C \rm =\{X;\ c_\M(X) \leq \alpha \}$ where $\M$ is some
class of metric spaces. If there exists a metric space $Y$ with
$c_\M(Y)>\alpha$, then  there exists $\delta<1$ such
that $R_\M(\alpha,n) < n^\delta$ for infinitely many $n$'s.

In the case of $\ell_2$ or ultrametrics much better bounds are
possible, showing that the bound in Theorem~\ref{thm:lowerlarge}
is almost tight. For ultrametrics this is a rather simple fact
\cite{bbm}. For embedding into $\ell_2$ this follows from bounds
for expander graphs, described later in more detail.

\begin{theorem}[near tightness]\label{thm:upperlarge} There exist  
absolute
constants $c,C>0$ such that for every $\alpha>2$ and every integer
$n$\/{\rm :}\/
$$
  R_{\um}(\alpha,n) \le R_2(\alpha,n) \le C n^{1-\frac{c}{\alpha}}.
$$
\end{theorem}

The behavior of $R_{\um}(\alpha,n)$ and $R_2(\alpha,n)$ exhibited
by the bounds in Theorems~\ref{thm:bfm86} and~\ref{thm:lowerlarge}
is very different. Somewhat surprisingly, we discover the
following phase transition:

\begin{theorem}[phase transition]\label{thm:phase} For every
$\alpha>1$ there exist constants $c,C,c',C',K>0$ depending only on
$\alpha$ such that $0<c'<C'<1$ and for every integer $n$\/{\rm :}\/

\begin{itemize}
\item[{\rm a)}] If $1 < \alpha < 2$ then  $c\log n\le
R_{\um}(\alpha,n) \le R_2(\alpha,n)\le 2\log_2 n +C$.

\item[{\rm b)}] If  $\alpha>2$ then $n^{c'}\le
R_{\um}(\alpha,n) \le R_2(\alpha,n)\le K\, n^{C'}$.
\end{itemize}
\end{theorem}

  Using bounds on the dimension with which any $n$ point
ultrametric is embeddable with constant distortion in $\ell_p$
\cite{blmn4} we obtain the following corollary:

\begin{corollary}[Ramsey-type theorems with low  
dimension]\label{thm:plowdim}
There exists $0<C(\alpha)<1$ such that for every $p \ge 1${\rm ,}
$\alpha>2${\rm ,} and every integer $n${\rm ,}
$$
R_{\ell_p^d}(\alpha,n) \ge n^{C(\alpha)},
$$
where $C(\alpha) \ge 1- \frac{c\log \alpha}{\alpha}${\rm ,} $d =
\bigl\lceil \lceil \frac{c'}{(\alpha-2)^2} \rceil C(\alpha)\log n
\bigr\rceil$, and $c,c'>0$ are universal constants.
\end{corollary}

This result is meaningful since, although $\ell_2$ isometrically
embeds\break into $L_p$ for every $1\le p\le \infty$, there is no known
$\ell_p$ analogue of the Johnson-Lindenstrauss dimension reduction
lemma \cite{jl} (in fact, the Johnson-\break Lindenstrauss lemma is known
to fail in $\ell_1$~\cite{bc}, \cite{ln}). These bounds are almost best
possible.

\medbreak
\scshape{Theorem 1.8} (\rm The Ramsey problem for finite dimensional normed  
spaces). 
{\it There exist absolute constants $C,c>0$ such that for any
$\alpha>2${\rm ,} every integer $n$ and every finite dimensional normed
space} $X${\rm ,}
  $$ R_{X}(\alpha,n) \leq  C n^{1-\frac{c}{\alpha}} (\dim X) \log  
\alpha. $$
\bigbreak

\setcounter{lemma}{8}

For completeness, we comment that a natural question, in our
context, is to bound the size of the  largest subspace of an arbitrary
finite metric space that is isometrically embedded in $\ell_p$. In
\cite{blmn5} we show that $R_p(n) = 3$ for every $1<p<\infty$ and
$n\ge 3$.

Finally, we note that one important motivation for this work is
the applicability of metric embeddings to the theory of
algorithms. In many practical situations, one encounters a large
body of data, the successful analysis of which depends on the way
it is represented. If, for example, the data have a natural metric
structure (such as in the case of distances in graphs), a low
distortion embedding into some normed space helps us draw on
geometric intuition in order to analyze it efficiently. We refer
to the papers \cite{bartal2}, \cite{feige}, \cite{llr} and the surveys
\cite{indyk}, \cite{linial} for some of the applications of metric
embeddings in Computer Science. More about the relevance of
Theorem~\ref{thm:lowerlarge} to Computer Science can be found in
\cite{blmn1} (see also \cite{bbm}, \cite{bm}).

\Subsec{Results for special classes of metric spaces}
We provide nearly tight bounds for concrete families of metric
spaces: expander graphs, the discrete cube, and high girth graphs.
In all cases the difficulty is in providing upper bounds on the
Euclidean Ramsey function.

  Let $G=(V,E)$ be a $d$-regular graph, $d\ge 3$, with
absolute multiplicative spectral gap $\gamma$ (i.e.\ the second
largest eigenvalue, in absolute value, of the \pagebreak adjacency matrix of
$G$ is less than $\gamma d$). For such expander graphs it is
known~\cite{llr}, \cite{matexpander} that $c_2(G)=\Omega_{\gamma,d}(\log
|V|)$ (here, and in what follows, the notation $a_n=\Omega(b_n)$
means that there exists a constant $c>0$ such that for all $n$,
$|a_n|\ge c|b_n|$. When $c$ is allowed to depend on, say, $\gamma$
and $d$ we use the notation $\Omega_{\gamma,d}$). In
Section~\ref{section:expanders} we prove the following:

\begin{theorem}[The metric Ramsey problem for expanders]\label{thm:expander}
Let $G=\break (V,E)$ be a $d$-regular graph{\rm ,} $d\ge 3$ with absolute
multiplicative spectral gap $\gamma<1$. Then  for every
$p\in[1,\infty)${\rm ,} and every $\alpha\ge 1${\rm ,}
$$
|V|^{1-\frac{C}{\alpha\log_d(1/\gamma)}}\le R_2(G;\alpha)\le
R_p(G;\alpha)  \le Cd|V|^{1-\frac{c\log_d(1/\gamma)}{p\alpha}},
$$
where $C,c>0$ are absolute constants.
\end{theorem}

The proof of the upper bound in Theorem~\ref{thm:expander}
involves proving certain Poincar\'e inequalities for power graphs
of $G$.

Let $\Omega_d=\{0,1\}^d$ be the discrete cube equipped with the
Hamming metric. It was proved by Enflo, \cite{enflocube}, that
$c_2(\Omega_d) = \sqrt{d}$. Both Enflo's argument, and subsequent
work of Bourgain, Milman and Wolfson~\cite{bmw}, rely on
nonlinear notions of type. These proofs strongly use the
structure of the whole cube, and therefore seem not applicable for
subsets of the cube. In Section~\ref{section:cube} we prove the
following strengthening of Enflo's bound:

\begin{theorem}[The metric Ramsey problem for the
discrete cube]\label{thm:cube}\hfill\break There exist absolute constants
$C,c$ such that for every $\alpha>1${\rm :}
$$ 2^{d\left(1-\frac{\log (C \alpha)}{\alpha^2}\right)} \leq  
R_2(\Omega_d;\alpha) \leq
C2^{d\left(1-\frac{c}{\alpha^2}\right)} .$$
\end{theorem}
\vglue8pt

The lower bounds on the Euclidean Ramsey function mentioned above
are based on the existence of large subsets of the graphs which
are within distortion $\alpha$ from forming an equilateral space.
In particular for the discrete cube this corresponds to a code of
large relative distance. Essentially, our upper bounds on the
Euclidean Ramsey function show that for a fixed size, no other
subset achieves significantly better distortions.

In \cite{lmn} it was proved that if $G=(V,E)$ is a $d$-regular graph,  
$d\ge
3$, with girth $g$, then $c_2(G)\ge c\frac{d-2}{d}\sqrt{g}$. In
Section~\ref{section:girth} we prove the following strengthening of this
result:

\begin{theorem}[The metric Ramsey problem for large girth  
graphs]\label{thm:girth} \hfill\break
Let $G= (V,E)$ be a $d$-regular graph{\rm ,} $d\ge 3${\rm ,} with girth $g$.
Then for every $1\le \alpha<\frac{\sqrt{g}}{6}${\rm ,}
       $$ R_2(G;\alpha) \leq C(d-1)^{-\frac{cg}{\alpha^2}}|V| ,$$
where $C,c>0$ are absolute constants. 
\pagebreak
\end{theorem}

The proofs of Theorem~\ref{thm:cube} and Theorem~\ref{thm:girth}
use the notion of Markov type, due to K.  Ball \cite{ball}. In
addition, we need to understand the algebraic properties of the
graphs involved (Krawtchouk polynomials for the discrete cube and
Geronimus polynomials in the case of graphs with large girth).

\section{Metric composition}\label{section:composition}

In this section we introduce the notion of metric composition,
which plays a basic role in proving both upper and lower bounds on
the metric Ramsey problem. Here we introduce this construction and
use it to derive some nontrivial upper bounds. The bounds
achievable by this method are generally not tight. For the Ramsey
problem on $\ell_p$, better upper bounds are given in
Sections~\ref{section:counting} and~\ref{section:expanders}. In
Section~\ref{section:lower bounds} we use metric composition in
the derivation of lower bounds.

\Subsec{The basic definitions}
\vglue-12pt
\begin{definition}[Metric composition]\label{def:metric-composition}
Let $M$ be a finite metric space. Suppose that there is a
collection of disjoint finite metric spaces $N_x$ associated with
the elements $x$ of $M$. Let $\N = \{ N_x \}_{x \in M}$. For
$\beta\geq 1/2$, the $\beta$-composition of $M$ and $\N$, denoted
by $C=M_\beta[\N]$, is a metric space on the disjoint union $\dot
{\boldsymbol{\cup}}_x N_x$. Distances in $C$ are defined as follows. Let $x,y \in
M$ and $u \in N_x, v \in N_y$; then:
$$ d_C(u,v)= \begin{cases} d_{N_x}(u,v) & x=y \\
   \beta \gamma d_M(x,y) & x\neq y,\end{cases}  $$
where $\gamma=\frac{\max_{z \in M} \diam(N_z)}{\min_{x\neq y \in
M} d_M(x,y)}$. It is easily checked that the choice of the factor
$\beta\gamma$ guarantees that $d_C$ is indeed a metric. If all the
spaces $N_x$ over $x\in M$ are isometric copies of the same space
$N$, we use the simplified notation $C=M_\beta[N]$.
\end{definition}

Informally stated, a metric composition is created by first
multiplying the distances in $M$ by $\beta \gamma$, and then
replacing each point $x$ of $M$ by an isometric copy of $N_x$.

A related notion is the following:

\begin{definition}[Composition closure] \label{def:comp}
Given a class $\M$ of finite metric spaces, we consider
$\comp_\beta(\M)$, its closure under $\ge \beta$-compositions.
Namely, this is the smallest class $\C$ of metric spaces that
contains all spaces in $\M$, and satisfies the following
condition: Let $M \in \M$, and associate with every $x \in M$ a
metric space $N_x$ that is isometric to a space in $\C$. Also, let
$\beta'\geq \beta$. Then $M_{\beta'}[\N]$ is also in $\C$.
\end{definition}
 
\Subsec{Generic upper bounds via metric composition}
We need one more  definition:

\begin{definition} A class $\cal C$ of finite metric spaces is called a
\emph{metric class} if it is closed under isometries. $\cal C$ is
said to be \emph{hereditary}, if $M\in \cal C$ and $N\subset M$
imply $N\in \cal C$. The class is said to be \emph{dilation
invariant} if $(M,d)\in \cal C$ implies that $(M,\lambda d)\in
\cal C$ for every $\lambda>0$.
\end{definition}

Let $\classeqv{\M}{\alpha} = \{ X; c_{\M}(X) \leq \alpha \}$ denote the
class of all metric spaces that $\alpha$-embed into some metric
space in $\M$. Clearly, $\classeqv{\M}{\alpha}$  is a hereditary,
dilation-invariant metric class.

We recall that $R_{\cal C}(X)$ is the largest cardinality of a
subspace of $X$ that is isometric to some metric space in the
class $\cal C$.

\begin{proposition} \label{prop:p0}
Let $\cal C$ be a hereditary{\rm ,} dilation invariant metric class of
finite metric spaces. Then{\rm ,} for every finite metric space $M$ and
a class $\N = \{ N_x \}_{x \in M}${\rm ,} and every $\beta \geq 1/2${\rm ,}
$$
R_{\cal C}(M_\beta[\N]) \le R_{\cal C}(M) \cdot \max_{x \in M}
R_{\cal C}(N_x).
$$
In particular{\rm ,} for every finite metric space $N${\rm ,}
$$
R_{\cal C}(M_\beta[N]) \le R_{\cal C}(M) R_{\cal C}(N).
$$
\end{proposition}

\Proof  Let $m=R_{\cal C}(M)$ and $k= \max_{x \in M} R_{\cal
C}(N_x)$. Fix any $X \subseteq \dot{\boldsymbol{\cup}}_x N_x$ with $|X|>mk$. For
every $z\in M$ let $X_z = X\cap N_z$. Set $Z = \{ z \in M; X_z
\neq \emptyset \}$. Note that $|X|=\sum_{z\in Z} |X_z|$ so that if
$|Z|\le m $ then there is some $y \in M$ with $|X_y|>k$. In this
case, the set $X_y$ consists of more than $k$ elements in $X$, the
metric on which is isometric to a subspace of $N_y$, and therefore
is not in $\cal C$. Since ${\cal C}$ is hereditary this implies
that $X\notin {\cal C}$.
  Otherwise, $|Z|>m$. Fix for each $z\in Z$
some arbitrary point $u_z\in X_z$ and set $Z' =\{u_z ;\ z\in Z\}$.
Now, $Z'$ consists of more than $m$ elements in $X$, the metric on
which is a $\beta \gamma$-dilation of a subspace of $M$, hence not
in $\cal C$. Again, the fact that ${\cal C}$ is hereditary implies
that $X\notin {\cal C}$.
\Endproof\vskip4pt 

In what follows let $R_\C(\A,n) = R_\C(\A;1,n)$. Recall that
$R_\C(\A;1,n) \ge t$ if and only if for every $X \in \A$ with $|X| = n$,
there is a subspace of $X$ with $t$ elements that is isometric to
some metric space in the class $\cal C$.

\begin{lemma} \label{lem:p1}
Let $\cal C$ be a hereditary{\rm ,} dilation invariant metric class of
finite metric spaces. Let $\A$ be a class of metric spaces{\rm ,} and
let $\delta\in(0,1)$. If there exists an integer $m>1$ such that
$R_\C(\A,m) \leq m^\delta${\rm ,} then for any $\beta \geq 1/2${\rm ,} and
infinitely many integers $n$\/{\rm :}\/
$$ R_\C(\comp_\beta(\A),n) \le n^\delta . $$
\end{lemma}

\Proof  Fix some $\beta \geq 1/2$.
Let $A \in \A$ be such that $|A| = m>1$ and $R_\C(A) \le m^\delta$.
Define inductively a sequence of metric spaces in
$\comp_\beta(\A)$ by: $A_1=A$ and $A_{i+1}=A_\beta[A_{i}]$.
Proposition \ref{prop:p0} implies that $R_\C(A_{i+1}) \le
R_\C(A_i)R_\C(A) \le R_\C(A_i) m^\delta$. It follows that
$R_\C(A_i) \le m^{i\delta} = |A_i|^\delta$.
\hfill\qed

\begin{lemma} \label{lem:p1.1}
Let $\cal C$ be a nonempty hereditary{\rm ,} dilation invariant metric
class of finite metric spaces. Let $\A$ be a class of finite
metric spaces{\rm ,} such that $R_\C(\A,m) < m$ for some integer $m$
\/{\rm (}\/i.e.{\rm ,} there is some space $A \in \A$ with no isometric copy in
$\C${\rm ).} Then there exists $\delta\in(0,1)${\rm ,} such that for any
$\beta \geq 1/2${\rm ,} and infinitely many integers $n${\rm :}
$$ R_\C(\comp_\beta(\A),n) \le n^\delta . $$
\end{lemma}
\vglue8pt

\Proof 
Let $m$ be the least cardinality of a space $A \in \A$ of with no
isometric copy in $\cal C$. Since $\cal C$ is nonempty and
hereditary, $m \ge 2$. Define $\delta$ by $m-1 = m^{\delta}$. Now
apply Lemma \ref{lem:p1}.
\Endproof\vskip4pt 

Lemma~\ref{lem:p1.1} can be applied to obtain nontrivial bounds
on various metric Ramsey functions.

\begin{corollary} \label{cor:p2}
Let $\cal C$ be a hereditary, dilation invariant metric class
which contains some{\rm ,} but not all finite metric spaces. Then there
exists a $\delta\in(0,1)${\rm ,} such that $R_{\cal C}(n) \leq
n^{\delta}$ for infinitely many integers $n$.
\end{corollary}
\Proof 
We use Lemma~\ref{lem:p1.1} with $\A = \comp_\beta(\A) =$ the
class of all metric spaces.
\Endproof\vskip4pt 

Let $\mathcal{M}$ be a fixed class of metric spaces and $\alpha\ge
1$. The following corollary follows when we apply
Corollary~\ref{cor:p2} with $\mathcal{C}= \classeqv{\M}{\alpha}$ as
defined above.

\begin{corollary} \label{cor:p2.1}
Let $\M$ be a metric class of finite metric spaces and $\alpha\geq
1$. The following assertions are equivalent\/{\rm :}\/
\begin{itemize}
 \ritem{{\rm a)}} There exists an integer $n${\rm ,} such that
$R_{\M}(\alpha,n) <n$. \medskip

\ritem{ {\rm b)}} There exists $\delta\in(0,1)${\rm ,} such that
$R_{\M}(\alpha,n) \leq n^{\delta}$ for infinitely many integers
$n$.\end{itemize}
\end{corollary}

For our next result, recall that a normed space $X$
is said to have cotype $q$ if there is a positive constant $C$
such that for every finite sequence $x_1,\dots\break\dots ,x_m\in X$,
$$
\left(\E \left\|\sum_{i=1}^m \varepsilon_i
x_i\right\|^2\right)^{1/2}\ge C\left( \sum_{i=1}^m
\|x_i\|^q\right)^{1/q},
$$
where $\varepsilon_1,\dots , \varepsilon_m$ are i.i.d. $\pm 1$
Bernoulli random variables. It is well known (see
\cite{milschechtbook}) that for $2\le q<\infty$, $\ell_q$ has
cotype $q$ (and it does not  have cotype $q'$ for any $q'<q$).

\begin{corollary} \label{corol:lb}
Let $ X$ be a normed space. Then the following assertions are
equivalent\/{\rm :}\/
\begin{itemize}
\ritem{a)} $X$ has finite cotype.

\ritem{b)} For any $\alpha>1${\rm ,} there exists
$\delta\in(0,1)$  such that for infinitely many integers $n${\rm ,} $R_{
X}({\alpha},n)\leq n^\delta$.
  
\ritem{c)} There exists $\alpha>1$ and an integer $n$ such
that $R_{ X}({\alpha},n)< n$.
\end{itemize}
\end{corollary}
\Proof 
To prove the implication ${\rm a)}\Longrightarrow {\rm b)}$, fix
$\alpha>1$. Now, since $X$ has finite cotype, there is an integer
$h$ such that $d(\ell_\infty^h,Z) > \alpha$ for every\break
$h$-dimensional subspace $Z$ of $X$, where $d(\cdot,\cdot)$ is the
Banach-Mazur distance. This implies that for some $\epsilon>0$, an
$\epsilon$-net $\cal E$ in the unit ball of $\ell_\infty^h$ does
not $\alpha$-embed into $X$. This follows from a standard argument
in nonlinear Banach space theory. Indeed, a compactness argument
would imply that otherwise $B_{\infty}^h$, the unit ball of
$\ell_\infty^h$, $\alpha$-embeds into $X$. By Rademacher's
theorem  (see for example \cite{benlin}) such an embedding must be
differentiable in an interior point of $B_{\infty}^h$. The
derivative, $T$, is a linear mapping which is easily seen to
satisfy $\|T\|\cdot \|T^{-1}\|\le \alpha$, so that
$d(\ell_\infty^h,Z) \le \alpha$ for the subspace
$Z=T(\ell_\infty^h)$. Apply Corollary~\ref{cor:p2.1} with $\M=X$,
and $n=|{\cal E}|$ to conclude that {\rm b)} holds.

The implication 
${\rm b)}\Longrightarrow {\rm c)}$ is obvious, so we turn to prove
that ${\rm c)}\Longrightarrow {\rm a)}$. Assume that $X$ does not
have finite co-type, and fix some $0<\epsilon<\alpha-1$. By the
Maurey-Pisier theorem (see \cite{maureypisier} or Theorem 14.1 in
\cite{djt}), it follows that for every $n$, $\ell_\infty^n$
$(\alpha-\epsilon)$-embeds into $X$. Since $\ell_\infty^n$
contains an isometric copy of every $n$-point metric space, we
deduce that for each $n$, $R_{X}(\alpha,n)=n$, contrary to our
assumption ${\rm c)}$.
\Endproof\vskip4pt 

We now need the following variation on the theme of metric
composition.

\begin{definition} \label{def:nearly-close}
A family of  metric spaces $\N$ is called  \emph{nearly closed
under composition}, if for every $\lambda > 1$, there exists some
$\beta \ge 1/2$ such that $c_\N(X) \le \lambda$ for every $X \in
\comp_\beta(\N)$. In other words,
$$ \comp_\beta(\N) \subseteq \classeqv{\N}{\lambda} . $$
\end{definition}

We have the following variant of Corollary~\ref{cor:p2.1}:

\begin{lemma} \label{lem:p2+}
Let $\cal M$ be a metric class of finite metric spaces and let
$\N$ be some class of finite metric spaces which is nearly closed
under composition. Assume that there is some space in $\N$ which
does not $\alpha$-embed into any space in $\cal M$. Then there
exists $\delta\in(0,1)${\rm ,} such that for every $1 \le \alpha' <
\alpha${\rm ,} $R_{\cal M}(\N;\alpha' ,n) \leq n^{\delta}$ for
infinitely many integers $n$.
\end{lemma}
\Proof 
Fix some $\alpha'<\alpha$ and let $\lambda = \alpha/\alpha'$. As
$\N$ is nearly closed under composition there exists $\beta \ge
1/2$ such that $\comp_\beta(\N) \subseteq \classeqv{\N}{\lambda}$.
This means that for every $Z \in \comp_\beta(\N)$ there exists
some $N \in \N$ that is $\lambda$-equivalent to~$Z$.

For every integer $p$ let $k(p) = R_\M(\N;\alpha',p)$. If
$|Z|=|N|=n$, then there is $X \subseteq N$ such that $c_\M(X) \le
\alpha'$ and $|X| \ge k(n)$. Let $Y \subseteq Z$ be the set
corresponding to $X$ under the $\lambda$-equivalence between $Z$
and $N$. Then, $|Y| = |X| \ge k(n)$ and by composition of maps,
$c_\M(Y) \le \lambda\alpha' = \alpha$. That is, every\break $n$-set $Z$
in $\comp_\beta(\N)$ contains a $k(n)$ subset $Y$ that
$\alpha$-embeds into a space in~$\M$; i.e.\ $Y \in
\classeqv{\M}{\alpha}$. In our notation, this means that $k(n) \le
R_\C(\comp_\beta(\N),n)$, where $\C = \classeqv{\M}{\alpha}$.

The assumption made in the lemma about $\N$ means that $R_\C(\N,m)
< m$ for some integer $m$. By Lemma~\ref{lem:p1.1} there exists
$\delta \in (0,1)$ such that for infinitely many integers $n$,
$$
R_\M(\N;\alpha',n) = k(n) \le R_\C(\comp_\beta(\N),n) < n^\delta,
$$
as claimed.
\Endproof\vskip4pt 

Next, we give a several results that demonstrate the applicability
of\break Lemma~\ref{lem:p2+}.

\begin{proposition} \label{prop:normed}
Let $(X,\|\cdot \|)$  be a normed space. The class $\M$ of finite
subsets of $X$ is nearly closed under composition.
\end{proposition}

\Proof 
Fix some $\lambda>1$. Let $Z \in \comp_\beta(\M)$ for some $\beta
> 1/2$ to be determined later.
We prove that $Z$ can be $\lambda$-embedded in $X$. The proof is
by induction on the number of steps taken in composing $Z$ from
spaces in $\M$. If $Z \in \M$ there is nothing to prove.
Otherwise, it is possible to express $Z$ as $Z = M_\beta[\N]$,
where $M \in \M$ and $\N = \{ N_z \}_{z \in M}$ such that each of
the spaces $N_z$ is in $\comp_\beta(\M)$ and can be created by a
shorter sequence of composition steps. By induction we assume that
there exists $\beta$ for which $N_z$ can be $\lambda$-embedded in
$X$. Fix for every $z \in M$, $\phi_z:N_z\to X$ satisfying:
$$
\forall u,v\in N_{z}, ~~~d_{N_z}(u,v)\le \|\phi_z(u)-\phi_z(v)\| \le
\lambda d_{N_z}(u,v) ,
$$ and for all $u \in N_z$, $\| \phi_z(u) \| \leq \lambda \diam(N_z)$  
(this
can be assumed by an appropriate translation).

Define $\phi: Z \to X$ as follows: for every $u \in Z$, let $z \in
M$ be such that $u \in N_z$, then $\phi(u)= \beta\gamma \cdot z
+\phi_z(u)$, where $\gamma=\frac{\max_z \diam(N_z)}{\min_{x\neq y
\in M} \| x -y \| }$.

We now bound the distortion of $\phi$. Assume $\beta > 2\lambda$.
Consider first $u,v \in N_z$ for some $z \in M$.
$$ d_Z(u,v)=d_{N_z}(u,v)\le\| \phi_z(u)- \phi_z(v) \|\le
\lambda d_{N_z}(u,v)= \lambda d_Z(u,v). $$

Now, let $u \in N_x, v \in N_y$, for $x \neq y \in M$,

\begin{align*}
\| \phi(u)- \phi(v) \| &\leq \beta \gamma \| x - y \|
   + \| \phi_x(u) - \phi_y(v) \| \\
&\leq \beta \gamma \| x - y \| + \lambda (\diam(N_x)+\diam(N_y)) \\
& \leq \gamma(\beta + 2\lambda) \| x - y \| = \frac{\beta +
2\lambda}{\beta} d_Z(u,v).
\end{align*}
Similarly,
\begin{align*}
\| \phi(u)- \phi(v) \| &\geq \beta \gamma \| x - y \|
   - \| \phi_x(u) - \phi_y(v) \| \\
&\geq \beta \gamma \| x - y \| - \lambda (\diam(N_x)+\diam(N_y)) \\
& \leq \gamma(\beta - 2\lambda) \| x - y \|  = \frac{\beta -
2\lambda}{\beta} d_Z(u,v).
\end{align*}

Hence if $\beta \ge 2\lambda\frac{\lambda+1}{\lambda-1}$, we have,
\vskip12pt
\hfill $
\displaystyle{\mathrm{dist}(\phi)\le \max \left\{\lambda,\frac{
\beta+2\lambda}{\beta-2\lambda}\right\} = \lambda .}
$ 
\Endproof\vskip12pt

Recall that a normed space $X$ is said to be $\lambda$ finitely
representable in a normed space $Y$ if for any finite dimensional
linear subspace $Z\subset X$ and every $\eta>0$ there is a
subspace $W$ of $Y$ such that $d(Z,W)\le \lambda+\eta$.

\begin{corollary}\label{coro:finrep} Let $X$ and $Y$ be normed spaces  
and
$\alpha>1$. The following are equivalent\/{\rm :}\/

\begin{itemize}
\ritem{1)} $X$ is not $\alpha$-finitely representable in
$Y$.

\ritem{2)} There are $\eta>0$ and $\delta\in (0,1)$ such
that $R_{Y}(X;\alpha+\eta,n) < n^\delta$ for infinitely many
integers $n$.

\ritem{3)} There is some $\eta>0$ and an integer $n$ such
that $R_{Y}(X;\alpha+\eta,n)<n$.
\end{itemize}
\end{corollary}

\Proof 
If $X$ is not $\alpha$-finitely representable in $Y$ then there is
a finite dimensional linear subspace $Z$ of $X$ whose Banach-Mazur
distance from any subspace of $Y$ is greater than $\alpha$. As in
the proof of Corollary \ref{corol:lb}, a combination of a
compactness argument and a differentiation argument imply that
there is a finite subset $S$ of $X$ which does not
$(\alpha+2\eta)$ embed in $Y$ for some $\eta>0$. Since the subsets
of $X$ are nearly closed under composition, by applying
Lemma~\ref{lem:p2+}, we deduce the implication ${\rm
1)}\Longrightarrow {\rm 2)}$.

The implication ${\rm 2)}\Longrightarrow {\rm 3)}$ is obvious, so
we turn to show ${\rm 3)}\Longrightarrow {\rm 1)}$. Let $A \subset
X$ be a finite subset that does not $\alpha+\eta$ embed in $Y$,
and let $Z$ be $A$'s linear span. Clearly $d(Z,W) > \alpha + \eta$
for any linear subspace $W$ of $Y$. It follows that $X$ is not
$\alpha$-finitely representable in $Y$.
\Endproof\vskip4pt 

Recall that a graph $H$ is called a {\em minor} of a graph $G$ if
$H$ is obtained from $G$ by a sequence of steps, each of which is
either a contraction or a deletion of an edge. We say that a
family $\F$ of graphs is minor-closed if it is closed under taking
minors. The Wagner conjecture famously proved by Robertson and
Seymour~\cite{rob_seym}, states that for any nontrivial
minor-closed family of graphs~$\F$, there is a {\em finite} set of
graphs, $\mathcal{H}$, such that $G\in \F$ if and only if no
member of $\mathcal{H}$ is a minor of $G$. We say then that $\F$
is characterized by the list $\mathcal{H}$ of forbidden minors.
For example, planar graphs are precisely the graphs which do not
have $K_{3,3}$ or $K_5$ as minors, and the set of all trees is
precisely the set of all connected graphs with no $K_3$ minor.

There is a graph-theoretic counterpart to composition. Namely, let
$G = (V,E)$ be a graph, and suppose that to every vertex $x \in V$
corresponds a graph $H_x=(V_x,E_x)$ with a marked vertex $r_x \in
V_x$, where the $H_x$ are disjoint. The corresponding {\em graph
composition}, denoted $G[\{ H_x \}_{x \in V} ]$, is a graph with
vertex set $\dot{\boldsymbol{\cup}}_{x \in V} V_x$, and edge set:
$$
E=\{[u,v];\ x\in V, \ [u,v]\in E_x\} \cup \{[r_{x},r_{y}];\ [x,
y]\in E\}.
$$
The composition closure of a family of graphs $\F$ can be defined
similarly to Definition~\ref{def:comp}, and family $\F$ is said to
be closed under composition if it equals its closure.

Recall that a connected graph $G$ is called bi-connected if it
stays connected after we delete any single vertex from $G$ (and
erase all the edges incident with it). The maximal bi-connected
subgraphs of $G$ are called its {\em blocks}.

We make the following elementary graph-theoretic observation:
\begin{proposition}
Let $H$ be a bi-connected graph \/{\rm (}\/with $\ge 3$ vertices\/{\rm )}\/ that is a
minor of a graph $G$. Then $H$ is a minor of a block of $G$.
\end{proposition}

\Proof 
Consider a sequence of steps in which edges in $G$ are being
shrunk to form $H$. If there are two distinct blocks $B_1, B_2$ in
$G$ that are not shrunk to a single vertex, then the resulting
graph is not bi-connected. Indeed, there is a cut-vertex $a$ in
$G$ that separates $B_1$ from $B_2$, and the vertex into which $a$
is shrunk still separates the shrunk versions of $B_1, B_2$. This
observation means that in shrinking $G$ to $H$, only a single
block $B$ of $G$ retains more than one vertex. But then $H$ is a
minor of $B$, as claimed.
\Endproof\vskip4pt

In the graph composition described above, each vertex $r_x \in
V_x$ is a cut vertex. Consequently, each block of the composition
is either a block of $G$ (the subgraph induced by the vertices
$\{r_x;x \in V\}$ is isomorphic with $G$) or of one of the $H_x$
(that is isomorphic with the subgraph induced on $V_x$). We
conclude:

\begin{proposition} \label{prop:graph_minor}
Let $\F$ be a minor\/{\rm -}\/closed family of graphs characterized by a
list of bi\/{\rm -}\/connected forbidden minors. Then $\F$ is closed under
graph composition.
\end{proposition}

Let $\F$ again be a family of undirected graphs. A metric space
$M$ is said to be {\em supported} on $\F$ if there exist  a graph
$G\in \F$ and positive weights on the edges of $G$ such that $M$
is the geodetic, or shortest path metric on a subset of the
vertices of the weighted $G$.

Here is the metric counterpart of
Proposition~\ref{prop:graph_minor}:

\begin{proposition} \label{prop:minor}
Let $\F$ be a minor\/{\rm -}\/closed family of graphs characterized by a
list of bi\/{\rm -}\/connected forbidden minors. Then the class of metrics
supported on $\F$ is nearly closed under composition.
\end{proposition}

\Proof 
Fix some $\lambda>1$. Let $\F'$ be the class of metrics supported
on $\F$. Let $X \in \comp_\beta(\F')$ for some $\beta
> 1/2$ to be determined later.
We prove that $X$ can be $\lambda$-embedded in $\F'$. The proof is
by induction on the number of steps taken in composing $X$ from
spaces in $\F'$. If $X \in \F'$ there is nothing to prove.

  Otherwise, there exists a weighted graph $G=(V,E,w)$ in
$\F$. For simplicity, we identify $G$ with a metric space in
$\F'$, equipped with the geodetic metric defined by its weights.
It is possible to express $X$ as $X = G_\beta[\HH']$, where $\HH'
= \{ H_z' \}_{z \in V}$ such that each of the metric spaces $H_z'$
is in $\comp_\beta(\F')$. By in-\break duction we assume that there exists
$\beta$ for which each $H_z'$ can be\break $\lambda$-embedded in $\F'$.
Therefore there exists a family of disjoint weighted graphs
$\{H_z=(V_z,E_z,w_z)\}_{z \in V}$, such that for every $z\in V$, there
is a noncontractive Lipschitz bijection, $\phi_z:H_z'\to V_z$,
satisfying for any $u,v \in H_z'$, $d_{H_z'}(u,v) \le
d_{H_z}(\phi_z(u),\phi_z(v)) \le \lambda d_{H_z'}(u,v)$.

Let $Y=G[\{H_z\}_{z \in V}]$ be the graph composition of the above
graphs. Define weights $w'$ on the edges of $Y$ as follows: For
any $z \in V$, $[u,v]\in E_z$, let $w'([u,v])=w_z([u,v]).$ For
$[x,y]\in E$, let $w'([r_{x},r_{y}]) = \beta \gamma w([x,y]),$
where $\gamma = \frac{\max_{z \in V} \diam(H_z')}{\min_{x\neq y
\in V} d_G(x,y)}$ (as in the definition of metric composition).
For simplicity,\break\vglue-11pt\noindent we identify $Y$ with the weighted graph defined
above as well as the geodetic metric defined by this graph. The
proof shows that if $\beta$ is large enough, then the geodetic
metric on the graph composition $Y$ is $\lambda$-equivalent (and
thus arbitrarily close) to the metric $\beta$-composition $X$.
   Proposition~\ref{prop:graph_minor}
implies that $Y$ belongs to $\F'$, which proves the claim.

Indeed, define the bijection $\phi:X\to \dot{\boldsymbol{\cup}}_{u \in V} V_u$
as follows: for $z \in V$, if $u \in H_z'$, then
$\phi(u)=\phi_z(u)$. The geodetic path between any two vertices
$u',v' \in V_z$ is exactly the same path as in $H_z$, since the
cost of every step outside of $V_z$ exceeds $\diam(H_z)$ (by
definition of $\gamma$). This implies that
\begin{eqnarray*}
d_X(u,v)=d_{H_z'}(u,v)& \le& d_{H_z}(\phi_z(u),\phi_z(v)) \\
&=&
d_{Y}(\phi(u),\phi(v)) \le \lambda d_{H_z'}(u,v) = \lambda
d_X(u,v).
\end{eqnarray*}

Also, the distance in the graph composition between $u' \in V_x$
and $v' \in V_y$ with $x \neq y \in V$, is at most $\beta\gamma
d_G(x,y)+ 2\lambda \max_z \diam(H_z) \leq \gamma(\beta+2\lambda)
d_G(x,y)$. It follows that for $u \in H_x'$ and $v \in H_y'$,
\begin{eqnarray*}
d_X(u,v)=\beta\gamma d_{G}(x,y) &\le& d_{Y}(\phi(u),\phi(v)) \\
& \le&
\gamma (\beta+2\lambda) d_G(x,y) =
\left(\frac{\beta+2\lambda}{\beta}\right) d_X(u,v).
\end{eqnarray*}

Hence if $\beta \ge \frac{2\lambda}{\lambda-1}$, we have,
$\mathrm{dist}(\phi)\le \max
\left\{\lambda,\frac{\beta+2\lambda}{\beta}\right\} = \lambda . $
\Endproof\vskip4pt

Recall that a Banach space $X$ is called super-reflexive if it
admits an equivalent uniformly convex norm. A finite-metric
characterization of such spaces was found by Bourgain
\cite{bourgaintrees}. Namely, $X$ is superreflexive if and only if
for every $\alpha>0$ there is an integer $h$ such that the
complete binary tree of depth $h$ doesn't $\alpha$-embed into $X$.
Let \text{{\rm TREE}} denote the set of metrics supported on
trees. Since any weighted tree is almost isometric to a subset of
a deep enough complete binary tree, we conclude using
Lemma~\ref{lem:p2+}.

\begin{corollary} Let $X$ be a Banach space. Then the following  
assertions
are equivalent\/{\rm :}\/

\medskip
\noindent{\rm a)} $X$ is super-reflexive.

\medskip
\noindent{\rm b)} For any $\alpha>1$ there exists $\delta<1$ such
that for infinitely many integers $n$\/{\rm ,}\/
$$
R_{X}(\text{{\rm TREE}}; \alpha,n)\leq n^\delta.
$$

\medskip
\noindent{\rm c)} For any $\alpha>1$ there exists an integer $n$
such that
$$
R_{X}(\text{{\rm TREE}}; \alpha,n) < n.
$$
\end{corollary}

\section{Metric Ramsey-type theorems}\label{section:lower bounds}

In this section we prove Theorem~\ref{thm:lowerlarge}; i.e., we
give an $n^{\Omega(1)}$
lower bound on $R_2(\alpha,n)$ for $\alpha>2$.

The proof actually establishes a lower bound on
$R_{\um}(\alpha,n)$. The bound on $R_2$ follows since ultrametrics
embed isometrically in $\ell_2$. The lower bound for embedding
into ultrametrics utilizes their representation as hierarchically
well-separated trees. We begin with some preliminary background on
ultrametrics and hierarchically well-separated trees in
Section~\ref{section:lower-preliminaries}. We also note that our
proof of the lower bound makes substantial use of the notions of
metric composition and composition closure which were introduced
in Section~\ref{section:composition}.

We begin with a description of the lemmas on which the proof of
the lower bound is based and the way they are put together to
prove the main theorem. This is done in
Section~\ref{section:lower-main}. Detailed proofs of the main
lemmas appear in
Sections~\ref{section:lower-weighted-comp}--\ref{section:lower-k-hst-comp}.
Most of the proof is devoted to the case where $\alpha$ is a
fixed, large enough constant. In
Section~\ref{section:lower-small-alpha}, we extend the proof to
apply for every $\alpha>2$.

\Subsec{Ultrametrics and hierarchically well-separated  
trees}\label{section:lower-preliminaries}
 Recall that an {\em ultrametric} is a metric space $(X,d)$ such that  
for every $x,y,z\in X$,
$$
d(x,z)\le \max\{d(x,y),d(y,z)\}.
$$

A more restricted class of metrics with an inherently hierarchical
structure plays a key role in the sequel. Such spaces have already
figured prominently in earlier work on embedding into ultrametric
spaces \cite{bartal1}, \cite{bbm}.

\begin{definition}[\cite{bartal1}]\label{def:hst}
For $k\geq 1$, a $k$-\emph{hierarchically well-separated tree}\break
($k$-HST) is a metric space whose elements are the leaves of a
rooted tree $T$. To each vertex $u\in T$ there is associated a
label $\Delta(u) \ge 0$ such that $\Delta(u)=0$ if and only if $u$ is a leaf
of $T$. It is required that if a vertex $u$ is a child of a vertex
$v$ then $\Delta(u)\leq \Delta(v)/k$ . The distance between two
leaves $x,y\in T$ is defined as $\Delta(\lca(x,y))$, where
$\lca(x,y)$ is the least common ancestor of $x$ and $y$ in $T$.

A $k$-HST is said to be \emph{exact} if $\Delta(u)=\Delta(v)/k$
for every two internal vertices where $u$ is a child of $v$.
\end{definition}

First, note that an ultrametric on a finite set and a (finite)
$1$-HST are identical concepts. Any $k$-HST is also a $1$-HST,
i.e., an ultrametric. However, when $k>1$, a $k$-HST is a stronger
notion which has a hierarchically clustered structure. More
precisely, a $k$-HST with diameter $D$ decomposes into subspaces
of diameter at most $D/k$ and any two points at distinct subspaces
are at distance exactly $D$. Recursively, each subspace is itself
a $k$-HST. It is this hierarchical decomposition that makes
$k$-HST's useful.

When we discuss $k$-HST's, we freely use the tree $T$ as in
Definition~\ref{def:hst}, \emph{the tree defining the} HST. An
internal vertex in $T$ with out-degree $1$ is said to be
\emph{degenerate}. If $u$ is nondegenerate, then $\Delta(u)$ is
the diameter of the sub-space induced on the subtree rooted by
$u$. Degenerate nodes do not influence the metric on $T$'s leaves;
hence we may assume that all internal nodes are nondegenerate
(note that this assumption need not hold for
\emph{exact} $k$-HST's). 

We need some more notation:
\begin{notation}
Let $\um$   denote the class of ultrametrics, and $k$-HST denote 
the class of $k$-HST's. Also let $\eq$ denote the class of
equilateral spaces.
\end{notation}

The following simple observation is not required for the proof,
but may help direct the reader's intuition. More complex
connections between these concepts do play an important role in
the proof.

\begin{proposition}\label{prop:k-hst=comp-eq}
The class of $k$\/{\rm -HST'}\/s is the $k$-composition closure of the class
of equilateral spaces{\rm ;} i.e.{\rm ,} $k\text{-\hst} = \comp_k(\eq)$.

In particular{\rm ,} the class of ultrametrics is the $1$\/{\rm -}\/composition
closure of the class of equilateral spaces\/{\rm ;} i.e.{\rm ,} $\um =
\comp_1(\eq)$.
\end{proposition}

We recall the following well known fact (e.g.~\cite{lemin}), that allows
us to reduce the Euclidean Ramsey problem to the problem of
embedding into ultrametrics:


\begin{proposition}\label{prop:hst-l2} Any ultrametric is isometrically
embeddable in $\ell_2$. In particular{\rm ,}
$$ R_2(\alpha,n) \geq R_{\um}(\alpha,n). $$
\end{proposition}

This proposition can be proved by induction on the structure of
the tree defining the ultrametric. It is shown inductively that
each rooted subtree embeds isometrically into a sphere with radius
proportional to the subtree's diameter, and that any two subtrees
rooted at an internal vertex are mapped into orthogonal subspaces.

When considering Lipschitz embeddings, the $k$-HST representation
of an ultrametric comes naturally into play. This is expressed by
the following variant on a proposition from \cite{bartal2}:

\begin{lemma}\label{lem:exact} For any $k>1${\rm ,} any ultrametric is  
$k$\/{\rm -}\/equivalent to an exact $k${\rm -HST}.
\end{lemma}

Lemma~\ref{lem:exact} is proved via a simple transformation of the
tree defining the ultrametric. This is done by coalescing
consecutive internal vertices, whose labels differ by a factor
which is less than $k$. The complete proof of
Lemma~\ref{lem:exact} appears in
Section~\ref{section:lower-um-k-hst}

We end this section with a proposition on embeddings into
ultrametrics, which is implicit in \cite{bartal1}. Although this
proposition is {\em not} used in the proofs, it is useful for
obtaining efficient algorithms from these theorems.

\begin{lemma}\label{lem:dist n} Every $n$\/{\rm -}\/point metric space is
$n$-equivalent to an ultrametric.
\end{lemma}

\Proof 
Let $M$ be an $n$-point metric space. We inductively construct an
$n$-point HST $X$ with $\diam(X)=\diam(M)$ and a noncontracting
$n$-Lipschitz bijection between $M$ and $X$.

Define a graph with vertex set $M$ in which $[u,v]$ is an edge if
and only if $d_M(u,v)< \frac{\diam(M)}{n}.$ Clearly, this graph is
disconnected. Let $A_1,\dots , A_m$ be the vertex sets of the
connected components. By induction there are HST's $X_1,\dots ,X_m$
with $\diam(X_i)=\diam((A_i,d_M)) < \diam(M)$ and bijections
$f_i:A_i\to X_i$ such that for every $u,v\in A_i$, $d_M(u,v)\le
d_{X_i}(f_i(u),f_i(v))\le |A_i|d_M(u,v) < n d_M(u,v)$. Let $T_i$
be the tree defining $X_i$. We now construct the HST $X$ whose
defining labelled tree $T$ is rooted at $z$. The root's label is
$\Delta(z) = \diam(M)$ and it has $m$ children, where the $i$th
child, $u_i$, is a root of a labelled tree isomorphic to $T_i$.
Since $\Delta(u_i) = \diam(X_i) < \diam(M) = \diam(X) =
\Delta(z)$, the resulting tree $T$ indeed defines an HST. Finally,
if $u \in A_i$ and $v \in A_j$ for $i \neq j$ then $d_M(u,v) \ge
\diam(M)/n$. Since $\diam(X) = \Delta(z) = \diam(M)$, the
inductive hypothesis implies the existence of the required
bijection.
\hfill\qed

\Subsec{An overview of the proof of  
Theorem~{\rm \ref{thm:lowerlarge}}}\label{section:lower-main}
In this section we describe the proof of the following theorem:

\begin{theorem}\label{thm:lower-ultrametric}
There exists an absolute constant $C>0$ such that for every
$\alpha > 2${\rm ,}
$$
R_\um(\alpha,n) \ge n^{1 - C\frac{\log \alpha}{\alpha}}.
$$
\end{theorem}

By Proposition ~\ref{prop:hst-l2}, the same bound holds true for
$R_2(\alpha,n)$.

We begin with an informal description and motivation. The main
lemmas needed for the proof are stated, and it is shown how they
imply the theorem. Detailed proofs for most of these lemmas appear
in subsequent subsections.

Our goal is to show that for any $\alpha>2$, every $n$ point
metric space $X$ contains a subspace which is $\alpha$-equivalent
to an ultrametric of cardinality $\ge n^{\xx(\alpha)}$, where
$\XX(\alpha)$ is independent of $n$. In much of the proof we
pursue an even more illusive goal. We seek large subsets that
embed even into $k$-HST's (recall that this is a restricted class
of ultrametrics). A conceptual advantage of this is that it
directs us towards seeking \emph{hierarchical} substructures
within the given metric space. Such structures can be described as
the composition closure of some class of metric spaces $\M$. A
metric space in $\comp_\beta(\M)$ is composed of a hierarchy of
dilated copies of metric spaces from $\M$, and the proof
iteratively finds such large structures. The class $\M$ varies
from iteration to iteration, gradually becoming more restricted,
and getting closer to the class $\eq$. When $\M$ is approximately
$\eq$ this procedure amounts to finding a $k$-HST (due to
Proposition~\ref{prop:k-hst=comp-eq}). It is therefore worthwhile
to consider a special case of the general problem, where $X \in
\comp_\beta(\M)$, and we seek a subspace of $X$ that is
$\alpha$-equivalent to a $k$-HST.

It stands to reason that if spaces in $\M$ have large Ramsey
numbers, then something similar should hold true also for spaces
in $\comp_\beta(\M)$. After all, if $\beta$ is large, then the
copies of dilated metric spaces from $\M$ are hierarchically
well-separated. This would have reduced the problem of estimating
Ramsey numbers for spaces in $\comp_\beta(\M)$ to the same problem
for the  class $\M$.

While this argument is not quite true, a slight modification of it
does indeed work. For the purpose of this intuitive discussion, it
is convenient to think of $\beta$ as large, in particular with
respect to $k$ and $\alpha$. Consider that $X$ is the
$\beta$-composition of $M \in \M$ and a set of $|M|$ disjoint
metric spaces $\{ N_i \}_{i \in M}$, $N_i \in \comp_\beta(\M)$.
Assume (inductively) that each $N_i$ contains a subspace $N_i'$
that is $\alpha$-equivalent to a $k$-HST $H_i$ of size
$|N_i|^\xx$. Find a subspace $M'$ of $M$ that is also
$\alpha$-equivalent to a $k$-HST $K$ and attach the roots of the
appropriate $H_i$'s to the corresponding leaves of $K$ (with an
appropriate dilation). This yields a $k$-HST $H$, and by the
separation property of compositions with large $\beta$, we obtain
a subspace $X'$ of $X$ which is $\alpha$-equivalent to $H$.
However, the size of the final subspace $X' = \dot{\boldsymbol{\cup}}_{i \in M'}
N_i'$ depends not only on the size of $M'$, the subspace we find
in $M$, but also on how large the chosen $N_i'$s are. Therefore,
the correct requirement is that $M'$ satisfies:
$$
  \sum_{i \in M'} |N_i|^\xx \ge \left(\sum_{i \in M} |N_i|\right)^\xx .
$$
This gives rise to the following definition:

\begin{definition}[The weighted Ramsey function]\label{def:XX}
Let $\M,\N$ be classes of metric spaces. Denote by $\XX_\M(\N,
\alpha)$ the largest $0 \le \xx \le 1$ such that for every metric
space $X \in \N$ and any weight function $w:X\to \R^{+}$, there is
a subspace $Y$ of $X$ that $\alpha$-embeds in $\M$ and satisfies:
\begin{equation}
\sum_{x\in Y}w(x)^\xx\ge \left(\sum_{x\in X}w(x)\right)^\xx .
\tag{$*$} \label{eq:weighted}
\end{equation}
When $\N$ is the class of all metric spaces, it is omitted from
the notation.
\end{definition}

In what follows the notion of a weighted metric space refers to a
pair $(X,w)$, where $X$ is a metric space and $w:X\to \R^{+}$ is a
weight function.

The following is an immediate consequence of
Definition~\ref{def:XX} (by using the constant weight
function $w(x) \equiv 1$).
\begin{proposition}\label{prop:chi-R}
$$ R_\M(\N;\alpha,n) \ge n^{\XX_\M(\N,\alpha)} . $$
In particular{\rm ,}
$$ R_\M(\alpha,n) \ge n^{\XX_\M(\alpha)} . $$
\end{proposition}
\vskip8pt

We note that it is possible to show, via the results of
Section~\ref{section:composition}, that in the setting of
embedding into composition classes, and in particular in our case
of $k$-HST's or ultrametrics, the last inequality in
Proposition~\ref{prop:chi-R} holds with equality for infinitely
many $n$'s.

The entire proof is thus dedicated to bounding the weighted Ramsey
function when the target metric class is the class of
ultrametrics. The proofs in the sequel produce embeddings into
$k$-HST's and ultrametrics. In this context, the following
conventions for $\XX_\M(\N,\alpha)$ are useful:
\begin{itemize}
\item $\XX_k(\N,\alpha) = \XX_{k\text{-\hst}}(\N,\alpha)$. In
particular, $\XX_k(\alpha) = \XX_{k\text{-\hst}}(\alpha)$.

\item $\XX(\N,\alpha) = \XX_1(\N,\alpha) = \XX_\um(\N,\alpha)$. In
particular, $\XX(\alpha) = \XX_{\um}(\alpha)$.
\end{itemize}

The following strengthening of Theorem~\ref{thm:lower-ultrametric}
is the main result proved in this section.

\demo{\scshape Theorem $3.7'$}
{\it There exists an absolute constant $C>0$ such that for every}
$\alpha > 2$,
$$
\XX(\alpha) \ge 1 - C\frac{\log \alpha}{\alpha} .
$$
\vglue12pt

Our goal can now be rephrased as follows:
  given an arbitrary weighted metric space $(X,w)$,
find a subspace of $X$, satisfying the weighted Ramsey condition
\eqref{eq:weighted} with $\XX(\alpha)$ as in Theorem~$3.7'$, that
is $\alpha$-equivalent to an ultrametric.

  Before continuing with the outline of the proof,
we state a useful property of the weighted Ramsey function. When
working with the regular Ramsey question it is natural to perform
a procedure of the following form: first find a subspace which is
$\alpha_1$-embedded in some ``nice" class of metric spaces, then
find a smaller subspace of this subspace which is $\alpha_2$
equivalent to our target class of metric spaces, thus obtaining
overall $\alpha_1\alpha_2$ distortion. If the first subspace has
size $n' \ge n^{\xx_1}$ and the second is of size $n'' \ge
n'{}^{\xx_2}$ then $n'' \ge n^{\xx_1\xx_2}$.

The weighted Ramsey problem has the same super-multiplicativity
property:

\begin{lemma}\label{lem:xx-product1}
Let $\M,\N,\PP$ be classes of metric spaces and $\alpha_1,\alpha_2
\ge 1$. Then
$$ \XX_\M(\PP,\alpha_1 \alpha_2) \ge \XX_\M(\N,\alpha_1) \cdot  
\XX_\N(\PP,\alpha_2)
. $$
\end{lemma}
\vskip8pt

The interpretation of this lemma (proved in
\S\ref{section:lower-weighted-comp}) is as follows: Suppose
that we are given a metric space in $\PP$ and we seek a subspace
that embeds with low distortion in $\M$, and satisfies
condition~\eqref{eq:weighted}. We can first find a subspace which
$\alpha_1$-embeds in $\N$ and then a subspace which  
$\alpha_2$-embeds in~$\M$. In the course of this procedure we
multiply the distortions and the $\xx$'s of the corresponding
classes.

The discussion in the paragraph preceding Definition~\ref{def:XX}
on how Ramsey-type properties of class $\M$ carry over to
$\comp_\beta(\M)$, leads to the following proposition: If for
every $X \in \M$ and every $w:X\to \R^+$ there is a subspace
$Y\subset X$, satisfying the weighted Ramsey
condition~\eqref{eq:weighted} with parameter $\xx$, which is
$\alpha$-equivalent to a $k$-HST, then the same holds true for
every $M \in \comp_\beta(\M)$. In our notation, we have the
following lemma (proved in
\S\ref{section:lower-weighted-comp}):

\begin{lemma}\label{lem:composition-xx1}
Let $\M$ be a class of metric spaces. Let $k \ge 1$ and $\alpha
\ge 1$. Then for any $\beta \ge \alpha k${\rm ,}
$$ \XX_k(\comp_\beta(\M),\alpha) = \XX_k(\M,\alpha) .
$$
In particular for $\beta \ge \alpha${\rm ,}
$$ \XX(\comp_{\beta}(\M),\alpha) = \XX(\M,\alpha) . $$
\end{lemma}
\vskip8pt

The following simple notion is used extensively in the sequel.

\begin{definition}\label{def:phi}
The \emph{aspect ratio} of a finite metric space $M$, is defined
as:
$$
\Phi(M)=\frac{\diam(M)}{\min_{x\neq y} d_M(x,y)}.
$$
When $|M|=1$ we use the convention $\Phi(M)=1$. We note that
$\Phi(M)$ can be viewed as $M$'s normalized diameter, or as its
Lipschitz distance from an equilateral metric space.
\end{definition}

Again, it is helpful to consider the $k$-HST representation of an
ultrametric~$Y$. In particular, notice that in this hierarchical
representation, the number of levels is $O(\log_k \Phi(Y))$. In
view of this fact, it seems reasonable to expect that when
$\Phi(X)$ is small it would be easier to find a large subspace of
$X$ that is close to an ultrametric. This is, indeed, shown in
Section~\ref{section:lower-aspect-ratio}.

\begin{definition}
The class of all metric spaces $M$ with aspect ratio $\Phi(M)\break \le
\Phi$, for some given parameter $\Phi$, is denoted by $\N(\Phi)$.
Two more conventions that we use are: For every real $\Phi \ge 1$,
\begin{itemize}
\item $\XX(\Phi,\alpha) = \XX(\N(\Phi),\alpha)$. Similarly
$\XX_k(\Phi,\alpha) = \XX_k(\N(\Phi),\alpha)$, and in general
where $\M$ is a class of metric spaces, $\XX_\M(\Phi,\alpha) =
\XX_\M(\N(\Phi),\alpha)$.

\item $\comp_\beta(\Phi) = \comp_\beta(\N(\Phi))$.
\end{itemize}
\end{definition}

The main idea in bounding $\XX(\Phi,\alpha)$ is that the metric
space can be decomposed into a small number of subspaces, the
number of which can be bounded by a function of $\Phi$, such that
we can find among these, subspaces that are far enough from each
other and contain enough weight to satisfy the weighted Ramsey
condition \eqref{eq:weighted}. Such a decomposition of the space
yields the recursive construction of a hierarchically
well-separated tree, or an ultrametric. This is done in the proof
of the following lemma. A more detailed description of the ideas
involved in this decomposition and the proof of the lemma can be
found in Section~\ref{section:lower-aspect-ratio}.

\begin{lemma}\label{lem:Phi1}
There exists an absolute constant $C'>0$ such that for every
$\alpha > 2$ and $\Phi\ge 1${\rm :}
$$
\XX(\Phi,\alpha)\ge 1-C'\frac{\log\alpha +
\log\log(4\Phi)}{\alpha}.
$$
\end{lemma}
\vskip8pt

Note that for the class of metric spaces with aspect ratio $\Phi
\le \exp(O(\alpha))$, Lemma~\ref{lem:Phi1} yields the bound stated
in Theorem {$3.7'$}.

Combining Lemma~\ref{lem:Phi1} with
Lemma~\ref{lem:composition-xx1} gives an immediate consequence on
$\beta$-composition classes: for $\beta \ge \alpha$, 

\begin{equation}\label{eq:comp-Phi}
\XX(\comp_\beta(\Phi),\alpha) = \XX(\Phi,\alpha) \ge
1-C'\frac{\log\alpha + \log\log(4\Phi)}{\alpha}.
\end{equation}
\vskip8pt

We now pass to a more detailed description of the proof of
Theorem~$3.7'$. Let $X$ be a metric space and assume that for some
specific value of $\alpha$ we can prove the bound in the theorem
(e.g., this trivially holds for $\alpha=\Phi(X)$ where we have
$\XX(X,\alpha)=1$).

Let $\hat{X}$ be an arbitrary metric space and let $X$ be a
subspace of $\hat{X}$ that is\break $\alpha$-equivalent to an
ultrametric, satisfying the weighted Ramsey
condition~\eqref{eq:weighted} with $\xx = \XX(\hat{X},\alpha)$. We
will apply the following ``distortion refinement" procedure:
  find a subspace of $X$ that is $(\alpha/2)$-equivalent to an
ultrametric, satisfying condition~\eqref{eq:weighted} with $\xx'
\ge (1-C''\frac{\log\alpha}{\alpha})$. This implies that
$\XX(\hat{X},\alpha/2) \ge
(1-C''\frac{\log\alpha}{\alpha})\XX(\hat{X},\alpha)$.
{Theorem~$3.7'$} now follows: we start with $\alpha=\Phi(\hat{X})$
and then apply the above distortion refinement procedure
iteratively until we reach a distortion below our target. It is
easy to verify that this implies the bound stated in the theorem.

The distortion refinement uses the bound in \eqref{eq:comp-Phi} on
$\XX(\comp_\beta(\Phi),\alpha'')$, in the particular case
$\alpha'' < \alpha/2$ and $\Phi \le \exp(O(\alpha))$. This is
useful due the following claim: if $X$ is $\alpha$-equivalent to
an ultrametric then it contains a subspace $X'$ which is
$(1+2/\beta)$-equivalent to a metric space $Z$ in
$\comp_\beta(\Phi)$, for $\Phi \le \exp(O(\alpha))$, and which
satisfies condition~\eqref{eq:weighted} with $\xx'' \ge
(1-2\frac{\log\alpha}{\alpha})\xx$. By \eqref{eq:comp-Phi} we
obtain a subspace $Z'$ of $Z$ which is $\alpha''$-equivalent to an
ultrametric. By appropriately choosing all the parameters, we see from Lemma~\ref{lem:xx-product1} that
there is a subspace
$X''$ of $X'$ which is $(\alpha/2)$-equivalent to an ultrametric,
and the desired bound on $\XX(\hat{X},\alpha/2)$ is achieved.

  The proof of the above claim is based on two lemmas relating  
ultrametrics,
$k$-HST's and metric compositions. Let $X$ be $\alpha$-equivalent
to an ultrametric $Y$. The subspace $X'$ is produced via a
Ramsey-type result for ultrametrics which states that every
ultrametric $Y$ contains a subspace $Y'$ which is
$\alpha'$-equivalent to a $k$-HST with $k > \alpha'$.
(Lemma~\ref{lem:exact} can be viewed as a non-Ramsey result of
this type when $k=\alpha'$.) Moreover, we can ensure that
condition~\eqref{eq:weighted} is satisfied for the pair $Y'\subset
Y$ with the bound stated below.

\begin{lemma}\label{lem:hst-sparse1}
For every $k \ge \alpha > 1${\rm ,}
$$ \XX_k(\um,\alpha) \ge 1 - \frac{\log (k/\alpha)}{\log \alpha} . $$
\end{lemma}
\vskip8pt

The proof of this lemma involves an argument on general tree
structures described in Section~\ref{section:lower-um-k-hst}.

Now, by Lemma~\ref{lem:xx-product1} we obtain a subspace $X'$ that
is $\alpha'\alpha$-equivalent to a $k$-HST for $k
> \alpha'$. If $k$ is large enough then the subtrees of the
$k$-HST impose a clustering of $X'$. That is, each subtree
corresponds to a subspace of $X'$ of very small diameter, whereas
the $\alpha$ distortion implies that the aspect ratio of the
metric reflected by inter-cluster distances is bounded by
$\alpha$. By a recursive application of this procedure we obtain a
metric space in $\comp_\beta(\alpha)$, with the exact relation
between $k$,$\alpha$, and $\beta$ stated in the lemma below. The
details of this construction are given in
Section~\ref{section:lower-k-hst-comp}.

\begin{lemma}\label{lem:k-hst-comp1}
For any $\alpha,\beta \ge 1${\rm ,} if a metric space $M$ is
$\alpha$-equivalent to an $\alpha\beta${\rm -HST} then $M$ is
$(1+2/\beta)$\/{\rm -}\/equivalent to a metric space in
$\comp_\beta(\alpha)$.
\end{lemma}

The distortion refinement process described above is formally
stated in the following lemma:

\begin{lemma}\label{lem:refinement1}
There exists an absolute constant $C''>0$ such that for every
metric space $\hat{X}$ and any $\alpha > 8${\rm ,}
$$ \XX\left(\hat{X},\frac{\alpha}{2}\right) \ge \XX(\hat{X},\alpha)  
\left( 1 - C''\frac{\log\alpha}{\alpha} \right) . $$
\end{lemma}

\vskip8pt

\Proof 
Fix a weight function $w:\hat{X}\to \R^+$, let $X$ be a subspace
of $\hat{X}$ that is $\alpha$ equivalent to an ultrametric $Y$,
and satisfies the weighted Ramsey condition~\eqref{eq:weighted}
with $\psi(\hat{X},\alpha)$.  Fix two numbers $\alpha',\beta \ge
1$ which will be determined later, and set $k=\alpha\alpha'\beta$.
Lemma~\ref{lem:hst-sparse1} implies that $Y$ contains a subspace
$Y'$ which is $\alpha'$-equivalent to a $k$-HST, and $Y'$
satisfies condition~\eqref{eq:weighted} with
$\XX_k(\um,\alpha')\ge 1-\frac{\log(k/\alpha')}{\log\alpha'}$. By
mapping $X$ into an ultrametric $Y$, and then mapping the image of
$X$ in $Y$ into a $k$-HST, we apply Lemma~\ref{lem:xx-product1},
obtaining a subspace $X'$ of $X$ that is
$\alpha'\alpha$-equivalent to a $k$-HST $W$, which satisfies
condition~\eqref{eq:weighted} with exponent $ \XX_k(\um,\alpha')
\cdot \XX(\hat{X},\alpha) \ge \left( 1- \frac{\log
(k/\alpha')}{\log \alpha'} \right) \XX(\hat{X},\alpha)$. Denote $\Phi =
\alpha'\alpha$. We have that $X'$ is $\Phi$-equivalent to a
$\Phi\beta$-HST and therefore by Lemma~3.16, $X'$
is $(1+2/\beta)$ equivalent to a metric space $Z$ in
$\comp_\beta(\Phi)$.
Now, we can use the bound in \eqref{eq:comp-Phi} to find a
subspace $Z'$ of $Z$ that is $\beta$-equivalent to an ultrametric,
and satisfies condition~\eqref{eq:weighted} with exponent
$\XX(\comp_\beta(\Phi),\beta)$. By mapping $X'$ into $Z \in
\comp_\beta(\Phi)$ and finally to an ultrametric, we apply
Lemma~\ref{lem:xx-product1} again, obtaining a subspace $X''$ of
$\hat{X}$ that is $\beta(1+2/\beta) =\beta + 2$ equivalent to an
ultrametric $U$, satisfying condition~\eqref{eq:weighted} with
exponent
\begin{eqnarray*}
\left(1-C'\frac{\log \beta + \log\log(4\Phi)}{\beta} \right)
\left(1-\frac{\log (\alpha\beta)}{\log \alpha'} \right)
\XX(\hat{X},\alpha).
\end{eqnarray*}

Finally, if we choose $\beta = \alpha/2-2$ and let $\Phi =
2^{2\alpha}$ (which determines $\alpha'$), we get that
$$
\XX\left(\hat{X},\frac{\alpha}{2}\right) \ge \left(1-9C'\frac{\log
\alpha}{\alpha} \right) \left(1-2\frac{\log \alpha}{\alpha}
\right) \XX(\hat{X},\alpha) \ge \left(1-C''\frac{\log
\alpha}{\alpha} \right) \XX(\hat{X},\alpha),
$$
for an appropriate choice of $C''$
\Endproof\vskip4pt

{Theorem~$3.7'$} is a straightforward consequence of
Lemma~\ref{lem:refinement1}:

\demo{Proof of Theorem~$3.7'$} By an appropriate choice of $C$  
we may assume that $\alpha >
8$. Let $X$ be a metric space and set $\Phi=\Phi(X)$. Recall that
$\XX(X,\Phi) = 1$. Let $m = \lfloor \log\alpha \rfloor$ and $M =
\lceil \log \Phi \rceil$. Lemma~\ref{lem:refinement1} implies that
$\XX(X,\alpha/2) \ge \XX(X,\alpha) -
C''\frac{\log\alpha}{\alpha}$, and so by an iterative application
of this lemma we get
\begin{eqnarray*}
\XX(X,\alpha)& \ge& \XX(X,2^m) \ge \XX(X,2^M) - C''\sum_{i=m+1}^{M}
\frac{i}{2^i} \\
&\ge& 1 - C''\sum_{i=m+1}^{\infty} \frac{i}{2^i} = 1 -
C''\frac{m+2}{2^m} \ge 1 - 6C''\frac{\log\alpha}{\alpha} .
\end{eqnarray*}
\vglue-28pt
\Endproof\vskip12pt

This completes the overview of the proof of Theorem~$3.7'$.
Sections~\ref{section:lower-weighted-comp}--\ref{section:lower-k-hst-comp}
contain the proofs of the lemmas described above.

  Additionally, in Section~\ref{section:lower-small-alpha} we
describe in detail how to achieve Ramsey-type theorems for
arbitrary values of $\alpha>2$. The main ideas that make this
possible are first, replacing Lemma~\ref{lem:Phi1} with another
lemma that can handle distortions $2+\epsilon$ and second,
providing a more delicate application of our Lemmas, using the
fact that we can find $k$-HST's with large $k$ ($\approx
1/\epsilon$) rather than just ultrametrics, to ensure that
accumulated losses in the distortion are small.

We end with a discussion on the algorithmic aspects of the metric
Ramsey problem. Given a metric space $X$ on $n$ points, it is
natural to ask wether we can find in polynomial time a subspace
$Y$ of $X$ with $n^\xx$ points which is $\alpha$-equivalent to an
ultrametric, for $\xx$ as in Theorem~\ref{thm:lower-ultrametric}.
It is easily checked that the proofs of our lemmas yield
polynomial time algorithms to find the corresponding subspaces.
Thus, the only obstacle in achieving a polynomial time algorithm,
is the fact, that the proof of {Theorem~$3.7'$} involves $O(\log
\Phi)$ iterations of an application of
Lemma~\ref{lem:refinement1}. We seek, however, a polynomial
dependence only on $n$. This is remedied as follows: It is easily
seen that using Lemma~\ref{lem:dist n} we can start from
$\XX(X,|X|) = 1$ rather than $\XX(X,\Phi(X))=1$. Thus we replace
the bound of $\Phi$ with $n$, and end up with at most $O(\log n)$
iterations of Lemma~\ref{lem:refinement1}. This implies a
polynomial time algorithm to solve the metric Ramsey problem.

\Subsec{The weighted metric Ramsey problem and its relation to  
metric composition}\label{section:lower-weighted-comp}
In this section we prove Lemmas~\ref{lem:composition-xx1} and
 \ref{lem:xx-product1}. We begin with
Lemma~\ref{lem:xx-product1}, which allows us to move between
different classes of metric spaces while working with the weighted
Ramsey problem.

\demo{\scshape Lemma 3.10}
{\it Let $\M,\N,\PP$ be classes of metric spaces and $\alpha_1,\alpha_2
\ge 1$. Then}
$$ \XX_\M(\PP,\alpha_1 \alpha_2) \ge \XX_\M(\N,\alpha_1) \cdot  
\XX_\N(\PP,\alpha_2)
. $$
 \vskip8pt

\Proof 
Let $\xx_1 = \XX_\M(\N,\alpha_1)$ and $\xx_2 =
\XX_\N(\PP,\alpha_2)$. Take $P \in \PP$ and a weight function
$w:P\to \R^+$. There are a subspace $P'$ of $P$ and an
$\alpha_2$-embedding $f:P'\to N$, where $N \in \N$, and
$$
\sum_{x\in P'}w(x)^{\xx_2}\ge \left(\sum_{x \in P}
w(x)\right)^{\xx_2}.
$$

  Similarly, for every weight function $w':N\to\R^+$ there exists a  
subspace $N'$ of $N$ and an
  $\alpha_1$-embedding $g:N'\to M$, where
$M \in \M$, and
$$
\sum_{y\in N'}w'(y)^{\xx_1}\ge \left(\sum_{y \in N}
w'(y)\right)^{\xx_1} .
$$
By letting $P'' = f^{-1}(N')$, and for $y \in N$, $w'(y) =
w(f^{-1}(y))^{\xx_2}$, we get that
$$
\sum_{x\in P''}w(x)^{\xx_1\xx_2}\ge \left(\sum_{x \in P'}
w(x)^{\xx_2} \right)^{\xx_1}  \ge \left(\sum_{x \in P}
w(x)\right)^{\xx_1\xx_2}.
$$
Define $h:P'' \to M$ by $h(x) = g(f(x))$; then $h$ is an
$\alpha_1\alpha_2$-embedding.
\Endproof\vskip4pt

Lemma~\ref{lem:composition-xx1} shows that the weighted Ramsey
function stays unchanged as we pass from a class $\M$ of metric
spaces to its composition closure.  To repeat:

\demo{\scshape Lemma 3.11}
{\it Let $\M$ be a class of metric spaces. Let $k \ge 1$ and $\alpha
\ge 1$. Then for any} $\beta \ge \alpha k$,
$$ \XX_k(\comp_\beta(\M),\alpha) = \XX_k(\M,\alpha) . $$
 \vskip8pt

\Proof 
Since $\M \subseteq \comp_\beta(\M)$, clearly
$\XX_k(\comp_\beta(\M),\alpha) \le \XX_k(\M,\alpha)$. In what
follows we prove the reverse inequality.

Let $\xx\! = \!\XX_k(\M,\alpha)$. Let $X \!\in\! \comp_\beta(\M)$. We
prove that for any $w\!:\!X\!\to\! \R^{+}$ there exists a subspace $Y$ of
$X$ and a $k$-HST $H$ such that $Y$ is $\alpha$-equivalent to $H$
via a noncontractive $\alpha$-Lipschitz embedding, and:
$$
\sum_{x\in Y}w(x)^\xx\ge \left(\sum_{x \in X} w(x)\right)^\xx .
$$

The proof is by structural induction on the metric composition. If
$X \in \M$ then this holds by definition of $\xx$. Otherwise,
let $M \in \M$ and $\N = \{ N_z \}_{z \in M}$ be such that $X =
M_\beta[\N]$.

By induction, for each $z \in M$, there is a subspace $Y_z$ of
$N_z$ that is\break $\alpha$-equivalent to a $k$-HST $H_z$, defined by
the tree $T_z$, via a noncontractive $\alpha$-Lipschitz
embedding, and

$$
\sum_{u\in Y_z}w(u)^\xx\ge \left(\sum_{u \in N_z} w(u)\right)^\xx
.
$$

For a point $z \in M$ let $w'(z) = \sum_{u \in N_z} w(u)$. There
exists a subspace $Y_M$ of $M$ that is $\alpha$-equivalent to a
$k$-HST $H_M$, defined by $T_M$, via a noncontractive
$\alpha$-Lipschitz embedding, and
$$
\sum_{z\in Y_M}w'(z)^\xx\ge \left(\sum_{z \in M} w'(z)\right)^\xx
= \left(\sum_{x \in X} w(x)\right)^\xx.
$$

Let $Y = \boldsymbol{\cup}_{z \in Y_M} Y_z$. It follows that
$$
\sum_{u\in Y}w(x)^\xx = \sum_{z \in Y_M}  \sum_{u \in Y_z}
w(u)^\xx  \ge \sum_{z \in Y_M} \left(\sum_{u \in N_z}
w(u)\right)^\xx \ge \left(\sum_{x \in X} w(x)\right)^\xx .
$$

We now construct a $k$-HST $H$ that is defined by a tree $T$, as
follows. Start with a tree $T'$ that is isomorphic to $T_M$ and
has labels $\Delta(u) = \beta\gamma \cdot \Delta_{T_M}(u)$ (where
$\gamma = \frac{\max_{z \in M} \diam(N_z)}{\min_{x \neq y \in M}
d_M(x,y)}$, as in Definition~\ref{def:metric-composition}). At
each leaf of the tree\break\vglue-12pt\noindent corresponding to a point $z \in M$, create a
labelled subtree rooted at $z$ that is isomorphic to $T_z$ with
labels as in $T_z$. Denote the resulting tree by $T$. Since we
have a noncontractive $\alpha$-embedding of $Y_z$ in $H_z$, it
follows that $\Delta(z) = \diam(H_z) \leq \alpha \diam(Y_z) \le
\alpha \diam(N_z)$. Let $p$ be a parent of $z$ in $T_M$. Since we
have a noncontractive $\alpha$-embedding of $Y_M$ into $H_M$, it
follows that $\Delta_{T_M}(p) \ge d_M(x,y)$ for some $x,y \in
Y_M$. Therefore $\Delta(p) \ge \beta\gamma \cdot \min \{ d_M(x,y);
x \neq y \in M \}$. Consequently, $\Delta(p)/\Delta(z) \ge
\beta/\alpha \ge k$. Since $H_M$ and $H_z$ are $k$-HST's, it
follows that $T$ also defines a $k$-HST.

It is left to show that $Y$ is $\alpha$-equivalent to $H$. Recall
that for each $z \in M$ there is a noncontractive Lipschitz
bijection $f_z:Y_z\to H_z$ that satisfies for every $u,v \in Y_z$,
$d_{Y_z}(u,v) \le d_{H_z}(f_z(u),f_z(v)) \le \alpha d_{Y_z}(u,v)$.
Define $f:Y\to H$ as follows. If $z \in M$ and $u\in N_z$, set
$f(u) = f_z(u)$. Then for $u,v \in M$ such that $u,v \in N_z$ we
have
\begin{eqnarray*}
d_Y(u,v) &=& d_{Y_z}(u,v) \le d_{H_z}(f_z(u),f_z(v))\\
& =&
d_{H}(f(u),f(v)) \le \alpha d_{Y_z}(u,v) = \alpha d_Y(u,v) .
\end{eqnarray*}

Additionally, we have a noncontractive Lipschitz bijection
$f_M:Y_M\to H_M$ that satisfies for every $x,y \in Y_M$,
$d_{Y_M}(x,y) \le d_{H_M}(f_M(x),f_M(y)) \le \alpha d_{Y_M}(x,y)$.
Hence for $x \neq y \in M$, and $u \in N_x, y \in N_y$,
\begin{eqnarray*}
d_Y(u,v) &=& \beta\gamma d_{Y_M}(x,y) \le \beta\gamma
d_{H_M}(f_M(u),f_M(v))\\
& = &d_{H}(f(u),f(v)) \le \alpha\beta\gamma
d_{Y_M}(x,y) = \alpha d_Y(u,v) .
\end{eqnarray*}\vskip-32pt 
\phantom{ha}\hfill\qed

\Subsec{Exploiting metrics with bounded aspect  
ratio}\label{section:lower-aspect-ratio}
In this section we prove Lemma~\ref{lem:Phi1}
(\S\ref{section:lower-main}). That is, we give lower bounds
on $\XX=\XX(\Phi,\alpha)$ which depend on the aspect ratio of the
metric space, $\Phi$.

The proof of the lemma starts by obtaining lower bounds for a
restricted class of weight functions $w$. These bounds are then
extended to general weights. The class of ``nice" weight functions
is itself divided into two classes. In one class we have a lower
bound on the minimal weight relative to the total overall weight,
and the other is the constant weight function. This is formally
defined as follows:

\begin{definition} Fix some $q\ge 1$. A sequence
$x=\{x_i\}_{i=1}^\infty$ of nonnegative real numbers will be
called $q$-{\em decomposable} if there exists $\omega > 0$ such
that:
$$
\{i\in \mathbb{N};\ x_i>0\}=\left\{i\in \mathbb{N};\ x_i\ge
\frac{1}{q}\sum_{j=1}^\infty x_j\right\}\bigcup\{i\in \mathbb{N};\
x_i=\omega\}.
$$
\end{definition}

We will prove the following lemma:

\begin{lemma}\label{lem:qPhi} Let $q\ge 2${\rm ,} and $t \ge 8$ be an  
integer. Let
$(M,d)$ be an $n$\/{\rm -}\/point metric space and let $w:M\to \R^{+}${\rm ,} be a
weight function such that $\{w(x)\}_{x\in M}$ is $q$\/{\rm -}\/decomposable.
Then there exists a subspace $N\subseteq M$ that is
$4t$\/{\rm -}\/equivalent to an ultrametric and satisfies\/{\rm :}\/ \vglue-13pt
$$
\sum_{x\in N}w(x)^{\xx}\ge \left(\sum_{x\in M}w(x)\right)^{\xx},
$$
where $\xx = [t\log (4q\Phi(M))]^{-2/t}$.
\end{lemma}

The proof of Lemma~\ref{lem:qPhi} uses a decomposition of the
metric space $M$ into a small number of subspaces. This type of
strategy has been used in several earlier papers in
combinatorics and theoretical computer science, but the argument
closest in spirit to ours is in \cite{bbm}. The idea is to consider
two diametrical \pagebreak points, split
the space into shells according to the distance from one of these
two points, and discard the points in one of the shells.
Intuitively, we would like to discard a shell with small weight.
The exact choice is somewhat more sophisticated, tailored to
ensure the weighted Ramsey condition~\eqref{eq:weighted}. The
other shells form two subsets of the space that are substantially
separated. By an appropriate choice of the parameters, we can
guarantee that the union of the inner layers has diameter smaller
than a constant factor of the diameter of the whole space, and
hence a smaller aspect ratio. The role of $q$-decomposable weights
is as follows: This argument works fairly well for uniform
weights, and a slight modification of it yields bounds as a
function of $q$ (in addition to~$\Phi$) when in the weighted case
only a few points carry each at least $\frac{1}{q}$
of the total weight. Here the argument splits according to the
diameter of the set of "heavy" points. If the diameter is small,
the previous argument is started from a point that resides far
away from the heavy points. This guarantees that none of the
"heavyweights" get eliminated in the above-described process. If
their diameter is proportional to that of the whole space, it is
possible to argue similarly to the uniform-weight case, except
that we now obtain better bounds, since we can make estimates in
terms of $q$ (rather than the cardinality of the space $n$).

The extension of Lemma \ref{lem:qPhi} to arbitrary weight
functions requires a lemma on numerical sequences. This lemma
allows us to reduce the case of general sequences of weights to
$q$-decomposable ones.

\begin{lemma}\label{lem:lb} Fix $q\ge 16$ and let  
$x=\{x_i\}_{i=1}^\infty$ be a
sequence of nonnegative real numbers. Denote $p=
1-\frac{\log_2\log_2 q}{\log_2 q}$. There exists a sequence
$y=\{y_i\}_{i=1}^\infty$ such that $y_i \leq x_i$ for all $i \ge
1${\rm ,} $ \sum_{i\ge 1} y_i^p \ge \left(\sum_{i\ge 1} x_i\right)^p, $
  and the sequence $\{y_i^p\}_{i=1}^\infty$ is $q$\/{\rm -}\/decomposable.
\end{lemma}

Together these lemmas imply our main lemma:

\vskip6pt{\scshape Lemma 3.14.} {\it For every
$\alpha > 2$ and every} $\Phi\ge 1$:
$$
\XX(\Phi,\alpha)\ge 1-C\frac{\log \alpha+\log\log4\Phi}{\alpha},
$$
{\it where $C$ is a universal constant.}

\vglue4pt
{\it Proof}.
Clearly we may assume that $\alpha \ge 32$. Let $X$ be a metric
space with $\Phi(X) \le \Phi$, and $w:X\to \R^{+}$ a weight
function. Set $t = \lfloor \alpha/4 \rfloor$. By applying
Lemma~\ref{lem:lb} to the sequence $\{w(x)\}_{x\in X}$ with
$q=2^t$, we obtain a weight function $w'$ such that $w'(x)\leq
w(x)$ for all $x\in X$, the sequence $\{w'(x)^p\}_{x\in X}$ is
$q$-decomposable, and
$$
\sum _{x\in X}w'(x)^p \geq \left(\sum_{x\in X} w(x)\right)^p,
$$
where $p= 1-\frac{\log_2 t}{t}$.\smallbreak

Let $\beta = [t\log (4q\Phi(X))]^{-2/t}$ and apply
Lemma~\ref{lem:qPhi} to the space $X$ and weights $w'' = w'{}^p$. We
obtain a subspace $Y$ which is $4t$-equivalent to an ultrametric,
such that
$$
\sum_{x\in Y}w(x)^{p\beta}\ge \sum_{x\in Y} w'(x)^{p\beta} \ge
\left(\sum_{x\in X}w'(x)^p\right)^{\beta} \ge \left(\sum_{x\in
X}w(x)\right)^{p\beta}.
$$
Therefore,
\begin{eqnarray*}
\XX(\Phi,\alpha) & \ge & p\beta \ge \left(1-\frac{\log_2
t}{t}\right)\left(1-\frac{4\log t}{t}-\frac{2\log\log
(4\Phi)}{t}\right) \\
&\ge&  1-C\frac{\log \alpha
+\log\log(4\Phi)}{\alpha} ,
\end{eqnarray*}
for an appropriate choice of $C$.
\Endproof\vskip4pt

We now pass to the proof of Lemma~\ref{lem:lb}.
  Let $x=\{x_i\}_{i=1}^\infty$ be a sequence of nonnegative real numbers
which isn't identically zero. Let $p\ge 0$. Recall that the
$(p,\infty)$ norm of $x$ is defined by
$\|x\|_{p,\infty}=\sup_{i\ge 1} i^{1/p}x^*_i$, where
$\{x_i^*\}_{i=1}^\infty$ is the nonincreasing rearrangement of
the sequence $(|x_i|)_{i=1}^\infty$. We will require the following
numerical fact:

\begin{lemma}\label{lem:pinfty} For every $x\in \ell_1$ as above and  
every $0<p<1${\rm :}
$$
\|x\|_{p,\infty}\ge \left(\frac{1-p}{2-p}\right)^{1/p}\cdot
\frac{\|x\|_1^{1/p}}{\|x\|_\infty^{(1-p)/p}}.
$$
\end{lemma}
\vskip8pt

\Proof 
We can assume without loss of generality that $\|x\|_1=1$ and
$\|x\|_\infty=x_1\ge x_2\ge\dots \ge 0$. Obviously
$\|x\|_{p,\infty}\ge x_1$ so that if $x_1\ge
[(1-p)/(2-p)]^{1/p}x_1^{-(1-p)/p}$ we are done. Assume therefore
that the reverse inequality holds, i.e., $x_1 < \frac{1-p}{2-p}$.
Set $\alpha=\|x\|_{p,\infty}^p$ and denote
$j=\lceil \alpha/x_1^p \rceil+1$. Note that for every $i\ge 1$,
$x_i \leq \left(\alpha/i\right)^{1/p}$. Therefore,
\begin{eqnarray*}
\sum_{i=1}^{j-1} x_i & \leq & (j-1) x_1
      \leq  \left\lceil \frac{\alpha}{x_1^p} \right\rceil x_1  \le
   x_1\left(\frac{\alpha}{x_1^p} + 1\right)
      =  \alpha x_1^{1-p} +x_1,
\end{eqnarray*}
and
\begin{eqnarray*}
\sum_{i=j}^\infty x_i
   & \leq & \sum_{i=j}^\infty \frac{\alpha^{1/p}}{i^{1/p}}
    \leq  \alpha^{1/p} \int_{j-1}^\infty z^{-1/p}dz \\
   & \leq & \alpha^{1/p} \frac{p}{1-p} \cdot
   \left\lceil \frac{\alpha}{x_1^p} \right\rceil^{-\frac{1-p}{p}}
   \le   
\alpha^{1/p}\frac{p}{1-p}\left(\frac{x_1^p}{\alpha}\right)^{\frac{1- 
p}{p}} = \frac{p}{1-p} \alpha x_1^{1-p}.
\end{eqnarray*}
By summing both inequalities and using the bound on $x_1$ we get
\begin{eqnarray*}
  \frac{1}{1-p} \alpha x_1^{1-p} +\frac{1-p}{2-p} \ge 1,
\end{eqnarray*}
which simplifies to give the required result.
\Endproof\vskip4pt

\demo{Proof of Lemma~{\rm \ref{lem:lb}}}
We may assume that $x$ is a nonincreasing sequence of
nonnegative real numbers and that $\|x\|_1=1$.

We will prove below that there exist indexes $0\le l\le b$ such
that $x_l^p\ge \frac{2}{q}$ and:
\begin{equation}\label{lb-goal}
  S = \sum_{i=1}^l x_i^p+(b-l)x_b^p\ge
1.
\end{equation}

If $b=l$ assume that $l$ is the minimal index for which
(\ref{lb-goal}) holds. It follows that in this case $S=
\sum_{i=1}^{l} x_i^p < 1+x_l^p \leq 2$. Similarly if $b>l$, fix
$l$ and assume that $b$ is the minimal index for which
(\ref{lb-goal}) holds. It follows that $ S = \sum_{i=1}^l
x_i^p+(b-l)x_b^p \leq \sum_{i=1}^l x_i^p+(b-1-l)x_{b-1}^p +x_b^p <
1+x_b^p \leq 2$.\vglue3pt

Define the sequence $\{y_i\}_{i=1}^{\infty}$ so that $y_i=x_i$ for
$i\le l$, $y_i=x_b$ for $l<i\le b$ and $y_i=0$ for $i>b$. It
follows that $y_i \leq x_i$ for all $i\ge 1$ and $\sum_{i\ge 1}
y_i^p = S \ge 1 = (\sum_{i\ge 1} x_i)^p$. Since for $j\le l$
$y_j^p = x_j^p \ge \frac{2}{q} \ge \frac{S}{q} = \frac{1}{q}
\sum_{i\ge 1} y_i^p$, for $l<j\le b$, $y_i^p = x_b^p$ and for
$j>b$, $y_i^p=0$, we get that $\{y_i^p\}_{i=1}^{\infty}$ is
$q$-decomposable.

It remains to prove (\ref{lb-goal}).
  Let $l\ge 0$ be the largest
integer for which $x_l^p\ge \frac{2}{q}$. If $\sum_{i=1}^l x_i^p
\ge 1$ we are done. Otherwise, consider the sequence
$z=(x_{l+1},x_{l+2},\dots )$. By the choice of $l$, for $i>l$,
$x_i < (2/q)^{1/p}$, and thus\break\vglue-10pt\noindent $\|z\|_{\infty} \leq (2/q)^{1/p}$.
Moreover,
$\frac{1-p}{(2-p)\|z\|_\infty^{1-p}}\ge\frac{\log_2\log_2
q}{2}\cdot \left(\frac{\log_2 q}{2}\right)^{(1-p)/p} \ge 1$, so by
applying Lemma \ref{lem:pinfty} to $z$ we get that
$\|z\|_{p,\infty}^p/\|z\|_1\ge 1$; i.e., there is an integer $b>l$
such that:
\vskip12pt
\hfill $
\displaystyle{(b-l)x_b^p\ge  \|z\|_1=1-\sum_{i=1}^lx_i\ge 1-\sum_{i=1}^lx_i^p.}
$ 
\Endproof\vskip12pt

We are now in position to prove the main technical lemma.

\demo{Proof of Lemma~{\rm \ref{lem:qPhi}}}
For simplicity denote $\beta(\Phi)=[t\log_2(4q\Phi)]^{-2/t}$. We
will prove by induction on $n$ that any $n$ point weighted metric
space $(M,d,w)$ contains a subspace $N\subset M$ such that:
$$
\sum_{x\in N}w(x)^{\beta(\Phi(M))}\ge \left(\sum_{x\in
M}w(x)\right)^{\beta(\Phi(M))},
$$
and for which there is a noncontractive, $4t$-Lipschitz embedding
of $N$ into an ultrametric $H$ with $\diam(H)=\diam(M)$.

In what follows for every $S\subset M$ we denote $w(S)=\sum_{x\in
S}w(x)$.

Let $(M,d,w)$ be a weighted $n$-point metric space such that $w$
is $q$-de\-composable. Without loss of generality we may assume that
$w(M)=1$ and $ \min_{x\neq y\in M} d(x,y)=1. $ Denote
$\Phi=\Phi(M)$. The latter assumption implies that
$\Phi=\diam(M)$. In what follows we denote for $r>0$ and $x\in M$,
$B(x,r)=\{y\in M;\ d(y,x)<r\}$. The proof proceeds by proving the
following claim:

\begin{claim}\label{claim:find i} There exist  $i\in
\{1,\dots ,t\}$ and $x_0\in M$ such that if   $A=\{x_0\}\cup
B\left(x_0,\frac{(i-1)\Phi}{4t}\right)$ and $B=M\setminus
B\left(x_0,\frac{i\Phi}{4t}\right)$ then\/{\rm :}\/
\begin{equation}\label{eq:cases}
\max\left\{\frac{w(A)^{\beta(\Phi/2)}}{[\max_{y\in A}
w(y)]^{\beta(\Phi/2)-\beta(\Phi)}},w(A)^{(\log_2
q)^{-1/(t-1)}}\right\}+w(B)\ge 1.
\end{equation}
\end{claim}

Before proving Claim \ref{claim:find i} we will show how it
implies the required result. Let $i,A,B$ be as in Claim
\ref{claim:find i}. Note that $A\neq \emptyset$ and
$\diam(A)<\frac{\Phi}{2}<\diam(M)$. In particular it follows that,
$|A|,|B|<n$ so that by the induction hypothesis there are
subspaces $A'\subset A$ and $B'\subset B$ such that
\begin{eqnarray*}
\sum_{x\in A'}w(x)^{\beta(\Phi(A))}&\ge& w(A)^{\beta(\Phi(A))}\\
\noalign{\noindent
 and}  \sum_{x\in B'}w(B)^{\beta(\Phi(B))}&\ge&
w(B)^{\beta(\Phi(B))}\ge w(B),
\end{eqnarray*}
HST's $X$ and $Y$ with $\diam(X)=\diam(A)$, $\diam(Y)=\diam(B)$,
and noncontractive embeddings $f:A'\to X$, $g:B'\to Y$ which are
$4t$-Lipschitz. Let $T$ be the tree defining $X$ and $u$ be its
root. Let $S$ be the tree defining $Y$ and $v$ be its root. Define
a tree $R$ as follows: its root is $r$ and the only two subtrees
emerging from it are isomorphic to $T$ and $S$. Label the root of
$R$ by setting $\Delta(r)=\diam(M)$, and leave the labels of $T$
and $S$ unchanged. Note that 
\begin{eqnarray*}
\Delta(r) = \diam(M) 
& \ge &
\max\{\diam(A),\diam(B)\}\\
&=&\max\{\diam(X),\diam(Y)\}=\max\{\Delta(u),\Delta(v)\},
\end{eqnarray*} 
so that with these definitions the leaves of $R$, $X\cup Y$, form
an HST with $\diam(X\cup Y)=\Phi=\diam(M)$. Define $h:A'\cup B'\to
X\cup Y$ by $h|_{A'}=f$ and $h|_{B'}=g$. If $a\in A'$ and $b\in
B'$ then $d(h(a),h(b))=\Phi\ge d(a,b)$. Hence $h$ is
noncontracting. Additionally:
$$
d(a,b)\ge d(b,x_0)-d(a,x_0)\ge \frac{\Phi
i}{4t}-\frac{\Phi(i-1)}{4t}=\frac{\Phi}{4t}=\frac{d_R(h(a),h(b))}{4t},
$$
so that $h$ is $4t$-Lipschitz.

Observe that since $\beta(\Phi)\le\beta(\Phi(A))\le (\log_2
q)^{-1/(t-1)}$ and $w(x)\le 1$ (point-wise),
$$
\sum_{x\in A'}w(x)^{\beta(\Phi)}\ge \sum_{x\in
A'}w(x)^{\beta(\Phi(A))}\ge w(A)^{\beta(\Phi(A))}\ge w(A)^{(\log_2
q)^{-1/(t-1)}}.
$$
Moreover, since $\Phi(A)\le \Phi/2$: 
\begin{eqnarray*}
\sum_{x\in A'}w(x)^{\beta(\Phi)} &\ge& \frac{1}{[\max_{y\in A}
w(y)]^{\beta(\Phi(A))-\beta(\Phi)}}\sum_{x\in A'}w(x)^{\beta(\Phi(A))}\\
&\ge& [\max_{y\in A}
w(y)]^{\beta(\Phi)}\left[\frac{w(A)}{\max_{y\in A}
w(y)}\right]^{\beta(\Phi(A))}\\
&\ge& [\max_{y\in A}
w(y)]^{\beta(\Phi)}\left[\frac{w(A)}{\max_{y\in A}
w(y)}\right]^{\beta(\Phi/2)}. \
\end{eqnarray*} 
We deduce that:
$$
\sum_{x\in A'}w(x)^{\beta(\Phi)}\ge
\max\left\{\frac{w(A)^{\beta(\Phi/2)}}{[\max_{y\in A}
w(y)]^{\beta(\Phi/2)-\beta(\Phi)}},w(A)^{(\log_2
q)^{-1/(t-1)}}\right\},
$$
so that by (\ref{eq:cases}),
$$
\sum_{x\in A'\cup B'}w(x)^{\beta(\Phi)}\ge \sum_{x\in
A'}w(x)^{\beta(\Phi)}+w(B)\ge 1,
$$
as required.

\demo{Proof of Claim {\rm \ref{claim:find i}}} The fact that
$w$ is $q$-decomposable implies that we can split $M=N_1 \dot \cup N_2$,
so that $w(x)\ge \frac{1}{q}$
for every $x\in N_1$ and there is $\omega>0$ such that $w(x)=\omega$ for
every $x\in N_2$. We distinguish between two cases:

\demo{Case {\rm 1}} $\diam_M(N_1)> \frac{\Phi}{2}$. In this
case there are $x_0,x_0'\in N_1$ such that $d(x_0,x_0')>
\frac{\Phi}{2}$. It follows in particular that $B(x_0,\Phi/4)\cap
B(x_0',\Phi/4)=\emptyset$ so that by interchanging the roles of
$x_0$ and $x_0'$, if necessary, we may assume that
$w(B(x_0,\Phi/4))\le \frac{w(M)}{2}=\frac{1}{2}$. Since $x_0\in
N_1$, $w(x_0)\ge \frac{1}{q}$. This implies that there exist
$i\in\{1,\dots ,t\}$ such that
$$ w\left(\{x_0\}\cup
B\left(x_0,\frac{(i-1)\Phi}{4t}\right)\right)^{(\log_2
q)^{-1/(t-1)}}\ge
w\left(B\left(x_0,\frac{i\Phi}{4t}\right)\right),
$$
since otherwise:  
\begin{eqnarray*}
\frac{1}{2}&\ge&
w\left(B\left(x_0,\frac{\Phi}{4}\right)\right)\\
&>&w\left( \{x_0\}\cup
B\left(x_0,\frac{(t-1)\Phi}{4t}\right)\right)^{(\log_2
q)^{-1/(t-1)}}\\
&>&\dots >w\left(
B\left(x_0,\frac{\Phi}{4t}\right)\right)^{(\log_2 q)^{-1}}\\
&\ge& w(x_0)^{(\log_2 q)^{-1}}\ge \frac{1}{q^{(\log_2 q)^{-1}}}=
\frac{1}{2},
\end{eqnarray*} 
which is a contradiction. Fixing such an index $i$ and defining
$A,B$ as in the statement of Claim \ref{claim:find i} we get that:
\begin{eqnarray*}
w(A)^{(\log_2 q)^{-1/(t-1)}} + w(B) 
&=& w\left(\{x_0\}\cup
B\left(x_0,\frac{(i-1)\Phi}{4t}\right)\right)^{(\log_2
q)^{-1/(t-1)}}\\[4pt]
&&+\left[1- 
w\left(B\left(x_0,\frac{i\Phi}{4t}\right)\right)\right]\ge1,
\end{eqnarray*} 
which proves (\ref{eq:cases}).

\demo{Case $2$} $\diam_M(N_1)\le \frac{\Phi}{2}$. In this
case take $x_0\in M$ such that $d(x_0,N_1)=\max_{x\in M}d(x,N_1)$.
We claim that this implies that $N_1\cap B(x_0,\Phi/4)=\emptyset$.
Indeed, otherwise it will follow that $d(x_0,N_1)<\Phi/4$ so that
by the choice of $x_0$, for every $x,y\in M$,
$$
d(x,y)\le
d(x,N_1)+d(y,N_1)+\diam(N_1)<2d(x_0,N_1)+\frac{\Phi}{2}<\Phi,
$$
which is a contradiction.

Set $m=|N_2|$ and denote for $i\in \{0,\dots ,t\}$:
$$
\epsilon_i=\frac{\left|\left(\{x_0\} \cup
B\left(x_0,\frac{i\Phi}{4t}\right)\right)\cap N_2\right|}{m}.
$$
Note that since $x_0\in N_2$, $m^{-1}=\epsilon_0\le
\epsilon_{t}\le 1$. We claim that this implies that there is some
$i\in \{1,\dots , t\}$ such that:
\begin{equation}\label{eq:condition}
\epsilon_{i-1}^{\beta(\Phi/2)}m^{\beta(\Phi/2)-\beta(\Phi)}\ge
\epsilon_{i}.
\end{equation}
Indeed, if we set $a=\log_2(2q\Phi)\ge 1$, then if there is no
such $i$ we have for every $i\in \{1,\dots ,t\}$:
$$
\epsilon_{i-1}<\left(\frac{\epsilon_i}{m^{\frac{1}{(ta)^{2/t}}- 
\frac{1}{[t(a+1)]^{2/t}}}}\right)^{(ta)^{2/t}}.
$$
Denote $b=m^{\frac{1}{(ta)^{2/t}}-\frac{1}{[t(a+1)]^{2/t}}}$ and
$c=(ta)^{2/t}$. The above inequality then becomes
$\epsilon_{i-1}<(\epsilon_i/b)^c$. Iterating this $t$ times we
get:
$$
\frac{1}{m}=\epsilon_0<\frac{\epsilon_t^{c^t}}{b^{c+c^2+\dots +c^t}}=\frac{\epsilon_t^{c^t}}
{b^{\frac{c}{c-1}(c^t-1)}}\le \frac{1}{b^{c^t-1}}.
$$
Thus,
$$
m^{(t^2a^2-1)\left[\frac{1}{(ta)^{2/t}}-\frac{1}{[t(a+1)]^{2/ 
t}}\right]}<m,
$$
but an application of the mean value theorem gives a
contradiction, since:
\begin{eqnarray*}
(t^2a^2-1)\left[\frac{1}{(ta)^{2/t}}-\frac{1}{[t(a+1)]^{2/ 
t}}\right]&\ge&\frac{t^2a^2}{2}\frac{2}{t^{1+2/t}(a+1)^{1+2/t}}
\\
&\ge& t^{1-2/t}\left(\frac{a}{a+1}\right)^2\ge
8^{3/4}\cdot\frac{1}{4}\ge 1.
\end{eqnarray*}

Choose an index $i\in \{1,\dots , t\}$ satisfying
(\ref{eq:condition}) and let $A,B$ be as in the statement of Claim
\ref{claim:find i} for this particular $i$. Observe that:
$$
B=M\setminus B\left(x_0,\frac{i\Phi}{4t}\right)\supset M\setminus
B\left(x_0,\frac{\Phi}{4}\right)\supset N_1,
$$
so that:  
\begin{eqnarray*}
\frac{w(A)^{\beta(\Phi/2)}}{[\max_{y\in A}
w(y)]^{\beta(\Phi/2)-\beta(\Phi)}}+w(B) &=&
\frac{(\omega\epsilon_{i-1}m)^{\beta(\Phi/2)}}{\omega^{\beta(\Phi/2)- 
\beta(\Phi)}}+w(N_1)+(1-\epsilon_i)m\omega\\
&=&\omega^{\beta(\Phi)}(\epsilon_{i-1}m)^{\beta(\Phi/2)}+w(N_1)+(1- 
\epsilon_i)m\omega\\
&\ge&
(m\omega)^{\beta(\Phi)}\epsilon_i+w(N_1)+(1-\epsilon_i)m\omega\\
&\ge& m\omega\epsilon_i+w(N_1)+(1-\epsilon_i)m\omega=w(M)=1.
\end{eqnarray*} 
This concludes the proof of Claim \ref{claim:find i}.
\hfill\qed

\Subsec{Passing from an ultrametric to a  
$k${\rm -HST}}\label{section:lower-um-k-hst}
In what follows we show that every ultrametric contains large
subsets which are embeddable in a $k$-HST with distortion $\alpha
< k$.

An unweighted version of the following result was proved in
\cite{bbm}. The bound for the weighted Ramsey function is a
straightforward modification of the proof in \cite{bbm}:
\begin{lemma}[\cite{bbm}] \label{lem:1-hst}
For every $k > \alpha > 1${\rm ,}
$$
\XX_k(\um,\alpha)\ge \frac{1}{\lceil \log_\alpha k \rceil} .
$$
\end{lemma}
\vglue8pt

If $k$ is large with respect to $\alpha$ then the bound of
\cite{bbm} provides a good approximation for $\XX_k(\um,\alpha)$.
In fact, this is how Lemma \ref{lem:1-hst} is used in
Section~\ref{section:lower-small-alpha}.

However, when $k$ is close to $\alpha$ the bound in
Lemma~\ref{lem:1-hst} is not good enough for proving our main
theorem.
  We obtain bounds for this range of parameters in the
following lemma (stated in \S\ref{section:lower-main} in
slightly weaker form).

\demo{\scshape Lemma 3.15}
{\it For every} $k > \alpha > 1$,
$$ \XX_k(\um,\alpha) \ge 1 -\frac{1}{\lceil \log_{k/\alpha} \alpha  
\rceil} . $$
\vskip12pt

  Before proving {Lemma 3.15} we require some lemmas
  concerning unweighted trees.

\begin{definition} Let $h>1$ be an integer and $i\in
\{0,\dots ,h-1\}$. We say that a rooted tree $T$ is $(i,h)$-{\em
periodically sparse} if for every $l\equiv i\ (\mod h)$, every
vertex at depth $l$ in $T$ is degenerate. $T$ is called $h$
periodically sparse if there exists $i\in \{0,\dots ,h-1\}$ for
which $T$ is $(i,h)$-periodically sparse.
\end{definition}

In what follows we always use the convention that a subtree $T'$
of a rooted tree $T$ is rooted at the root of $T$ and that the
leaves of $T'$ are also leaves of $T$. We denote by $\lvs(T)$ the
leaves of $T$.

\begin{lemma}\label{lem:sparse} Fix an integer $h>1$. Let $T$ be a  
finite rooted tree.
Then for any $w:\lvs(T)\to \R^{+}$ there exists a subtree of $T${\rm ,}
$T'${\rm ,} which is $h$ periodically sparse and\/{\rm :}\/
$$
\sum_{v\in \lvs(T')}w(v)^{\frac{h-1}{h}}\ge\left(\sum_{v\in
\lvs(T)}w(v)\right)^{\frac{h-1}{h}}.
$$
\end{lemma}
\vskip8pt

The techniques we use in the proof of Lemma~\ref{lem:sparse} are
similar to those used in the proof of Lemma~\ref{lem:1-hst} in
\cite{bbm}. It can also be derived from a result in \cite{bl2}
concerning influences in multi-stage games. These facts are also
closely related to an isoperimetric inequality of Loomis and
Whitney \cite{lw}.

\Proof  For every $i\in\{0,1,\dots ,h-1\}$ let $f_i(T)$ be the  
maximum
of $\sum_{v\in \lvs(T')}w(v)^{\frac{h-1}{h}}$ over all the
$(i,h)$-periodically sparse subtrees, $T'$, of $T$. We will prove
by induction on the maximal depth of $T$ that:
$$
\prod_{i=0}^{h-1}f_i(T)\ge \left(\sum_{v\in
\lvs(T)}w(v)\right)^{h-1},
$$
from which it will follow that $\max_{0\le i\le h-1}f_i(T)\ge
\left(\sum_{v\in \lvs(T)}w(v)\right)^{\frac{h-1}{h}}$, as
required.

For a tree $T$ of depth $0$, consisting of a single node $v$, we
have that $f_i(T) = w(v)^{\frac{h-1}{h}}$ and therefore
$$
\prod_{i=0}^{h-1}f_i(T)\ge \left(w(v)^{\frac{h-1}{h}}\right)^h =
w(v)^{h-1} .
$$
Assume that the maximal depth of $T$ is at least $1$, let $r$ be
the root of $T$ and denote by $v_1,\dots ,v_l$ its children.
Denote by $T_j$ the subtree of $T$ rooted at $v_j$. Observe that:
$$
f_0(T)\ge \max_{1\le j\le l} f_{h-1}(T_j),
$$
and for $i\in \{1,\dots , h-1\}$:
$$
f_i(T)=\sum_{j=1}^lf_{i-1}(T_j),
$$
By repeated application of H\"older's inequality:
$$
\sum_{j=1}^l\prod_{i=0}^{h-2}[f_i(T_j)]^{\frac{1}{h-1}}\le
\left(\prod_{i=0}^{h-2}\sum_{j=1}^lf_{i}(T_j)\right)^{\frac{1}{h-1}}.
$$
Therefore, by the induction hypothesis:
\begin{eqnarray*}
\prod_{i=0}^{h-1}f_i(T)&\ge& \max_{1\le j\le
l}f_{h-1}(T_j)\cdot\prod_{i=1}^{h-1}\sum_{j=1}^lf_{i-1}(T_j)\\
&\ge& \max_{1\le j\le
l}f_{h-1}(T_j)\cdot\left(\sum_{j=1}^l\left(\prod_{i=0}^{h 
-2}f_i(T_j)\right)^{\frac{1}{h-1}}\right)^{h-1} 
\\
&\ge& \left(\sum_{j=1}^l\left(\prod_{i=0}^{h 
-1}f_i(T_j)\right)^{\frac{1}{h-1}}\right)^{h-1}
\\
&\ge& \left(\sum_{j=1}^l \sum_{v\in\lvs(T_j)}w(v)\right)^{h-1} =
\left(\sum_{v\in \lvs(T)}w(v)\right)^{h-1}.
\end{eqnarray*} 
\vglue-28pt
\Endproof\vskip20pt

Before proving Lemma 3.15, we prove the following variant of
a proposition from \cite{bartal2}:

\demo{\scshape Lemma 3.5}
  {\it For any $k>1${\rm ,} any ultrametric is $k$\/{\rm -}\/equivalent to an exact} $k$-HST.
 
\Proof  Let $T$ be a labelled tree rooted at $r$. Define a
new labelled tree $T'$ as follows. Let $u$ be a minimal depth
vertex in $T$ that has a child $v$ for which $\Delta(u)\neq
k\Delta(v)$. Let $0\le i\in \mathbb{N}$ be defined via $k^i\le
\frac{\Delta(u)}{\Delta(v)}< k^{i+1}$. Relabel\break\vglue-12pt\noindent $v$ by setting
$\Delta'(v)=\frac{\Delta(u)}{k^{i}}\ge \Delta(v)$, and replace the
edge $[u,v]$ by a path of length $i$ whose labels decrease by a
factor $k$ at each step. Denote the tree thus obtained by $T'$. If
we start out with an HST $X$ with defining tree $T$, then the tree
$T'$ produced in this procedure defines a new HST. Iterating this
construction as long as possible, we arrive at a tree $\tilde{T}$
which defines an exact $k$-HST. To prove that we have distorted
the metric by a factor of at most $k$ observe that by
construction, for any $x,y\in X$,
$\lca_T(x,y)=\lca_{\tilde{T}}(x,y)$ and that for any $v\in T\cap
\tilde{T}$, $\Delta_T(v)\le \Delta_{\tilde{T}}(v)\le k
\Delta_T(v)$.
\Endproof\vskip4pt

{\it Proof of Lemma} 3.15.
Let $h = \lceil \log_{k/\alpha} \alpha \rceil$ and let $s =
k^{1/h}$. By Lemma 3.5, $X$ is $s$-equivalent to some exact
$s$-HST $Y$ via a noncontractive $s$-Lipschitz embedding. Let $T$ be
the tree defining $Y$. Lemma \ref{lem:sparse} yields a subtree $S$
of $T$ which is $(i,h)$-periodically sparse for some $i\in
\{0,\dots ,h-1\}$, such that
$$
\sum_{v\in S}w(g^{-1}(v))^{\frac{h-1}{h}}\ge \left(\sum_{x\in
X}w(x)\right)^{\frac{h-1}{h}}.
$$
By attaching a path of length $h-1-i$ to the root of $S$ we may
assume that $S$ is $(h-1,h)$-periodically sparse. Similarly, by
adding appropriate paths to the leaves of $S$ we may assume that
there is an integer $m$ such that all the leaves of $S$ are at
depth $mh$. Denote by $r$ the root of $S$. We change the tree $S$
as follows. For every integer $0\le j<m$ delete all the vertices
of $S$ whose depth is in the interval $[jh+1,(j+1)h-1]$ and
connect every vertex of depth $jh$ directly to all its descendants
of depth $(j+1)h$. Denote the tree thus obtained by $S'$ and
denote by $Y'$ the metric space induced by $S'$ on $Y$. It is
evident that $Y'$ is an exact $s^h$-HST. We claim that $Y'$ is
$s^{h-1}$ equivalent to a subspace of $X$ via a noncontractive
$s^{h-1}$ Lipschitz embedding. Indeed, fix $u,v\in Y'$ and let $w$
be their least common ancestor in $S$. If we denote by $q$ the
depth of $w$ in $S$ then the key observation is that since $S$ is
$(h-1,h)$-periodically sparse, $q\not\equiv (h-1)\ \pmod{h}$. We
can therefore write $q=i+jh$ for some $i\in \{0,\dots ,h-2\}$ and
$j\ge 0$. If we denote by $w'$ the least common ancestor of $u,v$
in $S'$ then by the construction, $w'$ is in depth $jh$ in $S$.
Hence $d_Y(u,v)=\frac{\Delta(r)}{s^{i+jh}}$ and
$d_{Y'}(u,v)=\frac{\Delta(r)}{s^{jh}}$, so that:
$$
d_Y(u,v)\le d_{Y'}(u,v)\le s^{i}d_Y(u,v)\le s^{h-2}d_Y(u,v).
$$
This shows that $Y'$ is $s^{h-2}$ equivalent to $Y$ via a
noncontractive $s^{h-2}$ Lipschitz embedding. Since $Y$ is $s$
equivalent to a subspace of $X$ via a noncontractive $s$
Lipschitz bijection we have that $Y'$ is $s^{h-1}$ equivalent to a
subspace of $X$.

Recall that $s^h = k$, and it remains to show that $s^{h-1} \le
\alpha$. Indeed by our choice of $h$, $h-1 \le \log_{k/\alpha}
\alpha$, or $\frac{1}{h-1} \ge \log_\alpha (k/\alpha)$. Therefore
$\alpha^{\frac{h}{h-1}} \ge k$, and so $s^{h-1} =
k^{\frac{h-1}{h}} \le \alpha$.
\hfill\qed

\Subsec{Passing from a $k${\rm -HST} to metric  
composition}\label{section:lower-k-hst-comp}
In this section we prove that if a metric space is close to a
$k$-HST then it is very close to a metric space in the composition
closure of a class of metric spaces with low aspect ratio.

\sdemo{Lemma 3.16}
{\it For any $\alpha,\beta \ge 1${\rm ,} if a metric space $L$ is
$\alpha$\/{\rm -}\/equivalent to a $\beta\alpha${\rm -HST} then $L$ is
$(1+2/\beta)$\/{\rm -}\/equivalent to a metric space in
$\comp_\beta(\alpha)$.}

\Proof 
Let $L$ be a metric space. Let $k=\beta\alpha$. Let $X$ be a
$k$-HST such that there is an $\alpha$ Lipschitz noncontractive
bijection $f:L\to X$. Namely, for every $x,y\in L$, $d_L(x,y)\le
d_X(f(x),f(y))\le \alpha d_L(x,y)$.

  Let $T$ be the tree defining
$X$.   For a vertex $u \in T$, let $T_u$ be the subtree of $T$
rooted at $u$. Let $X_u$ (a subspace of $X$) denote the HST
defined by $T_u$ and $L_u = f^{-1}(X_u)$. Then $\diam(L_u) \le
\diam(X_u) = \Delta(u)$.

Our goal is to build a metric space $Z \in \comp_\beta(\alpha)$
along with a noncontractive Lipschitz bijection $g:L\to Z$ which
satisfies for every $x,y\in L$, $d_L(x,y)\le d_Z(g(x),g(y))\le
\left(1+\frac{2}{\beta}\right) d_L(x,y)$. We prove this by
induction on the size of $L$. The inductive hypothesis needs to be
further strengthened with the requirement that $\diam(Z) \le
\diam(L) = \Delta$.

Let $r$ be the root of $T$, with $\Delta(r) = \Delta$. Let $C$
denote the set of children of $r$. By induction, there exists for
each child $u \in C$ a metric space $N_u \in \comp_\beta(\alpha)$
and a noncontractive Lipschitz bijection $g_u:L_u\to N_u$ which
satisfies for every $x,y\in L_u$, $d_{L_u}(x,y)\le
d_{N_u}(g_u(x),g_u(y))\le \left(1+\frac{2}{\beta}\right)
d_{L_u}(x,y)$. Also $\diam(N_u) \le \diam(L_u) = \Delta(u)$.

Define a metric space $M = (C,d_M)$ by setting for every $u \neq v
\in C$,
$$
d_M(u,v) = \max \{ d_L(x,y);\  x \in L_u, y \in L_v \} .
$$
Fix $u \neq v \in C$ and $x \in L_u, y\in L_v$. Since
$d_X(f(x),f(y)) = \Delta$, we have that $\Delta/\alpha \le
d_L(x,y) \le \Delta$. It follows that for every $u \neq v \in C$
and $x \in L_u, y\in L_v$,
$$
\frac{\Delta}{\alpha} \le d_L(x,y) \le d_M(u,v) \le \diam(L) =
\Delta.
$$
Therefore $\Phi(M) \le \alpha$ and $\diam(M) \le \Delta$. Also for
every $u,v,x,y$ as above,
\begin{eqnarray*}
d_L(x,y)\,\, \le\,\, d_M(u,v) & \le & d_L(x,y) + \diam(L_u) +  
\diam(L_v) \\
& = & d_L(x,y) + \Delta(u) + \Delta(v)
\le d_L(x,y) + 2\frac{\Delta}{k} \\
& \le & d_L(x,y) + 2\frac{\alpha d_L(x,y)}{\beta\alpha}
  \le  \left(1+\frac{2}{\beta}\right) d_L(x,y) .
\end{eqnarray*}

Now, we let
\begin{eqnarray*}
\gamma  =  \frac{\max_{u \in C} \diam(N_u)}{\min_{u,v \in C}
d_M(u,v)}, \quad\quad {\rm and} \quad \beta' = \frac{1}{\gamma}
\ge \frac{\Delta/\alpha}{\Delta/k} = \beta .
\end{eqnarray*}

Define $Z \in \comp_\beta(\alpha)$, by letting $Z =
M_{\beta'}[\N]$, where $\N = \{ N_u \}_{u \in C}$. Also define for
every $u \in C$ and $x \in X_u$,  $g(x) = g_u(x)$.

Let $u,v \in C$ and $x \in L_u, y\in L_v$. When $u=v$ the bound on
the distortion of $g$ follows from our induction hypothesis. For
$u \neq v$, $d_Z(g(x),g(y)) = \beta'\gamma d_M(u,v) = d_M(u,v)$,
which implies the required bound on the distortion of~$g$, and the
requirement $\diam(Z) \le \Delta$.
\hfill\qed

\Subsec{Distortions arbitrarily close to  
$2$}\label{section:lower-small-alpha}
Our goal in this section is to prove the following theorem:

\begin{theorem}\label{thm:lower-small-alpha} There is an absolute  
constant $c>0$ such that for any
$k \ge 1$ and $0<\epsilon<1${\rm ,} for any integer $n$\/{\rm :}\/
$$
R_{k\text{-\hst}}({2+\epsilon},n)\ge n^{\frac{c\epsilon}{\log(2
k/\epsilon)}}.
$$
In particular{\rm ,}
$$
R_{\um}({2+\epsilon},n)\ge n^{\frac{c\epsilon}{\log(2/\epsilon)}}.
$$
\end{theorem}
\vskip8pt

By Proposition ~\ref{prop:hst-l2} the same bound holds for
$R_2(2+\epsilon,n)$.

As in the case of large $\alpha$, we derive
Theorem~$3.26$ from the following stronger
claim.

\sdemo{Theorem $3.26'$}
{\it There is an absolute constant $c>0$ such that for any $k \ge 1$
and} $0<\epsilon<1$:
$$
\XX_k(2+\epsilon)\ge \frac{c\epsilon}{\log(2 k/\epsilon)}.
$$
\Enddemo

The proof of Theorem~$3.26'$ uses most of the techniques
developed for the case of large $\alpha$. The basic idea is 
first to  apply Theorem~{$3.7'$} to obtain some constant $\alpha'$ for
which there is a constant bound on $\psi(\alpha')$, e.g.\ $1/2$.
So, our goal is to find another subspace for which we can
improve the distortion from $\alpha'$ to $2+\epsilon$. Again, we
would like to exploit metric spaces with low aspect ratio $\Phi$.
For large $\alpha$ we could do this with $\Phi$ bounded with
respect to $\alpha$. Since we started with some constant
$\alpha'$ we can expect $\Phi$ to be constant as well.
However, the bound of Lemma~\ref{lem:Phi1} does not apply for
small values of $\alpha$. Thus, our first step is to obtain
meaningful lower bounds on $\XX_k(\Phi,2+\epsilon)$ for every
$\epsilon > 0$. This is done by giving a lower bound on
$\XX_{\eq}(\Phi,2+\epsilon)$, that is by finding a large
equilateral subspace, which is a special case of a $k$-HST. We can
now apply Lemma~\ref{lem:composition-xx1} to get lower bounds on
$\XX_k(\comp_\beta(\Phi),2+\epsilon)$. For large $\alpha$ we were
able to extend such bounds by finding a subspace close to a
$k$-HST, and therefore very close to a metric space in
$\comp_\beta(\Phi)$ via Lemma~3.16. In the
present case, "very close" means distortion $\approx 1+\epsilon$,
which implies that $k$ and $\beta$ must be in the range of
$1/\epsilon$. This is achieved by initially applying
Lemma~\ref{lem:1-hst} to get a bound on $\psi_k(\alpha')$.

We begin with a proof of the bound on $\XX_k(\Phi,2+\epsilon)$,
  which is based on bounds on
embedding into an equilateral space. We start with the following
result:

\begin{lemma}\label{lem:diameter}
Let $\alpha>2${\rm ,} $s\ge 2$ be real numbers and $t\geq 1$ be an
integer. Let $M$ be an $n$ point metric space. Then at least one
of the following two conditions holds\/{\rm :}\/
\begin{enumerate}
\item $M$ contains a subspace $N$ of size at least $s$ that is
$\alpha$\/{\rm -}\/equivalent to an equilateral space. \item $M$ contains a
subspace $N$ of size at least $n/s^t${\rm ,} such that $\diam(N) <
(\alpha/2)^{-t} \diam(M)$.
\end{enumerate}
\end{lemma}

\Proof 
By induction on $t$. Suppose that $M$ has no subspace of size $s$
that is $\alpha$ equivalent to an equilateral space. For $t=1$ let
$N_0 = M$. For $t>1$ we get by the induction hypothesis there is a
subspace $N_{t-1}\subseteq M$ which contains at least $n/s^{t-1}$
points and $\diam(N_{t-1}) \le (\alpha/2)^{-t+1}\diam(M)$.

Let $\{c_1,\dots , c_r\}$ be a maximal subset of $N_{t-1}$ such
that $$d(c_i,c_j)\geq \diam(N_{t-1})/\alpha$$ for $i \neq j$. Since
$\{c_1,\dots  c_r\}$ is $\alpha$ equivalent to an equilateral
space, our assumption implies that $r\leq s$. Let $C_i=N_{t-1}\cap
B(c_i,\diam(N_{t-1})/\alpha)$. By the maximality of $r$,
$\boldsymbol{\cup}_{i=1}^r C_i=N_{t-1}$, and so if we set $N_t$ to be the
largest $C_i$, we have that its cardinality is at least
$|N_{t-1}|/r\geq n/s^t$. Now:
\vglue12pt
\centerline{$
{\diam(N_t)\leq \diam(B(c_i,\diam(N_{t-1})/\alpha)) <
\frac{2}{\alpha} \diam(N_{t-1}) \le
\left(\frac{2}{\alpha}\right)^t \diam(M).}
$}
\vglue-8pt
\phantom{overthere}\Endproof\vskip12pt 

This implies a bound on the cardinality of a subspace that is
$\alpha$-equivalent to an equilateral space.

\begin{corollary}\label{coro:equi-phi} Fix $\alpha>2$ and an
integer $n\ge 4$. Let $M$ be a metric space of size $n$. Then{\rm ,}
$$
R_\eq(M;\alpha, n) \ge {\left(\frac{n}{2}\right)}^{{\left\lceil
\log_{\alpha/2} \Phi(M) \right\rceil} ^{-1}} \ge
n^{\frac{1}{2}{\left\lceil \log_{\alpha/2} \Phi(M) \right\rceil}
^{-1}}.
$$
\end{corollary}

\Proof  Apply Lemma \ref{lem:diameter} with $t=\left\lceil
\log_{\alpha/2} \Phi(M) \right\rceil$ and $s=(n/2)^{1/t}$. We
obtain a subspace $N$ of $M$. All we have to do is verify that
with these parameters the second condition in Lemma
\ref{lem:diameter} cannot hold. Indeed, otherwise $|N|\ge n/s^t=2$
so that $\diam(N)\ge \min_{x\neq y}d_M(x,y)$, and it follows that
$(\alpha/2)^t < \Phi$, which contradicts the choice of $t$.
\Endproof\vskip4pt 

We show next that Corollary~\ref{coro:equi-phi} implies bounds for
the weighted Ramsey problem, and so we can bound
$\XX_\eq(\Phi,\alpha)$ for any $\alpha>2$. Since an equilateral is
in particular a $k$-HST, we get a bound on $\XX_k(\Phi,\alpha)$.
To obtain this we need to extend the bound in
Corollary~\ref{coro:equi-phi} to hold for the weighted Ramsey
problem. To achieve this we make use of another lemma from
\cite{bbm}, which is similar in flavor to Lemma~\ref{lem:lb}:

\begin{lemma}[\cite{bbm}]\label{lem:bbm-bin-bal} Let  
$x=\{x_i\}_{i=1}^\infty$ be a
sequence of nonnegative real numbers. Then there exists a
sequence $y=\{y_i\}_{i=1}^\infty$ such that $y_i \leq x_i$ for all
$i \ge 1$ and\/{\rm :}\/
$$
\sum_{i\ge 1} y_i^{1/2} \ge \left(\sum_{i\ge 1} x_i\right)^{1/2}.
$$
Moreover{\rm ,} one of the following two cases holds true\/{\rm :}\/
\begin{enumerate}
\item For all $i > 2${\rm ,} $y_i = 0$.

\item There exists $\omega > 0$ such that for all $i \ge 1$ either
$y_i = \omega$ or $y_i = 0$.
\end{enumerate}
\end{lemma}

\begin{corollary}\label{coro:equi-chi} For any $k \ge 1${\rm ,} $\alpha>2$ and
$\Phi>
1${\rm ,}
$$
\XX_k(\Phi,\alpha)\ge \XX_{\eq}(\Phi,\alpha) \ge
\frac{1}{4}{{\left\lceil \log_{\alpha/2} \Phi \right\rceil}
^{-1}}.
$$
\end{corollary}
\vskip8pt

\Proof  Let $M$ be an $n$-point metric space with aspect ratio  
$\Phi(M)\le \Phi$.
Let $w:M\to \R^{+}$ be a weight function normalized so that
$\sum_{x\in M} w(x)=1$. Apply Lemma~\ref{lem:bbm-bin-bal} to the
sequence $\{w(x)\}_{x \in M}$ to obtain a sequence
$\{w'(x)\}_{x\in M}$ such that for all $x\in M$, $w'(x) \le w(x)$.
In addition, either (i) There are $u,v\in M$ such that
$w(u)^{1/2}+w(v)^{1/2} \ge w'(u)^{1/2}+ w'(v)^{1/2}\ge 1$, or
there is a subset $N\subseteq M$ such that for all $x \in N$,
$w(x) \ge w'(x) = \omega > 0$, and $|N|\omega^{1/2}\ge 1$. In the
first case the subset $\{u,v\}$ is isometric to an equilateral
space and we are done. In the second case, if $|N| \le 4$ then
$\omega \ge 1/16$. Hence, we can choose two points $u',v' \in N$
such that $w(u')^{1/4}+ w(v')^{1/4}\ge 2 \omega^{1/4} \ge 1$.
Again, $\{u',v'\}$ is isometric to an equilateral space.
Otherwise, $|N|
> 4$ and by Corollary \ref{coro:equi-phi} there is a subspace
$N'\subseteq N$ which is $\alpha$-equivalent to an equilateral
space and $ |N'|\ge |N|^{\frac{1}{2}{\left\lceil \log_{\alpha/2}
\Phi \right\rceil} ^{-1}}. $ Hence:
\begin{eqnarray*}
\sum_{x\in N'}w(x)^{\frac{1}{4}{\left\lceil \log_{\alpha/2} \Phi
\right\rceil} ^{-1}}&\ge& |N|^{\frac{1}{2}{\left\lceil
\log_{\alpha/2} \Phi \right\rceil} ^{-1}}
\omega^{\frac{1}{4}{\left\lceil
\log_{\alpha/2} \Phi \right\rceil} ^{-1}}\\
&=&(|N|\omega^{1/2})^{\frac{1}{2}{\left\lceil \log_{\alpha/2} \Phi
\right\rceil} ^{-1}}\ge 1.
\end{eqnarray*}
\vglue-20pt
\Endproof\vskip12pt

\demo{Proof of Theorem~$3.26'$}
Theorem~{$3.7'$} implies in particular that there is a constant
$\theta$ for which $\XX(\theta/2) \ge 1/2$. In other words, given
a metric space $X$, there exists a subspace $X'$ of $X$ which is
$(\theta/2)$-equivalent to an ultrametric $Y$ and satisfies the
weighted Ramsey condition~\eqref{eq:weighted} with $\psi=1/2$.

Let $\beta = 8k/\epsilon$ and $k' = \theta\beta$. It follows from
Lemma~\ref{lem:1-hst} that $Y$ contains a subspace $Y'$ which is
$2$-equivalent to a $k'$-HST and satisfies
condition~\eqref{eq:weighted} with $\XX_{k'}(\um,2) \ge \lceil
\log k' \rceil^{-1}$.  By mapping $X$ into an ultrametric $Y$ and
its image in $Y$ into a $k'$-HST, we can apply
Lemma~\ref{lem:xx-product1} to obtain a subspace $X''$ of $X$ that
is $(\theta/2)\cdot 2 = \theta$-equivalent to a $k'$-HST, and
satisfies condition~\eqref{eq:weighted} with
$$
\XX_{k'}(\theta) \ge \XX_{k'}(\um,2) \cdot \XX(\theta/2) \ge
\frac{1}{2\lceil \log k' \rceil} .
$$
Now, $X''$ is $\theta$-equivalent to a $\theta\beta$-HST and so
Lemma~3.16 implies that it is
$(1+2/\beta)$-equivalent to a metric space $Z$ in
$\comp_\beta(\theta)$. Therefore
$$
\XX_{\comp_\beta(\theta)}(1+2/\beta) \ge \XX_{k'}(\theta) \ge
\frac{1}{2\lceil \log k' \rceil}.
$$
Additionally, using Lemma~\ref{lem:composition-xx1} and the bound
of Corollary~\ref{coro:equi-chi}, we have that there is a constant
$c'$ such that
$$
\XX_k\left(\comp_\beta(\theta),2+\frac{\epsilon}{4}\right) =
\XX_k\left(\theta,2+\frac{\epsilon}{4}\right) \ge
\frac{c'\epsilon}{\log \theta}.
$$
It follows that $Z$ contains a subspace $Z'$ which is
$(2+\epsilon/4)$-equivalent to a $k$-HST.

By mapping $X$ into $Z \in \comp_\beta(\theta)$, and then its
image in $Z$ into a\break $k$-HST, we can apply
Lemma~\ref{lem:xx-product1} to obtain a subspace of $X$ which is\break
 $(2+\epsilon/4)(1+2/\beta) \le (2+\epsilon)$-equivalent to a
$k$-HST and which satisfies the weighted Ramsey
condition~\eqref{eq:weighted} with
$$
\XX_k(2+\epsilon) \ge
\XX_k\left(\comp_\beta(\theta),2+\frac{\epsilon}{4}\right) \cdot
\XX_{\comp_\beta(\theta)}\left(1+\frac{2}{\beta}\right) \ge
\frac{c'\epsilon}{2\log \theta \lceil \log(8\theta k/\epsilon)
\rceil},
$$
which implies the theorem by an appropriate choice of $c$.
\hfill\qed

\section{Dimensionality based upper bounds}\label{section:counting}

In this section we prove some upper bounds on the Ramsey function
of low dimensional spaces. In particular, these imply bounds on
the Euclidean Ramsey function $R_2(\alpha,n)$. In addition, these
bounds show that the lower bounds for low dimensional $\ell_p$
spaces from Corollary \ref{thm:plowdim} in the
introduction are nearly tight. Our upper bounds on $R_2(\alpha,n)$
for $\alpha<2$ improve the results of \cite{bfm} by showing that
for any $\alpha<2$, $R_2(\alpha,n)\leq 2\log_2 n +C(\alpha)$. The
bounds obtained on $R_2(\alpha,n)$ for $2 <\alpha \le \log n/\log
\log n$, are also possibly tight.

   The proof technique we employ here
originates from a counting argument by Bourgain
\cite{bourgainembedding} and later variants (see \cite{matbook}).
A different argument, based on geometric considerations, uses
expander graphs. Expander graphs, in fact yield the best upper
bound we have on $R_p(\alpha,n)$ for $\alpha \ge 2$ and all $p
\geq 1$. This is shown in Section~\ref{section:expanders}.

In this section we prove the following bounds:

\begin{theorem}\label{thm:lowdim}
Let $X$ be an $h$\/{\rm -}\/dimensional normed space and $n$ be an integer.
Then
\begin{itemize}
\item For any $1 < \alpha < 2${\rm ,} $ R_X(\alpha,n)  \leq  2\log_2 n +
2h \log_2\left( \frac{C}{2-\alpha}\right)$. \item For any $\alpha
\geq 2 ${\rm ,} $R_X(\alpha,n)  \leq C n^{1-c/\alpha} h\log\alpha $
\end{itemize}
where $c,C >0$ are some absolute constants.
\end{theorem}

Using the Johnson-Lindenstrauss dimension reduction
Lemma~\cite{jl} we\break derive the following bounds for the Euclidean
Ramsey function $R_2(\alpha,n)$.

\begin{corollary}\label{corol:basic bounds}
There are absolute constants $c,C>0$ such that for\break every integer
$n${\rm ,}
\begin{itemize}
\item For any $0 < \epsilon \leq 1${\rm ,}
  $ R_2(2-\epsilon,n)  \leq  2\log_2 n + C \frac{\log^2  
(2/\epsilon)}{\epsilon^2} $.

\item For any $2 \le \alpha \le \frac{\log n}{\log\log n}${\rm ,}
  $ R_2(\alpha,n)  \leq  Cn^{1-c/\alpha}$.
\end{itemize}
\end{corollary}

The counting argument presented below is based on the existence of
dense graphs for which all metrics defined on subgraphs are very
far from each other.

Let $A$ be a set of vertices in the graph $G=(V,E)$. We denote by
$E_A$ the set of edges in $G$ with both vertices in $A$, and the
cardinality of $E_A$ by $e_A$. The {\em density} of $A$ is
$\frac{e_A}{{|A| \choose 2}}$.\smallbreak

We first explain the relevance of large girth and high density to
our problem. Let $G=(V,E)$ be a large graph of girth $g$ with no
large sparse subgraphs. With every $H \subseteq E$ we associate a
metric on $V$ defined by
$\rho_H(u,v)=\min\left\{g-1,d_H(u,v)\right\},$ where $d_H$ is the
shortest path metric in the subgraph of $G$ with edge set $H$.
Below we show that among these $\rho_H$ are metrics that cannot be
embedded with small distortion in any low-dimensional normed
space.
\begin{lemma}\label{lemma:girth_and_density}
If there exists a graph $G=(V,E)$ of size $n$ with girth at least
$g${\rm ,} in which every set of $\ge s$ vertices has density at least
$q${\rm ,} then for every $h$\/{\rm -}\/dimensional normed space $X$ and every
real $1\le \alpha < g-1${\rm ,}
$$
R_{X}(\alpha,n) \leq \max \left\{s, \frac{1+2}{q}\left[h
\log_2\left( \frac{14\alpha g}{g-\alpha-1}\right) + \log_2
\left(\frac{n}{s}\right)\right] \right\}.
$$
\end{lemma}

\Proof 
To prove the theorem, we may certainly assume that
$R_{X}(\alpha,n) \ge s$. Let $k=R_{X}(\alpha,n)$; namely, every
$n$ point metric space contains a subset of size $k$ that
$\alpha$-embeds in $X$. In particular, for every $H \subseteq E$
there is a set of $k$ vertices $A_H\subseteq V$ such that
$(A_H,\rho_H)$ $\alpha$-embeds into $X$. Therefore there is a
certain set $A$ of $k$ vertices, that is suitable for many sets $H
\subseteq E$. That is, there is a class $\cal H$ of at least
$2^{|E|}/\tbinom{n}{k}$ subsets $H \subseteq E$ for which $A=A_H$,
and therefore  $(A,\rho_H)$ $\alpha$-embeds into $X$. Consider
$H_1, H_2 \in \cal H$ equivalent if $H_1 \cap E_A = H_2 \cap E_A$.
There are at most $2^{|E| - e_A}$ members in $\cal H$ that are
equivalent to a given set $H$. Consequently, there are at least
$2^{e_A}/\tbinom{n}{k}$ subsets $H \subseteq E$ which are mutually
inequivalent and for which $(A,\rho_H)$ $\alpha$-embeds into $X$.
Let $f_H:A\rightarrow X$ be such an embedding, i.e., for every
$u,v\in A$:
$$\frac{1}{\alpha}\rho_H(u,v) \le \|f_H(u)-f_H(v)\|_X\le \rho_H(u,v). $$
Since $\rho_H$ takes values in $\{0,1,\dots ,g-1\}$, by applying
an appropriate translation we may assume that $f_H(A)\subseteq
B_X(0,g)$. We now ``round" the images $f_H(A)$ to the points of a
$\delta$-net in $B_X(0,g)$, where $\delta$ will be determined
soon. Let ${\cal N}$ be a $\delta$-net of $B_X(0,g)$, and define
$\phi_H(v)$ to be the closest point in ${\cal N}$ to $f_H(v)$. We
claim that if $H_1,H_2\subseteq E$ are inequivalent, i.e., $H_1
\cap E_A\neq H_2 \cap E_A$, then $\phi_{H_1}\neq \phi_{H_2}$.
Indeed, we may assume that there are $u,v\in A$ such that
$(u,v)\in H_1\setminus H_2$. Since the girth of $G$ is at least
$g$, this implies that $\rho_{H_2}(u,v)= g-1$, whereas
$\rho_{H_1}(u,v)= 1$. Now, if $\phi_{H_1}(u)=\phi_{H_2}(u)$ and
$\phi_{H_1}(v)=\phi_{H_2}(v)$ then:
\begin{eqnarray*}
\frac{g-1}{\alpha}&=&\frac{\rho_{H_2}(u,v)}{\alpha}\\
&\le& \|f_{H_2}(u)-f_{H_2}(v)\|_X \le  
2\delta+\|\phi_{H_2}(u)-\phi_{H_2}(v)\|_X\\
&=&2\delta+\|\phi_{H_1}(u)-\phi_{H_1}(v)\|_X\le  
4\delta+\|f_{H_1}(u)-f_{H_1}(v)\|_X\\
&\le& 4\delta+\rho_{H_1}(u,v)=4\delta+1.
\end{eqnarray*}
We select $\delta=\frac{g-\alpha-1}{5\alpha}$ so that this becomes a
contradiction. It follows that each of the aforementioned
$2^{e_A}/\tbinom{n}{k}$ inequivalent sets $H \in \cal H$ gives
rise to a distinct function $\phi_H: A \to \N$.

By standard volume estimates, $|{\cal N}|\le
(\frac{2g}{\delta})^h$. Hence there are at most $|{\cal N}|^k\le
(\frac{2g}{\delta})^{kh}$ distinct functions from $A$ to ${\cal
N}$. Consequently,
$$
\left(\frac{2g}{\delta}\right)^{kh} \ge
\frac{2^{e_A}}{\tbinom{n}{k}} \ge
\frac{2^{q\tbinom{k}{2}}}{\tbinom{n}{k}}.
$$
By estimating $\binom{n}{k}\leq (\frac{ne}{k})^k \le
(\frac{ne}{s})^k$ we have that
$$
h\log_2 \left(\frac{2g}{\delta}\right) \ge \frac{(k-1)q}{2} -
\log_2\left(\frac{en}{s}\right),
$$
which yields the claimed bound on $k$.
\Endproof\vskip4pt 
Such graphs do exist as we now show:
\begin{lemma}\label{lemma:exist_high_girth}
For every integer $g \ge 4${\rm ,} there exist graphs $G=(V,E)$ of
arbitrarily large order $n$ and girth at least $g$ in which every
set $A \subseteq V$ of cardinality at least $n^{1-\frac{1}{8g}}$
has density at least $n^{-1+\frac{1}{2g}}$.
\end{lemma}

\Proof 
This is a standard construction from random graph theory. Let $N
\ge C^g$ be an arbitrarily large integer, where $C$ is an
appropriately chosen constant. Let $\eta=\frac{1}{4g}$. Pick a
random graph in $G(N,p)$ where $p=2\cdot N^{-1+2\eta}$. We claim: (i) With probability $\ge
\frac{1}{2}$ this graph has fewer than $\frac{N}{2}$ cycles of length $< g$, and (ii) With almost
certainty, every set of cardinality $\ge
(\frac{N}{2})^{1-\frac{\eta}{2}}$ has density $\ge
(\frac{N}{2})^{-1+2\eta}$. The theorem now follows by taking a
graph with these two properties and removing $\frac{N}{2}$
vertices, including at least one vertex from each cycle of length
$< g$. The resulting graph has $\frac{N}{2}$ vertices, it has no
short cycles, and satisfies the density condition.

The expected number of cycles of length $<g$ is
$$
\sum_{i=3}^{g-1} \frac{1}{2i} p^{i} N(N-1)\dots (N-i+1)
\leq\tfrac{1}{6} \sum_{i=3}^{g-1}(pN)^{i} \leq (pN)^{g-1} =
\left(2\cdot N^{\frac{1}{2g}}\right)^{g-1} \leq \frac{N}{4}.
$$
In the last inequalities we use the facts that $pN\geq 2$, $N \ge
C^{g}$, and $C\geq 4$. It follows that with probability $\ge
\frac{1}{2}$, there are no more than $N/2$ cycles shorter than
$g$.

The expected density in every set of vertices is, of course, $p$.
To estimate the deviation, we use the Chernoff bound:
$$\Pr[e_A \leq \tfrac{1}{2} \tbinom{|A|}{2}p]
\leq e^{-\tbinom{|A|}{2} \frac{p}{8}}.$$ Thus, the probability
that there exists a set of cardinality $\ge k$ and density $\le
p/2$ does not exceed $2^N \exp(-\tbinom{k}{2}p/8)$. For $k =
(\tfrac{N}{2})^{1-\frac{\eta}{2}}$, $p=2\cdot N^{-1+2\eta}$, and the assumption  $N\geq C^g$, this
is  easily
seen to be $o(1)$. The claim follows.
\Endproof\vskip4pt 

Theorem~\ref{thm:lowdim} now follows easily:

\demo{Proof of Theorem~{\rm \ref{thm:lowdim}}}
The claim for $\alpha \ge 2$ is obtained combining
Lemma~\ref{lemma:girth_and_density} and
Lemma~\ref{lemma:exist_high_girth} with $g=\lceil \alpha+2
\rceil$. As $1/q \le s \le n^{1-\frac{1}{8\alpha}}$, we obtain
that for an appropriate choice of constant $C$, $$R_X(\alpha,n)
\leq Cn^{1-\frac{1}{8\alpha}} (h\log\alpha +\log
(n^{\frac{1}{8\alpha}})).
$$
 We choose $c=1/16$, so that $
n^{1-c/\alpha} = n^{1-\frac{1}{16\alpha}} \geq
n^{1-\frac{1}{8\alpha}} \log (n^{\frac{1}{8\alpha}})$. The claim
for $\alpha \geq 2$ now follows.

For the case $\alpha < 2$ we use, instead of
Lemma~\ref{lemma:exist_high_girth}, the (trivial) analogous
statement for the complete graph $K_n$.  {That is}, we apply
Lemma~\ref{lemma:girth_and_density}  with $g=3$, $s=2$, and $q=1$.
\Endproof\vskip4pt

We are now ready to prove the promised upper bounds on
$R_2(\alpha,n)$.

\demo{Proof of Corollary~{\rm \ref{corol:basic bounds}}}
The result follows from the Johnson-Lindenstrauss dimension
reduction lemma~\cite{jl} for $\ell_2$. Let $\alpha\geq 1$, and
let $k =R_2(\alpha,n)$; i.e., every $n$-point metric space $M$
contains a $k$-point subspace that $\alpha$-embeds into $\ell_2$.
By~\cite{jl}, for any $0 < \delta \leq 1$, this subspace
$\alpha(1+\delta)$-embeds into $\ell_2^h$ with $h\leq \frac{C \log
k}{\delta^2}$. Hence, $R_2(\alpha,n) \leq R_{\ell_2^h}
(\alpha(1+\delta),n)$. Now apply Theorem~\ref{thm:lowdim}. The
claim for $\alpha \ge 2$ follows by taking $\delta=1$. For
$\alpha=2-\epsilon$ we set $\delta = \epsilon/4$. Then
$$ k = R_2(2-\epsilon,n) \leq  
R_{\ell_2^h}\left(2-\frac{\epsilon}{2},n\right) \leq
  2\log_2 n + \frac{C \log k \log
\left(\tfrac{2}{\epsilon}\right)}{\epsilon^2}  ,$$ which implies
the bound in the proposition.
\Endproof\vskip4pt 

Another interesting consequence of Theorem~\ref{thm:lowdim} are
upper bounds for metrics defined by planar graphs. This may be
interesting in view of the fact that the target metrics in our
lower bounds in Section~\ref{section:lower bounds} are
ultrametrics (and thus planar).

\begin{theorem}\label{thm:planar}
Let ${\mathcal F}$ be a family of graphs{\rm ,} none of which contains a
fixed minor $H$  on $r$ vertices. Then for every integer $n$ and
every $ \alpha \geq 1$\/{\rm :}\/
$$
R_{\cal F}(\alpha,n)\le C r^3 n^{1-c/\alpha} \log^2 n\log\alpha,
$$
where $c,C>0$ are universal constants.
\end{theorem}

\Proof  It is implicit in \cite{rao} that the
Euclidean embedding that Rao constructed is also a good low
dimensional embedding into $\ell_\infty$. More precisely, if $F\in
{\mathcal F}$ is a graph on $n$ points then it embeds with
distortion $C$ into $\ell_\infty^h$, with $h\le C r^3\log^2 n$.
Thus $R_{\cal F}(\alpha,n) \leq R_{\ell_\infty^h} (C \alpha, n)$.
The result now follows from Theorem~\ref{thm:lowdim}.
\Endproof\vskip4pt

\section{Expanders and Poincar\'e inequalities}\label{section:expanders}

In this section we prove lower bounds for the metric Ramsey
function in the case of expanders. The proof is based on
generalizations of Poincar\'e inequalities used by Matous\'{e}k to
prove lower bounds on the Euclidean distortion of expander graphs.
To obtain these inequalities we pass to a power of the graph and
  delete vertices with small degree. The argument shows that large
subsets of expanders contain large sub-subsets which satisfy an
appropriate Poincar\'e inequality (see Lemma \ref{lem:newpoin}
below). First, we recall some basic concepts on graphs.

Let $G=(V,E)$ be a $d$-regular graph, and let $A$ be its adjacency
matrix, i.e.\ $A_{uv}=1$ if $[u,v]\in E$ and $A_{uv}=0$ otherwise.
Let $\lambda_1 \ge \lambda_2 \ge \dots  \ge \lambda_n$ be the
eigenvalues of $A$. It is easy to observe that $\lambda_1=d$.
Also, $\mathrm{trace}(A)=0$, so that  $\lambda_n < 0$. We occasionally
write $\lambda_i(G)$ to specify that $G$ is the graph under
consideration. We define $G$'s multiplicative spectral gap as:
$$
\gamma(G)=\frac{\lambda_2(G)}{d}.
$$

We also define the absolute multiplicative spectral gap of $G$ as:
$$
\gamma_+(G)=\frac{\max_{i\ge
2}|\lambda_i(G)|}{d}=\frac{\max\{\lambda_2(G),-\lambda_n(G)\}}{d}.
$$

In what follows we will use the following standard estimate:
$\gamma_+(G)\ge 1/d$. To verify it observe that
$nd=\mathrm{trace}(A^2)=\sum_{i=1}^n\lambda_i(G)^2\le
d^2+(n-1)[d\gamma_+(G)]^2$ and use the fact that $1\le d\le n-1$.
We remark that this elementary bound is weaker than the
Alon-Boppana bound~\cite{alon}, but it is sufficient for our
purposes, and holds for all $d$ (while the Alon-Boppana bound only
holds for small enough~$d$).

The main statement of this section is:
\begin{theorem}\label{thm:exp} Let $G=(V,E)$ be a $d$\/{\rm -}\/regular graph{\rm ,}  
$d\ge 3$.
Let $\gamma = \gamma_+(G)$. Then for every $p,\alpha\ge 1$\/{\rm :}\/
$$ R_p(G;\alpha) \leq Cd|V|^{1-c\frac{\log_d(1/\gamma)}{p \alpha}}, $$
where $C,c$ are absolute constants.
\end{theorem}

Given $S,T\subseteq V$, we denote by $E(S,T)$ the set of directed
edges between vertices in $S$ and $T$; i.e,
$$
E(S,T)=\{(u,v)\in S\times T; [u,v]\in E\}.
$$
We also denote by $E(S)$ the set of edges in the subgraph induced
by $G$ on $S$; i.e.,
 $$
E(S)=\left\{\{u,v\};\ u,v\in S,\ [u,v]\in E\right\}.
$$
With this notation, $|E(S)|=\frac{|E(S,S)|}{2}$.

The ``Expander Mixing Lemma'' \cite{alon} states:
\begin{lemma}[Expander mixing lemma] Let $G=(V,E)$ be a $d$-regular  
graph \/{\rm (}\/which may have loops and\/{\rm /}\/or parallel edges\/{\rm ).}\/
Then for every $S,T\subseteq V${\rm ,}
$$
\left||E(S,T)|-\frac{d|S||T|}{|V|}\right|\le
\gamma_+(G)d\sqrt{|S||T|}.
$$
In particular{\rm ,}
$$
\left|\frac{2|E(S)|}{|S|}-\frac{|S|}{|V|}d\right|\le \gamma_+(G)d.
$$
\end{lemma}

\begin{lemma} \label{lem:expander-large-subset}
Let $G=(V,E)$ be a $d$\/{\rm -}\/regular graph{\rm ,} $d\ge 3$. Let $\gamma =
\gamma_+(G)$. Then for any $B\subset V$ satisfying $|B|\geq 8
\gamma |V|${\rm ,} there exists $C\subset B$ such that $|C|\geq |B|/3${\rm ,}
and for any $u\in C${\rm ,}
$$ d \frac{|B|}{8|V|} \leq \mathrm{deg}_C(u) \leq d \frac{4|B|}{|V|}.$$
\end{lemma}
\vskip8pt

\Proof 
Denote $k=|B|$. By the expander mixing lemma,
$$ |E(B)| \le
\frac{d k^2}{2n}\left(1+\frac{1}{8}\right) \le \frac{d k^2}{n}.
$$
Set $B'=\left\{v\in B;\ \mathrm{deg}_B(v)\le (4d k)/n\right\}$.
Since the graph induced by $G$ on $B$ contains $k-|B'|$ vertices
of degree greater than $(4dk)/n$, it follows that:
$$
\frac{d k^2}{n}\ge |E(B)|\ge \frac{k-|B'|}{2}\cdot\frac{4d k}{n},
$$
so that $|B'|\ge k/2$. Again, by the expander mixing lemma,
$$
\frac{2|E(B')|}{k} \geq \frac{2|E(B')|}{2|B'|}  \geq
\frac{dk}{4n}.
$$

We now apply an iterative procedure which produces a sequence
$B'=B_0\supset B_1\supset B_2\supset \dots $ as follows: if
$\min_{v\in B_i}\mathrm{deg}_{B_i}(v)\le \frac{d k}{8n}$ then
$B_{i+1}$ is obtained\break\vglue-12pt\noindent from $B_i$ by throwing away a vertex $u\in
B_i$ with $\mathrm{deg}_ {B_i}(u)=\min_{v\in
B_i}\mathrm{deg}_{B_i}(v)$. Otherwise $B_{i+1}=B_{i}$. This
procedure eventually ends, and we are left with a subset
$C\subseteq B'$. Since at each step we delete at most $\frac{d
k}{8n}$ edges from $B'$, we have that:
$$
|E(C)|\ge |E(B')|-\frac{d k^2}{8n}\ge \frac{d k^2}{4n}-\frac{d
  k^2}{8n}=\frac{d k^2}{8n}.
$$

Note that the graph induced on $C$ by $G$ has minimal degree at
least $\frac{d k}{8n}$. To estimate $|C|$ we apply the expander
mixing lemma to get that:
$$
\frac{2|E(C)|}{|C|}\le d\left(\frac{|C|}{n}+\gamma\right).
$$
Thus
$$
\frac{d|C|^2}{2n}+\frac{\gamma d |C|}{2}\ge |E(C)|\ge \frac{d
k^2}{8n}.
$$
Since $k\ge 8\gamma n$, $\frac{\gamma d |C|}{2}\le \frac{d
k^2}{16n}$. Hence:
$$
\frac{d |C|^2}{2n}\ge \frac{d k^2}{16n},
$$
so that $|C|\ge \frac{k}{3}$.
\Endproof\vskip4pt 

The following is the Poincar\'e inequality used in the proof of
Theorem~\ref{thm:exp}.

\begin{lemma}\label{lem:newpoin} Let $G=(V,E)$ be a $d$\/{\rm -}\/regular graph{\rm ,}  
$d\ge 3$. Let $\gamma =
\gamma_+(G)$. Then for any $B\subset V$ satisfying $|B|\geq 8
\gamma |V|${\rm ,} there exists $C\subset B$ such that $|C|\!\geq\! |B|/3$
and the following holds true\/{\rm :}\/ For any $p\!\ge\! 1$ and for every
$f\!:\! C\!\to\! \ell_p$\/{\rm :}\/
$$
\sum_{u,v\in C}\|f(u)-f(v)\|_p^p\le \frac{(32p)^p |V|}{d}
\sum_{[u,v]\in E(C)}\|f(u)-f(v)\|_p^p.
$$
\end{lemma}
\vskip8pt

The proof of Lemma~\ref{lem:newpoin} proceeds by first proving a
slightly stronger version of it for $p=2$ and then extrapolating
to the general case via the following lemma based on an
extrapolation argument which was used in \cite{matexpander}.  Its
proof   is delayed to the end of the section.

\begin{lemma}[Extrapolation lemma for Poincar\'e  
inequalities]\label{lem:extrapolation} Let
$G=\break (V,E)$ be a graph with maximal degree at most $\Delta$. Fix
$p\ge 1$ and let $A>0$ be a constant such that for every $f:V\to
\R${\rm ,}
\begin{eqnarray}\label{eq:poin}
\sum_{u,v\in V}|f(u)-f(v)|^p\le (A p)^p \frac{|V|}{\Delta}
\sum_{[u,v]\in E}|f(u)-f(v)|^p.
\end{eqnarray}
Then for every $0<q\le p$ and for every $f:V\to \ell_q${\rm ,}
$$
\sum_{u,v\in V}\|f(u)-f(v)\|_q^q\le (A p)^p \frac{|V|}{\Delta}
\sum_{[u,v]\in E}\|f(u)-f(v)\|_q^q.
$$
Additionally{\rm ,} for every $p<q<\infty$ and every $f:V\to \ell_q$\/{\rm :}\/
$$
\sum_{u,v\in V}\|f(u)-f(v)\|_q^q\le \left(4A q\right)^q
\frac{|V|}{\Delta} \sum_{[u,v]\in E}\|f(u)-f(v)\|_q^q.
$$
\end{lemma}
\vskip8pt

\demo{Proof of Lemma~{\rm \ref{lem:newpoin}}}
Denote $n=|V|$ and $k=|B|$.
  By Lemma~\ref{lem:expander-large-subset}, there exists $C\subset B$
with $|C|\geq k/3$, such that the induced subgraph of $G$ on $C$,
has minimal degree at least $kd/8n$ and maximal degree at most
$\Delta = 4kd/n$.

  We first prove that the following inequality hold true for every  
$f:V\to \ell_2$:
\begin{eqnarray}\label{eq:l2-poin}
  \sum_{u,v\in C}\|f(u)-f(v)\|_2^2&\le& \frac{32n}{d} \sum_{[u,v]\in  
E(C)}\|f(u)-f(v)\|_2^2\\
& =& \frac{32\cdot 4
k}{\Delta} \sum_{[u,v]\in E(C)}\|f(u)-f(v)\|_2^2.\nonumber
\end{eqnarray}

By summation we may clearly assume that $f:C\to \R$. By
translation we may assume that $\sum_{v\in C} f(v)=0$. Extend $f$
to $V$ by letting $f(u)=0$ for $u\notin C$. Now,
\begin{align*}
\sum_{u,v\in C}[f(u)-f(v)]^2
&=2(|C|+1)\sum_{v\in C}f(v)^2-2\sum_{u,v\in C}f(u)f(v)\\
&= 2(|C|+1)\sum_{v\in V}f(v)^2-\biggl(\sum_{v\in C}f(v)\biggr)^2\\
&=2(|C|+1)\sum_{v\in V}f(v)^2\le 4 k\sum_{v\in V}f(v)^2.
\end{align*} 

Since $\sum_{v\in V}f(v)=0$ we can use the spectral gap of $H$ to
get that: 
\begin{eqnarray*}
\sum_{[u,v]\in E(C)}[f(u)-f(v)]^2 &=&2\sum_{v\in
C}\mathrm{deg}_{E(C)}(v)f(v)^2-
2\sum_{[u,v]\in E(C)}f(u)f(v)\\
&\ge& 2\sum_{v\in V}\frac{d k}{8n}f(v)^2-2\sum_{[u,v]\in E}f(u)f(v)\\
&=&
\frac{d k}{4n}\sum_{v\in V}f(v)^2-\langle A^t f, f\rangle\\
&\ge& \frac{d k}{4n}\sum_{v\in V}f(v)^2-\gamma d \sum_{v\in V}f(v)^2\\
&\ge& \left(\frac{k}{4n}-\gamma \right)d \sum_{v\in V}f(v)^2\ge
\frac{k}{8n}d \sum_{v\in V}f(v)^2,
\end{eqnarray*} 
which implies inequality~\eqref{eq:l2-poin}. The Poincar\'e
inequalities for $p \ge 1$ now follow immediately from
inequality~\eqref{eq:l2-poin} via Lemma~\ref{lem:extrapolation},
and by substituting the value of $\Delta$.
\Endproof\vskip4pt

{\it Proof of Theorem} \ref{thm:exp}.
The proof proceeds by showing that for every $B\subseteq V$
satisfying $c_p(B)\le \alpha$,
$$
|B|\le 100d|V|^{1-\frac{\log_d(1/\gamma)} {2561p \alpha}}.
$$

Set $k=|B|$ and $n=|V|$. Define:
$$
t=\left\lfloor\frac{\log(8n)}{2560p \alpha\log
d+\log\left(\frac{1}{\gamma}\right)}\right\rfloor.
$$
Note that $t\le \mathrm{diam}(G)$ since it is well known that
$\mathrm{diam}(G)\ge \log_d(n)$. We may also assume that $t\ge 1$,
since otherwise, using the fact that $1/\gamma\le d$, we get that
$n<e^{2561p\alpha\log d}$, in which case the required result holds
vacuously.

Denote by $A$ the adjacency matrix of $G$. Let $H$ be the
multi-graph with adjacency matrix $A^t$. In other words, the
number of $H$-edges between two vertices in $V$ is the number of
distinct paths of length $t$ joining them (and it is $0$ if no
such path exists). The multi-graph $H$ is $d^t$ regular and by the
spectral theorem, $\gamma_+(H)=\gamma^t$. We may assume that $k\ge
8\gamma^t n$ (otherwise the conclusion of the theorem in trivial).

  It follows from
Lemma~\ref{lem:newpoin} that there exists $C\subset B$ such that
$|C|\geq k/3$, and for every $f:C\to \ell_p$:
$$
\sum_{u,v\in C}\|f(u)-f(v)\|_p^p\le \frac{(32p)^p n}{d^t}
\sum_{[u,v]\in E_H(C)}\|f(u)-f(v)\|_p^p.
$$

Let $f:B\to \ell_p$ be an embedding such that for all $u,v\in B$,
$\frac{d_G(u,v)}{\alpha}\le \|f(u)-f(v)\|_p\le d_G(u,v)$. Then:
$$
\sum_{[u,v]\in E_H(C)}\|f(u)-f(v)\|_p^p\le \sum_{[u,v]\in
E_H(B)}d_G(u,v)^p\le |E_H(B)|t^p\le \frac{d^t k^2 t^p}{n},
$$
where the last inequality follows from an application of the
expander mixing lemma:
$$ |E_H(B)| \le
\frac{d^t k^2}{2n}\left(1+\frac{1}{8}\right) \le \frac{d^t
k^2}{n}.
$$

Let $s=\left\lfloor\log_d\left(\frac{k}{12}\right)\right\rfloor$.
We may clearly assume that $s>1$. The number of vertices of
distance at most $s$ from a given vertex $v_0\in G$ is bounded by:
$$
1+d+\dots +d^s\le 2d^s\le \frac{k}{6}\le \frac{|C|}{2}.
$$
Hence:
\begin{eqnarray} \label{eq:exp1}
  \sum_{u,v\in C}\|f(u)-f(v)\|_p^p&\ge& \frac{1}{\alpha^p}\sum_{u,v\in
C}d_G(u,v)^p\ge \frac{|C|^2s^p}{2\alpha^p}\\
&\ge&
\frac{k^2s^p}{18\alpha^p}\ge
\frac{k^2}{80\alpha^p}\left[\log_d\left(\frac{k}{12}\right)\right]^p.\nonumber
\end{eqnarray}
Plugging this into the Poincar\'e inequality we get that:
$$
\frac{k^2}{80\alpha^p}\left[\log_d\left(\frac{k}{12}\right)\right]^p\le
\frac{(32p)^p n}{d^t}\cdot\frac{d^t k^2 t^p}{n},
$$
which gives
$$
\log_d\left(\frac{k}{12}\right)\le 2560 p \alpha t\Longrightarrow
k\le 12 d^{2560 p \alpha  t}.
$$
Since $t\le \frac{\log (8n)}{2560 p \alpha\log
d+\log\left(\frac{1}{\gamma}\right)}$, it follows that:
\begin{eqnarray*}
k&\le& 12 \exp\left[\frac{\log (8n)\cdot 2560 p \alpha\log d}{2560 p
\alpha\log d+\log\left(\frac{1}{\gamma}\right)} \right]\\
&\le&
100n^{1-\frac{\log(1/\gamma)}{2560 p \alpha\log d
+\log(1/\gamma)}}\le 100n^{1-\frac{\log_d(1/\gamma)} {2561 p
\alpha}},
\end{eqnarray*}
where we have used once more the estimate $\log(1/\gamma)\le \log
d$.
\Endproof\vskip4pt

It remains to prove Lemma~\ref{lem:extrapolation}.

\demo{Proof of Lemma~{\rm \ref{lem:extrapolation}}} The case $0<q\le  
p$ is
simple. Coordinate-wise summation of (\ref{eq:poin}) shows that
for every $f:V\to \ell_p$:
$$
\sum_{u,v\in V}\|f(u)-f(v)\|_p^p\le (A p)^p \frac{|V|}{\Delta}
\sum_{[u,v]\in E}\|f(u)-f(v)\|_p^p.
$$
Since $\ell_2$ is isometric to a subspace of $L_p$, it follows
that for every $f:V\to \ell_2$,
$$
\sum_{u,v\in V}\|f(u)-f(v)\|_2^p\le (A p)^p \frac{|V|}{\Delta}
\sum_{[u,v]\in E}\|f(u)-f(v)\|_2^p.
$$
Since $(\R, |x-y|^{q/p})$ is isometric to a subset of $\ell_2$
(\cite{ww}, \cite{dezalau}), the required inequality follows.

We now pass to the case $p<q$. In this case the following standard
numerical inequality holds true for every $a,b\in \R$ (see Lemma 4
in \cite{matexpander}):
\begin{equation}\label{eq:numerical}
\big| |a|^{q/p}\mathrm{sign}(a)-|b|^{q/p}\mathrm{sign}(b)\big|\le
\frac{q}{p}|a-b|\left(
|a|^{\frac{q}{p}-1}+|b|^{\frac{q}{p}-1}\right).
\end{equation}

It is suffices to prove the claims coordinate-wise, i.e.\ for
functions $f\!:\!V\!\to\! \mathbb{R}$. Fix some $f\!:\!V\!\to\! \R$. By continuity
there is some $c\in \R$ such that:
$$
\sum_{v\in V}|f(v)+c|^{q/p}\mathrm{sign}(f(v)+c)=0.
$$
Hence, by replacing $f$ with $f+c$ we may assume that:
$$
\sum_{v\in V}|f(v)|^{q/p}\mathrm{sign}(f(v))=0.
$$
Now: 
\begin{eqnarray*}
\sum_{v\in V} |f(v)|^q&=& \sum_{v\in
V}\left||f(v)|^{q/p}\mathrm{sign}(f(v))-\frac{1}{|V|}\sum_{u\in
V}|f(u)|^{q/p}\mathrm{sign}(f(u))\right|^p\\
&\le& \frac{1}{|V|}\sum_{u,v\in V}\big|
|f(u)|^{q/p}\mathrm{sign}(f(u))-|f(v)|^{q/p}\mathrm{sign}(f(v))\big|^p\\
&\le& \frac{(A p)^p}{\Delta} \sum_{[u,v]\in E}\big|
|f(u)|^{q/p}\mathrm{sign}(f(u))-|f(v)|^{q/p}\mathrm{sign}(f(v))\big|^p\\
&\le& \frac{(A q)^p}{\Delta}\sum_{[u,v]\in E}
|f(u)-f(v)|^p\left(|f(u)|^{\frac{q}{p}-1}+|f(v)|^{\frac{q}{p} 
-1}\right)^p,
\end{eqnarray*} 
where in the last two steps we have used (\ref{eq:poin}) and
(\ref{eq:numerical}), respectively. An application of H\"older's
inequality gives that:
\begin{multline*}
\sum_{[u,v]\in E}
|f(u)-f(v)|^p \cdot \left(|f(u)|^{\frac{q}{p}-1}+|f(v)|^{\frac{q}{p} 
-1}\right)^p\\
 \le  \left(\sum_{[u,v]\in
E}|f(u)-f(v)|^q\right)^{p/q}\left(\sum_{[u,v]\in
E}\left(|f(u)|^{\frac{q}{p}-1}+|f(v)|^{\frac{q}{p} 
-1}\right)^{\frac{pq}{q-p}}\right)^{1-\frac{p}{q}}.
\end{multline*} 
Using the assumption on the maximal degree we get that:
\begin{eqnarray*}
\sum_{[u,v]\in E} \left(
|f(u)|^{\frac{q}{p}-1}+|f(v)|^{\frac{q}{p}-1}\right)^{\frac{pq}{q-p}}
&\le& 2^{\frac{qp}{q-p}-1}\sum_{[u,v]\in E}(|f(u)|^q+|f(v)|^q)\\
& \le&
2^{\frac{qp}{q-p}}\Delta \sum_{v\in V}|f(v)|^q.
\end{eqnarray*} 
Summarizing, we have shown that:
$$
\sum_{v\in V} |f(v)|^q\le \frac{(2A q)^p}{\Delta^{\frac{p}{q}}}
\left(\sum_{[u,v]\in E}|f(u)-f(v)|^q\right)^{p/q}\left(\sum_{v\in
V}|f(v)|^q\right)^{1-\frac{p}{q}}.
$$
This inequality simplifies to:
$$
\sum_{v\in V} |f(v)|^q \le \frac{(2A q)^q}{\Delta} \sum_{[u,v]\in
E}|f(u)-f(v)|^q.
$$
We conclude by noting that:
\vskip14pt
\hfill $
\displaystyle{\sum_{u,v\in V}|f(u)-f(v)|^q\le 2^q|V|\sum_{v\in V}|f(v)|^q.}
$ 
\Endproof\pagebreak

We now show that the interplay between the Euclidean distortion
and the cardinality in Theorem~\ref{thm:exp}, for $p=2$, is tight,
up to the dependence on $d$ and~$\gamma$. We require an upper
estimate for the diameter of an $n$-point expander. It is well
known that the diameter is $O(\log n)$, but here we will be a
little bit more accurate.

We need the following bound on the diameter of expander graphs
\cite{chung}:

\begin{proposition}
Let $G=(V,E)$ be an $n$\/{\rm -}\/vertex{\rm ,} $d$\/{\rm -}\/regular graph. Denote
$\gamma=\gamma_+(G)$. Then the diameter of $G$ is at most
$\log_{1/\gamma}n+1$.
\end{proposition}

\begin{proposition} \label{prop:uniform-expander}
Let $G=(V,E)$ be an $n$\/{\rm -}\/vertex{\rm ,} $d$\/{\rm -}\/regular graph{\rm ,} $d\ge 3${\rm ,} and
set $\gamma=\gamma_+(G)$. Then{\rm ,} there is an absolute constant
$C>0$ such that for any $\alpha>1${\rm ,}
$$ R_2(G;\alpha) \ge R_{\eq}(G;\alpha) \ge
n^{1- \frac{C}{\alpha\log_{d}(1/\gamma)}}.$$
\end{proposition}

\vskip8pt

\Proof  Iteratively, extract a point $x\in V$
together with a ball of radius $r=\diam(G)/\alpha$ around $x$.
Each such ball contains at most $d+d(d-1)+\dots +d(d-1)^{\lfloor
r\rfloor}\le 3(d-1)^{r+1}$ points, and thus we can repeat this
process at least $n/(3(d-1)^{r+1})$ times, and get the desired set.
Its size is at least
\begin{eqnarray*} \frac{n}{3(d-1)^{r+1}}&=  &
\frac{n}{3}(d-1)^{-\left(\frac{\diam(G)}{\alpha}+1\right)} \\
&\geq&
   \frac{n}{3(d-1)^2}(d-1)^{-\frac{\log_{1/\gamma}n}{\alpha}}=
   \frac{1}{3(d-1)^2}n^{1- \frac{1}{\alpha\log_{(d-1)}(1/\gamma)}}.
\end{eqnarray*}
\vglue-23pt
\Endproof\vskip4pt

\section{Markov type, girth and hypercubes}
\label{section:markov}

Markov type was defined in~\cite{ball} and was applied
in~\cite{lmn} to obtain lower bounds for the Euclidean distortion
of regular graphs with large girth. This concept plays a key role
in our analysis of the metric Ramsey problem for the discrete cube
and graphs with large girth. Let $(X,d)$ be a metric space. We
shall say that $\{M_k\}_{k=0}^{\infty}$ is a stationary
time-reversible Markov chain on $X$ if there are
$x_1,\dots ,x_n\in X$, an $n\times n$ stochastic matrix $A$ and a
stationary distribution $\pi=(\pi_1,\dots ,\pi_n)$ of $A$ such
that for every $i,j$, $\pi_i A_{ij}=\pi_j A_{ji}$,
$\{M_k\}_{k=0}^\infty$ is a Markov chain with transition matrix
$A$ and $M_0$ is distributed according to $\pi$. $(X,d)$ is said
to have Markov type $p> 0$ with constant $C$ if for any stationary
time-reversible Markov chain on $X$, and for any time $s$:
  $$ \mathbb{E}[d(Z_s,Z_0)^p] \leq C^p s \mathbb{E}[d(Z_1,Z_0)^p]
  .$$

In~\cite{ball} (see also~\cite{lmn}) it was shown that Hilbert
space has Markov type $2$ with constant $1$. Actually, these
references deal with the special case in which $A$ is symmetric
and $\pi$ is the uniform distribution on the states
$x_1,\dots ,x_n$, but the proof is easily seen to carry over to
stationary time-reversible Markov chains.

\Subsec{Graphs with large girth}\label{section:girth}
For later applications, it will be convenient to introduce a
notion of ``Euclidean distortion at small distances'' as follows.
Let $(X,d_X), (Y,d_Y)$ be metric spaces and $s>0$. For every
injective $f:X\to Y$ define:
$$
\mathrm{dist}^{(s)}(f)=\left(\sup_{0<d_X(x,y)\le
s}\frac{d_Y(f(x),f(y))}{d_X(x,y)}\right)\cdot
\left(\sup_{0<d_X(x,y)\le
s}\frac{d_X(x,y)}{d_Y(f(x),f(y))}\right),
$$
and:
$$
c^{(s)}_Y(X)=\inf\left\{\mathrm{dist}^{(s)}(f);\ f:X\to Y\right\}.
$$
As before, we write $c_2^{(s)}(X)=c_{\ell_2}^{(s)}(X)$.
\vskip2pt

Let $G=(V,E)$ be a graph. In what follows we denote by $\delta(G)$
the average degree of $G$, i.e.\ $$
\delta(G)=\frac{\sum_{v\in
V}\mathrm{deg}(v)}{|V|}=\frac{2|E|}{|V|}.
$$

We begin with the following strengthening of a result from
\cite{lmn}.
\begin{theorem}\label{th:average}
Let $G=(V,E)$ be a graph with girth $g$ and average degree
$\delta=\delta(G)${\rm ;} then for every integer $1<s<g/2${\rm ,}
$c^{(s)}_2(G)\geq \frac{\delta-2}{\delta}\sqrt{s}$. In particular{\rm ,}
$$
c_2(G)\ge\frac{\delta-2}{\delta}\sqrt{\left\lfloor
\frac{g}{2}\right\rfloor -1}.
$$
\end{theorem}
\vskip8pt

\Proof 
  Assume first that $G$ is connected. Consider the reversible Markov  
chain
$\{Z_k\}_{k=0}^\infty$ that corresponds to the canonical random
walk on $G$. Recall that $\pi_v=\deg(v)/(\delta n)$ is a
stationary distribution of this Markov chain.

For every $1<s<g/2$,
\begin{align*}
\E[d_G(Z_s,Z_0)]&\ge  \E_{v\in V} \left[ \frac{\deg(v)-1}{\deg(v)}
(\E[d_G(Z_{s-1},Z_0)| Z_{s-1}=v]+1)   \right.\\
  & \qquad \qquad
\left. + \frac{1}{\deg(v)}(\E[d_G(Z_{s-1},Z_0)| Z_{s-1}=v]-1) \right] \\
  &=
  \E_{v\in V} \left[ \frac{\deg(v)-2}{\deg(v)} + \E[d_G(Z_{s-1},Z_0)|  
Z_{s-1}=v]
   \right] \\
   &= \E[d_G(Z_{s-1},Z_0)] +1 - \sum_v \pi_v \frac{2}{\deg(v)}\\
   &= \E[d_G(Z_{s-1},Z_0)] +1 - \sum_v \frac{\deg(v)}{\delta n}  
\frac{2}{\deg(v)}\\
   &= \E[d_G(Z_{s-1},Z_0)] +\frac{\delta-2}{\delta}.
  \end{align*}
By induction $\E[d_G(Z_s,Z_0)] \geq s \frac{\delta-2}{\delta}$.
Therefore
$$ \E[d_G(Z_s,Z_0)^2] \geq [\E d_G(Z_s,Z_0)]^2 \geq s^2
\left(\frac{\delta-2}{\delta}\right)^2.$$
  On the other hand, since Hilbert space has Markov type $2$ with
  constant $1$,
  $$ \E[d_G(Z_s,Z_0)^2]\leq c^{(s)}_2(G)^2 s\E[d_G(Z_1,Z_0)^2]=  
c_2^{(s)}(G)^2 s .$$
So $c^{(s)}_2(G)\geq \frac{\delta-2}{\delta}\sqrt{s}$.

If $G$ is disconnected, there is a connected component $C$ of $G$
in which the average degree is at least $\delta = \delta(G)$. The
theorem follows by applying the above proof to the connected graph
spanned by $C$.
\Endproof\vskip4pt 

Let $G=(V,E)$ be a $d$-regular graph, and let $A$ be its adjacency
matrix, and $\lambda_1 \ge \lambda_2 \ge \dots  \ge \lambda_n$ be
the eigenvalues of $A$. Recall (see
Section~\ref{section:expanders}) that the multiplicative spectral
gap of $G$ is $ \gamma(G)=\frac{\lambda_2(G)}{d}$ and the absolute
multiplicative spectral gap of $G$ is
$\gamma_+(G)=\frac{\max\{\lambda_2(G),-\lambda_n(G)\}}{d}.$ For $S
\subseteq V$, let $ E(S)=\left\{\{u,v\};\ u,v\in S,\ [u,v]\in
E\right\}.$ Recall that the Expander Mixing Lemma implies that
$$
\frac{2|E(S)|}{|S|}\ge d\left[\frac{|S|}{|V|}-\gamma_+(G)\right].
$$
This statement motivates the following useful definition:
\begin{definition}[self mixing parameter]
Let $G=(V,E)$ be a $d$-regular graph. The {\it self-mixing parameter} of
$G$ is defined as:
$$
\mu(G)=\max\left\{\frac{|S|}{|V|}-\frac{2|E(S)|}{d|S|}; S\subseteq
V\right\}.
$$
\end{definition}

The Expander Mixing Lemma implies that $\mu(G)\le \gamma_+(G)$. We
have in fact the following estimate:
\begin{lemma}\label{le:selfmix} Let $G=(V,E)$ be a  $d$\/{\rm -}\/regular
$n$\/{\rm -}\/vertex graph{\rm ,} let $A$ be $G$\/{\rm '}\/s adjacency matrix and let
$d=\lambda_1 \ge \dots  \ge \lambda_n$ be its eigenvalues. Then\/{\rm :}\/
$$
\mu(G)\le\frac{-\lambda_n}{d}.
$$
\end{lemma}
\vskip8pt
\Proof 
Let $w_1,\dots  ,w_n$ be an orthonormal system of eigenvectors for
$A$ with $Aw_i=\lambda_iw_i$ for $i=1,\dots ,n$. Let ${\bf 1}={\bf
1}_V $ be the all-ones vector; then $w_1 =
\frac{1}{\sqrt{|V|}}{\bf 1}$. Let ${\bf 1}_S$ be the indicator of
some subset $S \subseteq V$. Then,  
\begin{eqnarray*}
2|E(S)|&=&\bigl\langle A{\bf 1}_S,{\bf 1}_S\bigr\rangle\\
&=& \Bigl\langle A\sum_{i=1}^{n}\langle{\bf 1}_S,w_i\rangle w_i,
\sum_{i=1}^{n}\langle{\bf 1}_S,w_i\rangle w_i\Bigr\rangle\\
&=&\sum_{i=1}^{n}\langle{\bf 1}_S,w_i\rangle^2\lambda_i\\
&\ge& \langle{\bf 1}_S,w_1\rangle^2\lambda_1
+\lambda_n\sum_{i=1}^{n} \langle{\bf 1}_S,w_i\rangle^2\\
&=&\frac{|S|^2}{n}\cdot d + \lambda_n |S|.
\end{eqnarray*} 
\vglue-24pt
\Endproof\vskip12pt

\begin{lemma}\label{le:subset}
Let $G=(V,E)$ be a $d$\/{\rm -}\/regular graph with girth $g$ and put
$\mu=\mu(G)$. Fix $B\subseteq V${\rm ,} $1\le s<g/2$ and denote
$\alpha=c_2^{(s)}(B,d_G)$. Assume that $\alpha^2<s$. Then{\rm ,}
$$
|B|\le \mu
|V|+\frac{2|V|}{d\left(1-\frac{\alpha}{\sqrt{s}}\right)}.
$$
\end{lemma}
\vskip8pt

\Proof 
Set $|B|=k$ and $|V|=n$. By the definition of the self mixing
parameter,
$$
2|E(B)|\ge \frac{dk^2}{n}-\mu dk.
$$
Consider the graph on $B$ induced by $G$ (i.e.\ the edges are the
edges of $G$ which are also in $B\times B$). Its girth is not less
than $g$ and its average degree is:
$$
\delta=\delta(B)=\frac{2|E(B)|}{k}\ge \frac{dk}{n}-\mu d.
$$
Moreover, since $G$ has girth $g$ and $s<g/2$, if $d_B(u,v)\le s$,
for some two vertices $u,v\in B$, then $d_B(u,v)=d_G(u,v)$.
Consequently, $c_2^{(s)}(B,d_B)\le c_2^{(s)}(B,d_G)=\alpha$. An
application of Theorem~\ref{th:average} yields:
$$
\alpha\ge \left(1-\frac{2}{\delta}\right)\sqrt{s},
$$
so that,
$$
\frac{dk}{n}-\mu d\le \delta\le
\frac{2}{1-\frac{\alpha}{\sqrt{s}}},
$$
which gives:
\vskip12pt
\hfill $
\displaystyle{k\le \mu n +\frac{2n}{d\left(1-\frac{\alpha}{\sqrt{s}}\right)}.}
 $
\Endproof\vskip12pt

Let $G=(V,E)$ be a graph and $1\le t\le \mathrm{diam}(G)$ be an
integer. We define the $t$-distance graph of $G$ as $G^{(t)}=(V,
E^{(t)})$ where $[u,v]\in E^{(t)}$ if and only if $d_G(u,v)=t$. We
collect below some properties of $G^{(t)}$ (part ${\rm 6}$ of the
lemma below will not be applied in the sequel, and is included
here for possible future reference).

\begin{lemma}\label{le:gt} Let $G=(V,E)$ be a $d$\/{\rm -}\/regular graph{\rm ,} $d\ge  
3${\rm ,} with
girth $g$ and let $1\le t <g/2$. Then\/{\rm :}\/
\begin{itemize}
\ritem{1)} $G^{(t)}$ is a $d(d-1)^{t-1}$ regular graph.

\ritem{2)} The girth of $G^{(t)}$ is not less than $g/t$.

\ritem{3)} For every $u,v\in V${\rm ,}
$$
d_{G^{(t)}}(u,v)<\frac{g}{2t}\Longrightarrow
d_{G^{(t)}}(u,v)=\frac{d_G(u,v)}{t}.
$$

\ritem{4)} For every $B\subseteq V$ and $1\le
s<\frac{g}{2t}${\rm ,} $c_2^{(s)}(B,d_{G^{(t)}})\le c_2(B,d_G)$.

 \ritem{5)} If $t$ is even then $\mu(G^{(t)})\le
8(d-1)^{-t/4}$.

\ritem{6)} If $t$ is odd then $\gamma(G^{(t)})\le
8e^{-(1-\gamma(G))t/8}$. \end{itemize}
\end{lemma}

\Proof  Since $G$ has girth $g$ and $t<g/2$, the number of  
vertices with distance $t$
from a given vertex is the number of leaves of a $d$-regular tree
of depth $t$, which is $d(d-1)^{t-1}$. This proves ${\rm 1)}$. The
statements ${\rm 2)}$ and ${\rm 3)}$ are also simple consequences
of the fact that $G$ has girth $g$. Assertion ${\rm 4)}$ follows
immediately from assertion ${\rm 3)}$.

To prove assertion ${\rm 5)}$, note that the adjacency matrix of
$G^{(t)}$, $A^{(t)}$, is the $t$-distance matrix of $G$; i.e., $A^{(t)}_{uv}=1$ if $d_G(u,v)=1$ and $0$
otherwise. If we denote by $A$ the adjacency matrix of $G$ then there exists a polynomial
$P_t$ of degree $t$ such that $A^{(t)}=P_t(A)$ (this is the
so-called Geronimus polynomial. The properties of the
polynomials used here can be found e.g.\ in
\cite{biggs}, \cite{lmn}). The polynomial $P_t$ has degree $t$; all its
roots are real and reside in the interval
$[-2\sqrt{d-1},2\sqrt{d-1}]$. An explicit trigonometric expression
for $P_t$ is:
$$
P_t(2\sqrt{d-1}\cos \theta)=
(d-1)^{t/2-1}\frac{(d-1)\sin((t+1)\theta)-\sin((t 
-1)\theta)}{\sin\theta}.
$$
Finally, if $t$ is even, then $P_t$ is an even function and if $t$
is odd, then $P_t$ is an odd function. The spectral theorem shows
that $\{P_t(\lambda_i(G))\}$ are the eigenvalues of~$A^{(t)}$.

We turn to estimate the smallest eigenvalue of $A^{(t)}$, for $t$
even. This eigenvalue must be negative, but $P_t$ is positive
outside the interval $[-2\sqrt{d-1},2\sqrt{d-1}]$. In other words,
if $P_t(x) < 0$, then $x\in [-2\sqrt{d-1},2\sqrt{d-1}]$ and
$x=2\sqrt{d-1}\cos \theta$ for some $\theta \in [-\pi,\pi]$.
Therefore,
$$
P_t(x)=(d-1)^{t/2-1}\frac{(d-1)\sin((t+1)\theta)-\sin((t 
-1)\theta)}{\sin\theta}.
$$
Using the elementary estimate $|\sin r\alpha|\le r|\sin \alpha|$
for $\alpha\in [-\pi,\pi]$ and $r\ge 1$, it follows that:
$$
|P_t(x)|\le d(t+1)(d-1)^{t/2-1}.
$$
Hence, by Lemma~\ref{le:selfmix},
$$
\mu(G^{(t)})\le \frac{d(t+1)(d-1)^{t/2-1}}{d(d-1)^{t-1}}=
(t+1)(d-1)^{-t/2}\le 8(d-1)^{-t/4}.
$$

To prove assertion ${\rm 6)}$ we distinguish between two cases:

\demo{Case one} $\lambda_2(G^{(t)})=P_t(\lambda_2(G))$.
In this case we apply the mean value theorem and find some $a\in
(\lambda_2(G),\lambda_1(G))$ such that:
\begin{eqnarray*}
\log\left[\frac{1}{\gamma(G^{(t)})}\right]&=&\log\left[\frac{P_t(d)}{P_t 
(\lambda_2(G))}\right]\\
&=&[d-\lambda_2(G)]\frac{P_t'(a)}{P_t(a)}\\
&=&[1-\gamma(G)]d\sum_{i=1}^t\frac{1}{a-y_i}
\end{eqnarray*} 
where $y_i$ are the roots of $P_t$. Therefore,
$$
\log\left[\frac{1}{\gamma(G^{(t)})}\right] \ge
[1-\gamma(G)]t\frac{d}{d+2\sqrt{d-1}}\ge \frac{[1-\gamma(G)]t}{2},
$$
as claimed.

\demo{Case two} $\lambda_2(G^{(t)})=P_t(\lambda_i(G))$
for some $i\ge 3$. We claim that $\lambda_i$ must be in the
interval $[-2\sqrt{d-1},2\sqrt{d-1}]$. Recall that all the zeros
of $P_t$ are in this interval. It is impossible that $\lambda_i <
-2\sqrt{d-1}$, since $P_t < 0$ there ($t$ is odd). Also,
$\lambda_i > 2\sqrt{d-1}$ is impossible, since $P_t$ is increasing
on $[2\sqrt{d-1},\infty)$ and $\lambda_2(G)\ge \lambda_i(G)$,
whereas $P_t(\lambda_i(G))>P_t(\lambda_2(G))$. Therefore,
$\lambda_i(G)\in [-2\sqrt{d-1},2\sqrt{d-1}]$, and as in the proof
of ${\rm 5)}$, we estimate:
$$
P_t(\lambda_i(G))\le d(t+1)(d-1)^{t/2-1}.
$$
Hence:
$$
\gamma(G^{(t)})=\frac{P_t(\lambda_i(G))}{d(d-1)^{t-1}}\le
(t+1)(d-1)^{-t/2}\le (t+1)\cdot 2^{-t/2},
$$
which implies the required result.
\Endproof\vskip4pt 

We can now prove an upper bound for the Ramsey problem for graphs
with large girth.

\begin{theorem}\label{th:girthramsey}
Let $G=(V,E)$ be a $d$\/{\rm -}\/regular graph{\rm ,} $d\ge 3$ with girth $g$. Let
$1\le \alpha< \frac{\sqrt{g}}{6}$. There is an absolute constant
$c>0$ such that
$$R_2(G;\alpha)\leq 12(d-1)^{- c\frac{g}{\alpha^2}}|V|.$$
\end{theorem}
 \vskip8pt

\Proof 
The proof proceeds by showing that for every $B\subseteq V$ such
that $c_2(B,d_G)\le \alpha$, the following estimate holds:
$$
|B|\le 12(d-1)^{- \frac{g}{64 \alpha^2}}|V|.
$$

Let $t$ be the unique even defined by $\frac{g}{8\alpha^2}-2\le t
< \frac{g}{8\alpha^2}$. Put $s=4\alpha^2$. Now, since
$s<\frac{g}{2t}$, part ${\rm 4)}$ of Lemma~\ref{le:gt} implies
that:
$$
c_2^{(s)}(B,d_{G^{(t)}})\le c_2(B,d_G)\le \alpha.
$$
By Lemma~\ref{le:gt}, $\mathrm{girth}(G^{(t)}) \ge g/t$. Also,
$s<\frac{g}{2t}$, so we can apply Lemma~\ref{le:subset} to
$G^{(t)}$. Combined with assertion ${\rm 5)}$ of Lemma~\ref{le:gt}
we deduce: 
\begin{eqnarray*}
|B|&\le&  
\mu(G^{(t)})|V|+\frac{2|V|}{d(d-1)^{t-1}\left(1- 
\frac{\alpha}{\sqrt{4\alpha^2}}\right)}\\
&\le&
\left[8(d-1)^{-t/4}+4(d-1)^{-t}\right]|V|\\
&\le& 12 (d-1)^{-t/4}|V|\\
&\le& 12
(d-1)^{-\frac{1}{4}\left(\frac{g}{8\alpha^2}-2\right)}|V|\le
12(d-1)^{- \frac{g}{64\alpha^2}}|V|.
\end{eqnarray*} 
 \vglue-34pt \phantom{over}\hfill\qed
\vglue8pt 

\Subsec{The discrete cube}\label{section:cube}
The solution for the metric Ramsey problem for the discrete cube
is also based on the notion of Markov type. The discrete cube has
a small girth, and so other ideas are called for. Our analysis
utilizes another family of orthogonal polynomials --- the Krawtchouk
polynomials which appear in many studies related to the discrete
cube.

Let $k\le d$. The degree-$k$ Krawtchouk polynomial for the
$d$-dimensional cube is:
$$
K^{(d)}_k(x)=\sum_{j=0}^k (-1)^j{x \choose j}{d-x\choose k-j}.
$$

Again we need an estimate for the smallest value that this
polynomial takes.

\begin{lemma}\label{le:kra} Let $1\le k\le \frac{d}{2}$ be even. Then\/{\rm :}\/
$$
K_k^{(d)}(x)\ge -\left(\frac{64k}{d}\right)^{k/2}{d\choose k}.
$$
\end{lemma}

\Proof  It is known (see for example \cite{lev}) that all $k$  
zeros of $K_k^{(d)}$ are
real and belong to the interval:
$$
\left[\frac{d}{2}-\sqrt{(k-1)(d-k+2)},\frac{d}{2}+\sqrt{(k-1)(d- 
k+2)}\right].
$$

Since $k$ is even, $K_k^{(d)}$ is symmetric around $d/2$ (i.e.\ $K_k^{(d)}(x)=
K_k^{(d)}$\break\vglue-9pt\noindent $\cdot (d-x)$). It
is also easily checked that the leading coefficient of $K_k^{(d)}$ is
$\frac{(-2)^k}{k!}$.\break\vglue-10pt\noindent  So
$K_k^{(d)}(x)\le 0$ only for $x$ in the above interval. Let $z_1,
z_2, \dots  , z_k$ be the zeros\break\vglue-10pt\noindent   of $K_k^{(d)}(x)$,
$$
K_k^{(d)}(x)=\frac{2^{k}}{k!}\prod_{i=1}^k(z_i-x),
$$
and since $x,z_1,\dots ,z_k$ are all in an interval of length
$2\sqrt{(k-1)(d-k+2)}\le 4\sqrt{k(d-k)}$, it follows that for $x$
in the interval above:
$$
|K_k^{(d)}(x)| \le \frac{2^{k}}{k!}(4\sqrt{k(d-k)})^k \le
\left(\frac{64k}{d}\right)^{k/2}\cdot {d \choose k}.
$$
To verify this inequality, note that after clearing equal terms,
it reduces to $[d(d-k)]^{k/2} \le d(d-1)\dots (d-k+1)$. This
follows by multiplying the inequality $d(d-k) \le (d-j)(d+j-k+1)$
over $j=0,\dots ,k/2$.
\Endproof\vskip4pt 

Let $\Omega_d=\{0,1\}^d$ be the graph of the $d$-dimensional cube.
(Two vectors are adjacent if and only if they differ in exactly
one coordinate.) As before, we consider the $t$-distance graph on
the cube $\Omega_d^{(t)}$. It is well known (e.g.~\cite{dels}) and
easy to show \footnote{ To show this, recall the $2^d$ Walsh
functions $\{W_S~|~S \in \Omega_d\}$ that are defined via
$W_S(T)=(-1)^{\langle S,T\rangle }$. It is not hard to see that they form a
complete set of eigenfunctions for $\Omega_d^{(t)}$ and the
eigenvalue corresponding to $W_S$ is $K_t^{(d)}(n-2|S|)$. } that
the eigenvalues of the graph $\Omega_d^{(t)}$ are the numbers
$K_t^{(d)}(i)$ for $i=0,\dots ,d$ where the $i$-th eigenvalue
appears with multiplicity ${d \choose i}$. This graph is ${d
\choose t}$-regular and so its largest eigenvalue is ${d \choose
t}$. Lemmas~\ref{le:selfmix} and~\ref{le:kra} now yield an
estimate for the self-mixing parameter $\mu(\Omega_d^{(t)})$.

\begin{lemma}\label{le:tdist} For
every even integer $1\le t< d/2${\rm ,}
$$
\mu(\Omega_d^{(t)})\le \left(\frac{64 t}{d}\right)^{t/2}.
$$
\end{lemma}
\vskip8pt

To prove the main result of this section, we need an additional
estimate.

\begin{lemma}\label{le:a} Let $t,x,d$ be integers such that
$2t \le x\le d/2$. Then\/{\rm :}\/
$$
\sum_{j=t/3}^t {x \choose j}{d-x \choose t-j}\le 2
\left(\frac{150x}{d}\right)^{t/3}{d\choose t}.
$$
\end{lemma}
\vskip8pt

\Proof 
Clearly, ${d\choose t} = \sum_j {x \choose j}{d-x \choose t-j}$,
so it suffices to consider the range $x \le \frac{d}{150}$. In
other words, we assume $\frac{t}{3} \le j \le t \le \frac{x}{2}
\le \frac{d}{300}$. In this range, the\break\vglue-11pt\noindent terms decrease
geometrically, ${x \choose j+1}{d-x \choose t-(j+1)} \le
\frac{1}{10}{x \choose j}{d-x \choose t-j}$. It therefore suffices\break\vglue-9pt\noindent 
to show that ${x \choose t/3}{d-x \choose 2t/3} \le
\left(\frac{150x}{d}\right)^{t/3}{d\choose t}.$ Recall the
following elementary and well-known estimates of binomial
coefficients: For every $1\le k \le n$,
$$
\left(\frac{n}{k}\right)^k\le {n \choose k} \le
\left(\frac{en}{k}\right)^k.
$$
We plug this into the expression and simplify, to conclude that
the inequality holds.
\hfill\qed

\begin{theorem}\label{th:cubebounds} There are absolute constants  
$C,c,c'>0$ such that
for every integer $d$ and for every $\alpha\ge 1${\rm ,}
$$ 2^{\left(1-\frac{\log (c'\alpha)}{\alpha^2}\right)d}
\le R_2(\Omega_d;\alpha) \le
C2^{\left(1-\frac{c}{\alpha^2}\right)d}.
$$
\end{theorem}
\vskip8pt

\Proof 
We start with the lower bound. An easy fact from coding theory,
called the Gilbert-Varshamov bound, \cite{macsloane}, states that
there exists a subset $B \subseteq \Omega_d$ such that all $u
\neq v\in B$, are at distance $\ge \frac{d}{\alpha^2}$, and:
$$
|B|\ge \frac{2^d}{\sum_{m\le d/\alpha^2} {d\choose m}} \ge
2^{\left(1-\frac{\log (c'\alpha)}{\alpha^2}\right)d},
$$
where the last inequality follows from standard estimates on
binomial coefficients. Note that for every $u,v\in \Omega_d$,
$\|u-v\|_2=\sqrt{\rho(u,v)}$, where $\rho$ stands for the Hamming
distance. But for every distinct $u,v\in B$, $\frac{d}{\alpha^2}
\le \rho(u,v) \le d$, so that  $\frac{\sqrt{d}}{\alpha} \le
\frac{\rho(u,v)}{\|u-v\|_2} \le \sqrt{d}$. Consequently,
$c_2(B)\le \alpha$.

To motivate the proof of the upper bound, let us sketch a proof
based on Markov type for (a weakening of) the classical
fact~\cite{enflocube} that $c_2(\Omega_d) \ge a \sqrt{d}$ for some
absolute $a > 0$ (in fact, $c_2(\Omega_d)= \sqrt{d}$). The random
walk on $\Omega_d$ almost surely drifts with constant speed from
its point of origin for time $> a' d$ for some absolute $a' > 0$.
This is true because a constant fraction of the coordinates stay
unchanged for this duration. On the other hand, the fact that
Hilbert space has Markov-type 2 implies that the corresponding
walk on an image of $\Omega_d$ will typically drift only
$O(\sqrt{d})$ away from its origin. This discrepancy implies a
metric distortion $\ge c \sqrt{d}$, as claimed. The spirit of the
proof we present is similar. We only have a subset $B \subseteq
\Omega_d$, so we consider (a dense connected component of) the
graph $\Omega_d^{(t)}$. The main technical effort is in estimating
the typical rate of drift from the walk's origin.

We wish to show, then, that if $B \subseteq \Omega_d$ satisfies
$c_2(B) \le \alpha$, then $|B| \le C
2^{\left(1-\frac{c}{\alpha^2}\right)d}$. Let $n=2^d$ and $k=|B|$.
We seek an upper bound on $k$. As in the proof of
Theorem~\ref{th:girthramsey}, we investigate the random walk on
the distance $t$ graph of the graph in question, namely
$\Omega_d^{(t)}$. We define $t$ as the even integer nearest to
$\frac{d}{K\alpha^2}$, where $K$ is a suitably large absolute
constant to be specified later. It can be verified that
$2^{\left(1-\frac{c}{\alpha^2}\right)d} \ge
2n\left(\frac{32t}{d}\right)^{t/2}$ for this choice of $t$.
Therefore, we may assume that:
$$
\frac{k}{n}\ge 2\left(\frac{64t}{d}\right)^{t/2},
$$
or else the required upper bound on $k$ already holds.

Denote by $E_t(B)$ the number of unordered pairs of points of
distance $t$ in~$B$. In terms of the graph $\Omega_d^{(t)}$ this
is: $E_t(B)=|E_{\Omega_d^{(t)}}(B)|$. By Lemma~\ref{le:selfmix}:
$$
2E_t(B)\ge \frac{{d\choose t}k^2}{n}-\mu(\Omega_d^{(t)}){d\choose
t}k \ge \frac{{d\choose
t}k^2}{n}-\left(\frac{64t}{d}\right)^{t/2}{d\choose t}k,
$$
so that:
$$
\delta=\delta_{\Omega_d^{(t)}}(B)=\frac{2E_t(B)}{k} \ge {d\choose
t}\left[\frac{k}{n}-\left(\frac{64t}{d}\right)^ {t/2}\right]\ge
{d\choose t}\cdot\frac{k}{2n}.
$$
There is a connected component $C$ of the subgraph of
$\Omega_d^{(t)}$ spanned by $B$ that has average degree $\delta'
\ge \delta$, i.e.,
$$
\delta'=\delta_{\Omega_d^{(t)}}(C)\ge \delta\ge {d\choose
t}\cdot\frac{k}{2n}.
$$
Let $\{Z_r\}_{r=0}^{\infty}$ be the random walk on
$\Omega_d^{(t)}$ restricted to $C$. We start the walk at the
stationary distribution, viz.,
$$
P(Z_0=v)=\frac{ \mathrm{deg}_C^t(v)}{\delta' |C|},
$$
where $\mathrm{deg}_C^t(v)$ is the degree of vertex $v$ in
$\Omega_d^{(t)}$ restricted to $C$ (i.e.\ the number of elements of
$C$ with Hamming distance $t$ to $v$).

Suppose that our random walk starts from $S \in C$ and reaches,
after some time, a vertex $T$ with $x=\rho(S,T)$. Say that we next
step from $T$ to $W$. We seek an \pagebreak upper bound on the probability
that $\rho(S,W) \le x + \frac{t}{3}$. The total number of
neighbors $W$ of $T$ in $\Omega_d^{(t)}$ for which this holds is
$$
A(x)=\sum_{j \ge \lceil t/3\rceil} {x \choose j}{d-x \choose t-j}.
$$
By Lemma~\ref{le:a}, $A(x)\le 2
\left(\frac{150x}{d}\right)^{t/3}{d \choose t}$ when $2t\le x\le
\frac{d}{2}$.

Now for every possible walk, $\rho(Z_r,Z_0)\le rt$ holds for every
integer $r>1$. For times $2\le r\le \frac{d}{2t}$ we are able to
show that the walk tends to drift at least $a'' t$ per step away
from its origin for some absolute $a'' > 0$.
\begin{eqnarray*}
\E[\rho(Z_{r+1}&,& Z_0)]\\[3pt]
&\ge&\E\Bigg[\frac{
\mathrm{deg}_C^t(Z_r)-A(\rho(Z_r,Z_0))}{\mathrm{deg}_C^t(Z_r)}
\left(\rho(Z_r,Z_0)+\frac{t}{3}\right) \\[3pt]
&&\phantom{\E\Bigg[} +\frac{A(\rho(Z_r,Z_0))}{\mathrm{d 
eg}_C^t(Z_r)}\left(
\rho(Z_r,Z_0)-t\right)\Bigg]\\[3pt]
&=&\E[\rho(Z_r,Z_0)]+\frac{t}{3}-\frac{4t}{3}
\E\left[\frac{A(\rho(Z_r,Z_0))}{\mathrm{deg}_C^t(Z_r)}\right]\\[3pt]
&\ge&
\E[\rho(Z_r,Z_0)]+\frac{t}{3}- 
\frac{8t}{3}\left(\frac{150rt}{d}\right)^{t/3}{d
\choose t}
\E\left[\frac{1}{\mathrm{deg}_C^t(Z_r)}\right]\\[3pt]
&=&\E[\rho(Z_r,Z_0)]+\frac{t}{3}- 
\frac{8t}{3}\left(\frac{150rt}{d}\right)^{t/3}{d\choose
t}
\sum_{v\in  
C}\frac{1}{\mathrm{deg}_C^t(v)}\cdot\frac{\mathrm{deg}_C^t(v)}{\delta'|C 
|}\\[3pt]
&=&
\E[\rho(Z_r,Z_0)]+\frac{t}{3}- 
\frac{8t}{3\delta'}\left(\frac{150rt}{d}\right)^{t/3}{d\choose t}\\[3pt]
&\ge&
\E[\rho(Z_r,Z_0)]+\frac{t}{3}- 
\frac{16t}{3}\left(\frac{150rt}{d}\right)^{t/3}\frac{n}{k}.
\end{eqnarray*} 
As in the proof of Theorem~\ref{th:girthramsey}, we now contrast
this estimate with the fact that Hilbert space has Markov type
$2$. Namely, that for every $r$,
$$
\alpha^2 r t^2 \ge c_2(B)^2 r t^2\ge c_2(C)^2 r
\E[\rho^2(Z_1,Z_0)]\ge \E[\rho^2(Z_r,Z_0)]\ge [\E\rho(Z_r,Z_0)]^2.
$$
Consequently,
\begin{eqnarray*}
\alpha t\sqrt{r} \ge \E[\rho(Z_r,Z_0)] & \ge &
\frac{rt}{3}-\frac{16nt}{3k}\left(\frac{150t}{d}\right)^{t/3}
\sum_{r \ge j \ge 1}j^{t/3}\\
&\ge&
\frac{rt}{3}-\frac{16n}{k}\left(\frac{150t}{d}\right)^{t/3}
r^{\frac{t}{3}+1}.
\end{eqnarray*} 
We set $r=\lceil 36\alpha^2\rceil$ and, as stated above, choose
$t$ as the even integer nearest to $\frac{d}{5500 \alpha^2}$ to
conclude the proof of the desired result.
\hfill\qed

\begin{remark} It is known that for every $1\le p\le 2$ the metric space
$(\ell_p, \|x-y\|_p^{p/2})$ embeds isometrically into $\ell_2$
(see \cite{ww}). It follows that $\ell_p$ has  Markov type $p$
with constant $1$. We can therefore apply the above arguments and
conclude that Theorem~\ref{th:girthramsey} and
Theorem~\ref{th:cubebounds} remain true when dealing with
embeddings into $\ell_p$, $1<p\le 2$. The only necessary
modification is that in the upper bound on the Ramsey function
$\alpha^2$ should be replaced by $\alpha^{p/(p-1)}$.
\end{remark}

{\it Acknowledgments.} The authors would like to express their
gratitude to Guy Kindler, Robi Krauthgamer, Avner Magen and Yuri
Rabinovich for some helpful discussions.

\bigskip
\bigskip

\noindent Yair Bartal, Institute of Computer Science, Hebrew
University, Jerusalem 91904, Israel. \\{\bf yair@cs.huji.ac.il}

\medskip
\noindent Nathan Linial, Institute of Computer Science, Hebrew
University, Jerusalem 91904, Israel. \\ {\bf nati@cs.huji.ac.il}

\medskip
\noindent Manor Mendel, Institute of Computer Science, Hebrew
University, Jerusalem 91904, Israel. \\ {\bf
mendelma@cs.huji.ac.il}

\medskip
\noindent Assaf Naor, Theory Group, Microsoft Research, One
Microsoft Way 113/2131, Redmond WA 98052-6399, USA. \\ {\bf
anaor@microsoft.com}

\bigskip
\bigskip
\noindent 2000 AMS Mathematics Subject Classification: 52C45,
05C55, 54E40, 05C12, 54E40.

\end{document}